\documentclass[11pt]{article}

\usepackage{fullpage}
\usepackage{amsfonts}
\usepackage{amsthm}
\usepackage{enumitem}
\usepackage{todonotes}
\usepackage{xcolor}
\usepackage{multirow}
\usepackage{fix-cm}
\usepackage{placeins}

\usepackage{amsmath,amssymb,mathtools}  
\usepackage{graphicx}                   
\usepackage{booktabs}                   
\usepackage{caption,subcaption}         
\usepackage{siunitx}                    
\usepackage{adjustbox}                  
\usepackage[T1]{fontenc}
\usepackage[utf8]{inputenc}
\usepackage{parskip}                    
\usepackage{comment}
\usepackage{cases}
\usepackage{xfrac}
\usepackage[hidelinks]{hyperref}
\hypersetup{
    colorlinks=true,
    linkcolor={red!50!black},
    citecolor={blue!50!black},
    urlcolor={blue!80!black}
}
\newcommand\rev[1]{#1}

\newcommand{\bU}{\mathbf{U}}
\newcommand{\bX}{\mathbf{X}}
\newcommand{\bV}{\mathbf{V}}

\newcommand{\bu}{\mathbf{u}}
\newcommand{\bg}{\mathbf{g}}
\newcommand{\blf}{\mathbf{f}}
\newcommand{\bs}{\mathbf{s}}
\newcommand{\R}{\mathbb{R}}
\newcommand{\fra}{{{1}/{2}}}

\newcommand{\hL}{h_{\!L}}
\newcommand{\hR}{h_{\!R}}
\newcommand{\hLmin}{\hleft^{\min}}
\newcommand{\hLmax}{\hleft^{\max}}
\newcommand{\hRmin}{\hright^{\min}}
\newcommand{\hRmax}{\hright^{\max}}
\newcommand{\bmu}{\boldsymbol{\mu}}  
\newcommand{\NL}{N_L}
\newcommand{\NR}{N_R}
\newcommand{\Nx}{N_x}
\newcommand{\dx}{\Delta x}
\newcommand{\dt}{\Delta t}

\DeclareMathOperator{\spann}{span}

\newcommand{\tensor}[1]{\mathcal{#1}} 
\newcommand{\Sig}{\boldsymbol{\Sigma}}
\newcommand{\ChiL}{\chi_{p,L}}        
\newcommand{\ChiR}{\chi_{p,R}}        
\newcommand{\Ccore}{\mathbf{C}}       
\newcommand{\balpha}{\boldsymbol{\alpha}} 

\newcommand{\LtwoLtwo}{L^2 (L^2)}
\newcommand{\LtwoHone}{L^2(H^1)}
\newcommand{\Esup}{E^{\sup}}
\newcommand{\Eint}{E^{\rm avr}}
   
\newcommand{\epsloc}{\epsilon_{\text{loc}}}
\newcommand{\rh}{\ell_h}
\newcommand{\rqq}{\ell_{q}}

\newcommand{\hleft}{\hL}
\newcommand{\hright}{\hR}
\newcommand{\cleft}{c}
\newcommand{\cm}{c_m}
\newcommand{\xdam}{x_{\text{dam}}}

\newcommand{\OmegaD}{\Omega}

\title{A parametric tensor ROM for the shallow water dam break problem}

\author{
Md Rezwan Bin Mizan \thanks{Dept. of Mathematics, University of Houston, Houston, TX 77204,
mbinmiza@cougarnet.uh.edu} \quad
Maxim Olshanskii \thanks{Dept. of Mathematics, University of Houston, Houston, TX 77204, maolshanskiy@uh.edu} \quad
Ilya Timofeyev\thanks{Dept. of Mathematics, University of Houston, Houston, TX 77204, itimofey@cougarnet.uh.edu}
}
\date{}

\begin{document}
\maketitle

\begin{abstract}
We develop a variant of a tensor reduced-order model (tROM) for the parameterized shallow-water dam-break problem. This hyperbolic system presents multiple challenges for model reduction, including a slow decay of the Kolmogorov $N$-width of the solution manifold, shock formation, and the loss of smooth solution dependence on parameters. These issues limit the performance of traditional Proper Orthogonal Decomposition based ROMs. Our tROM approach, based on a low-rank tensor decomposition, builds a parameter-to-solution map from high-fidelity snapshots and constructs localized reduced bases via a local POD procedure. We apply this method to \rev{1D dry-bed and wet-bed problems and 2D wet-bed problem with topography and bottom friction}, showing that the non-interpolatory variant of the tROM, combined with Chebyshev sampling near critical parameter values, 
effectively captures parameter-dependent behavior and significantly outperforms standard POD-ROMs. This is especially evident in the wet-bed case, where POD-ROMs exhibit poor resolution of shock waves and spurious oscillations. 

\end{abstract}

Keywords: tensor reduced-order model; parametric dependence; shallow-water equations; dam break problem

\section{Introduction}
\label{sec:intro}
The development of accurate and efficient Reduced-Order Models (ROMs) that capture parametric variations in numerical simulations of time-dependent partial differential equations (PDEs) remains a major challenge in scientific computing. This task is typically problem-dependent and often necessitates adaptations of established techniques. The challenge becomes particularly acute when the parameter vector is high-dimensional. Among the most prominent approaches are traditional Proper Orthogonal Decomposition (POD)-based ROMs, which have seen many successful applications~\cite{benner2015survey}.

However, the efficiency of many traditional ROMs is theoretically constrained by the decay rate of the Kolmogorov $N$-width of the solution manifold for the parametric problem of interest. The Kolmogorov $N$-width quantifies how well a manifold can be approximated by an $N$-dimensional linear subspace. While sharp estimates of the Kolmogorov $N$-width are rare~\cite{melenk2000n,lassila2013generalized}, its equivalence to the best $N$-term approximation in greedy POD~\cite{binev2011convergence} suggests exponential or algebraic decay for parabolic and elliptic problems with sufficiently smooth parametric dependence; see, e.g.,~\cite{cohen2010convergence,buffa2012priori,cohen2016kolmogorov,hesthaven2016certified,bachmayr2017kolmogorov}.

In contrast, hyperbolic problems often exhibit much slower decay rates. For example, for linear problems with discontinuous initial data, the authors in~\cite{ohlberger2015reduced,greif2019decay} observed a decay rate of $O(N^{-1/2})$. This slower decay implies that a relatively high-dimensional subspace is required for accurate solution approximation. This partially explains why POD-based ROMs for hyperbolic problems are comparatively less developed; see, however,~\cite{podromswe2012,abgrall2019,abgrall2018,chan2020}.

The situation becomes even more complex for nonlinear hyperbolic problems due to the formation of shocks, which require accurate tracking. The Shallow Water Equations (SWE), addressed in this study, provide a classical example of such a problem. To overcome the limitations of linear reduction techniques, machine learning (ML) methods have been employed to construct ROMs for hyperbolic systems; see, e.g.,~\cite{abgrall2019,abgrall2018,chan2020}. ML-based reduced models for the SWE include~\cite{timofeyev2025,pinnswe-bihlo2022,pinnswe-qi2024}. However, ML approaches have their own limitations, particularly the lack of a mature theoretical foundation.

A more mathematically grounded alternative is the tensor ROM (tROM) approach~\cite{mamonov2022interpolatory,mamonov2024tensorial,olshanskii2025approximating}, developed recently to address the limitations of POD-based ROMs for parametric PDEs. This method is especially well-suited to time-dependent PDEs with high-dimensional parameter spaces~\cite{mamonov2022interpolatory}. The tROM is a projection-based method that performs high-fidelity simulations across different parameter values to construct an $(N+D)$-dimensional “snapshot” tensor in a low-rank tensor format, where $N$ is the physical dimension and $D$ is the number of parameters. This tensor represents a discrete approximation of the parameter-to-solution map. By operating directly on these low-rank structures, one can compute a parameter-dependent reduced basis for any new parameter value. Further details are given in Section~\ref{sec:tROM}.

In this work, we employ a Low-Rank Tensor Decomposition (LRTD) in  the  Tucker format to build a tROM and distinguish between interpolatory and non-interpolatory variants. For interpolatory tROMs, a local reduced basis is computed via low-rank interpolation of the parameter-to-solution map. In the non-interpolatory case, interpolation is replaced by a local POD procedure. As we shall demonstrate, the non-interpolatory tROM performs better.

We apply the tensor ROM approach to a hyperbolic problem  -- specifically, the one-dimensional SWE dam-break problem with both dry-bed and wet-bed initial conditions \rev{as well as two-dimensional SWE equations with topography,  Manning friction and parametrized initial conditions}. The initial conditions are characterized by the parameters $(\hleft, \hright)$, representing the water heights to the left and right of the dam, respectively. In this case, the parameter space has dimension $D = 2$, \rev{whereas for the 2D problem additional parameters are required to describe the topography and the friction coefficient, resulting in $D=5$}.
The dry-bed case ($\hright = 0$) leads to a relatively simple solution involving a rarefaction wave. In contrast, the wet-bed case ($\hright > 0$) is more challenging due to the formation of a shock wave downstream of the dam. For this case, the shock persists for all times $t > 0$, and the derivative of the solution with respect to $\hright$ near the shock becomes unbounded as $\hright \to 0$. This makes the construction of projection-based ROMs particularly difficult, as the quality of approximation hinges on the regularity of the parameter-to-solution map~\cite{mamonov2025priori}. Thus, even though the physical and parameter dimensions of the problem are low, \emph{the hyperbolic nature of the problem, shock formation and the lack of regularity near $\hright = 0$ \rev{pose} a significant challenge for ROM development}.

We demonstrate that the non-interpolatory tROM, combined with a Chebyshev \rev{distribution of sampling points for the parameter governing the water depth in the bed}, provides an effective strategy. Although the resulting local reduced spaces have higher dimensions than those typically required for parabolic or elliptic problems, this is consistent with the known slow decay of the Kolmogorov $N$-width in hyperbolic settings. We also show that this method significantly outperforms standard POD-ROMs, especially in the wet-bed case, where POD-ROMs struggle to resolve shock waves and tend to generate spurious high-frequency oscillations.

The remainder of the paper is organized as follows. In Section~\ref{sec:method}, we present the governing equations, describe the projection-based ROM approach, and provide an overview of tROM and POD-ROM techniques. For further details on tROMs, we refer the reader to~\cite{mamonov2022interpolatory,mamonov2024tensorial}. Section~\ref{sec:experiments} presents numerical results for both dry-bed and wet-bed cases, including an analysis of local basis generation thresholds, variation in parameter $\hright$, and a performance comparison between interpolatory and non-interpolatory tROMs as well as between POD-ROM and tROM. \rev{Section~\ref{sec:2d_experiments} collects results for the 2D problems parametrized with five parameters.} Finally, in Appendix~\ref{sec:ap}, we discuss the solution regularity of the wet-bed dam-break problem near the dry-bed limit.

\section{Methods and Setup}
\label{sec:method}
In this section, we briefly describe the shallow-water equations and the formalism for developing Reduced-Order Models. In particular, we outline the standard POD-ROM approach and discuss a tensor ROM. We consider two test cases - dam break for the 1D and 2D shallow
water equations.

\subsection{1D Shallow Water Equations}
\label{sec:1DSWE}

We consider the one-dimensional (1D) shallow water equations (SWEs) to model the dam-break problem, a classical benchmark for hyperbolic conservation laws. The dam-break scenario involves a sudden release of water, generating shock and rarefaction waves, making it an ideal test case for numerical methods and reduced-order modeling of hyperbolic systems.

\textbf{Governing Equations.} In conservative form, without source terms, the 1D SWEs are expressed as:
\begin{equation}
  \partial_t\bu + \partial_x\blf = \mathbf{0},
  \label{eq:swe}
\end{equation}
with
\[
  \bu =\begin{bmatrix} h \\ q \end{bmatrix},
  \qquad
  \blf =\begin{bmatrix}
          q \\[2pt]
          q^{2}/h + \fra\,g\,h^{2}
        \end{bmatrix},
\]
where \(h(x,t)\) is the water depth,
\(q(x,t)=h\,u\) the discharge,
\(u(x,t)\) the depth-averaged velocity, and \(g\) gravity.

\textbf{Full-order model.}
Full-order model (FOM) solutions to equation~\eqref{eq:swe} are computed using a standard finite volume discretization outlined below.
\rev{In time, we apply a second-order} accurate explicit Heun (RK2) scheme,
while spatial fluxes use the local \rev{Lax--Friedrichs} (LLF) formula
\begin{equation}
  \blf_{i+\fra} =
  \frac12\bigl[\blf(\bu_i)+\blf(\bu_{i+1})\bigr] -
  \frac{\lambda_{i+\fra}}{2}\bigl(\bu_{i+1}-\bu_i\bigr),
  \label{eq:llf}
\end{equation}
with maximum wave speed 
\(
\lambda_{i+\fra}=
\max\bigl(|u_i|+c_i,\;|u_{i+1}|+c_{i+1}\bigr),\;
c_i=\sqrt{g\,h_i}.
\)
A constant time step \(\dt\) is chosen such that the CFL condition
\(\max_i\lambda_{i+\fra}\,\dt/\dx<1\) holds; simulations advance
to a prescribed final time $T$. \rev{The computational domain is \([0, L_x]\) with $L_x=100$ and is discretized using a uniform grid with \(\Nx\) interior points plus two ghost points (one at each end).} The spatial step size is \(\dx=L_x/\Nx\), with interior points at \(x_i=(i-1/2)\dx\) with \((i=1,\dots,\Nx)\), and ghost points at \(x_0=-\dx/2,\;x_{\Nx+1}=L_x+\dx/2\).

\paragraph{Parametrized initial conditions.}
Placing a depth discontinuity at \(\xdam=L_x/2\) gives
\begin{equation}
  h(x,0)=
  \begin{cases}
    \hL, & x<\xdam,\\[2pt]
    \hR, & x\ge \xdam,
  \end{cases}
  \qquad
  q(x,0)=0,
  \label{eq:db_ic}
\end{equation}
with \(\hL>\hR\ge0\).
These depths form the parameter vector
\(\bmu=(\hL,\hR)\in \mathcal{P}\subset\R^{2}\). Our parameter domain is a box \(\mathcal{P}=[\hLmin,\hLmax]\times[\hRmin,\hRmax]\). The setups with $\hR=0$ and $\hR>0$ are known as the \emph{dry-bed} and \emph{wet-bed} dam break problem, respectively.    

\paragraph{Boundary conditions.}
We impose outflow boundary conditions by copying the nearest interior state to each ghost cell, equivalent to \(\partial_x h=\partial_x q=0\) at \(x=0\) and \(x=L_x\).

\paragraph{Parameters sampling.}
Snapshots are collected on a two-parameter grid.  
The dam depths are sampled uniformly:
\begin{equation}
  \hL^i = \hLmin +
          \frac{(\hLmax-\hLmin)(i-1)}{\NL-1},
  \qquad i=1,\dots,\NL,
  \label{eq:hL_nodes}
\end{equation}
while for the downstream bed depths we use \(\NR\) uniformly distributed  or  Chebyshev nodes:
\begin{equation}
  \hR^j = \fra\!\Bigl[\bigl(\hRmin+2\hRmax\bigr) -
          \bigl(2\hRmax-\hRmin\bigr)
          \cos\!\Bigl(\tfrac{\pi j}{2(\NR-1)}\Bigr)\Bigr],
  \qquad j=1,\dots,\NR-1.
  \label{eq:hR_nodes}
\end{equation}
For every \((\hL,\hR)\) pair we store snapshots of \(h\) and \(q\) at
all \(\Nx\) interior cells at regular time intervals. We use $(\hLmin, \hLmax)=(10, 28)$ and $(\hRmin, \hRmax) = (0, 8)$ in this paper, which is significantly larger compared to other studies of the parametrized dam-break problem known from the literature~\cite{ chang2011numerical, magdalena2022numerical,muchiri2024numerical,peng20121d, podromswe2012,  zoppou2000numerical}.
\rev{We test different values of $N_L$ and $N_R$ in this paper. Particular values of sampled parameters for the Figures presented in this paper are reported in the appendix \ref{sec:ap2}.
}

\subsection{1D SWE Reduced-Order Modeling}

Reduced-order modeling constructs low-dimensional approximations of high-dimensional systems to enable efficient simulations while preserving the dominant dynamics of the full-order model. For the 1D dam-break problem governed by the SWEs, the ROMs facilitate rapid exploration of parametric variations in initial conditions, such as left and right water heights, denoted by the parameter vector \( \boldsymbol{\mu} = (\hL, \hR) \). This section outlines the principles of ROM, introduces  tROM using the LRTD with Tucker format, and details the interpolatory and non-interpolatory tROM approaches employed. Additionally, the POD ROM is discussed for comparison, and implementation details highlight the offline--online strategy that enhances computational efficiency.


A projection ROM approximates the discrete solutions
\(h^n,q^n\in\R^{\Nx}\) at each time level by
\begin{equation}\label{eq:rom_approx}
  h^n \;\approx\; \bV_h^{\rh}(\bmu)\,\balpha_h^n,
  \qquad
  q^n \;\approx\; \bV_q^{\rqq}(\bmu)\,\balpha_q^n,
\end{equation}
where \rev{the columns of}
\(\bV_h^{\rh}(\bmu)\in\R^{\Nx\times\rh},\bV_q^{\rqq}(\bmu)\in\R^{\Nx\times\rqq}\)
are orthonormal bases in $\rh$-dimensional and $\rqq$-dimensional ROM spaces, respectively,  possibly dependent on
\(\bmu\) where \(\rh\ll\Nx\) and \(\rqq\ll\Nx\). The  coefficients \(\balpha_h^n\in\R^{\rh}\) and
\(\balpha_q^n\in\R^{\rqq}\) \rev{are recovered as solutions of the projected problem and so they}
evolve in time \rev{and depend on the $\bmu$ and initial conditions}.
Galerkin projection of the FOM onto these subspaces is given by
\begin{equation}\label{eq:rom_proj}
  \begin{aligned}
    (\bV_h^{\rh})^T
    \Bigl(
      \frac{h^{\,n+1}-h^{\,n}}{\dt}
      + \frac{1}{\dx}
        \bigl(
          \blf_{h,i+\fra}^{\,n}
          -\blf_{h,i-\fra}^{\,n}
        \bigr)
    \Bigr) &= \mathbf{0},\\[4pt]
    (\bV_q^{\rqq})^T
    \Bigl(
      \frac{q^{\,n+1}-q^{\,n}}{\dt}
      + \frac{1}{\dx}
        \bigl(
          \blf_{q,i+\fra}^{\,n}
          -\blf_{q,i-\fra}^{\,n}
        \bigr)
    \Bigr) &= \mathbf{0},
  \end{aligned}
\end{equation}
where \(\blf_{h,i+\fra}^{\,n}\) and
\(\blf_{q,i+\fra}^{\,n}\) are the interface fluxes
evaluated midway through the RK-2 step.
Initial conditions are also needs to be projected,
\(
  \balpha_h^0 = (\bV_h^{\rh})^T h^0,\;
  \balpha_q^0 = (\bV_q^{\rqq})^T q^0.
\)

The core challenge lies in constructing \emph{parameter-specific} \rev{bases} \(\bV_h^{\rh}(\bmu)\) and \(\bV_q^{\rqq}(\bmu)\). The tROM exploits the 
\rev{space--time--parameter} structure of snapshot data, offering superior adaptability compared to the global POD bases.

\subsubsection{1D SWE Tensor ROM}\label{sec:tROM}

Tensor ROM enhances projection-based ROM by organizing snapshot data into a higher-order tensor, preserving the intrinsic structure across spatial, temporal, and parametric dimensions. Unlike classical POD, which flattens snapshots into a matrix, tROM employs LRTD to compress data and generate parameter-adaptive bases, making it well-suited for problems with a parametric variability, such as the dam-break problem with its shock fronts.

In tROM, the goal is to approximate solutions in low-dimensional subspaces \rev{spanned by the columns of} \( \bV_h^{\rh}(\bmu), \bV_q^{\rqq}(\bmu) \subset \mathbb{R}^{\Nx} \), where \( \rh,\rqq \ll \Nx \), tailored to a parameter \( \bmu \in \mathcal{D} \). The snapshot set, comprising FOM solutions \( \{ h(\mathbf{x}, t_n; \bmu_k), q(\mathbf{x}, t_n; \bmu_k) \}_{n,k} \) for time steps \( t_n \) and parameter samples \( \bmu_k \in \mathcal{D} \), is organized into 4-D tensors \( \tensor{Q}_h, \tensor{Q}_q \in \mathbb{R}^{\Nx \times N_{\mu_1} \times N_{\mu_2} \times N_T} \), where \( N_{\mu_1}, N_{\mu_2} \) are the number of samples for \(\hL\) and \(\hR\), and \( N_T \) is the number of temporal snapshots.

\paragraph{LRTD}

We adopt the Tucker decomposition, a multilinear LRTD method, to approximate the snapshot tensors. For the water-depth tensor, the approximation is
\begin{equation}\label{eq:Tucker}
\tensor{Q}_h \approx \widetilde{\tensor{Q}}_h
  = \tensor{G}_h
    \times_1 \mathbf{W}_h
    \times_2 \Sig_h^{(1)}
    \times_3 \Sig_h^{(2)}
    \times_4 \mathbf{V}_h,
\end{equation}
where \(\tensor{G}_h \in \mathbb{R}^{\tilde{N}_x \times \tilde{N}_{\mu_1} \times \tilde{N}_{\mu_2} \times \tilde{N}_T}\) is the core tensor, \(\mathbf{W}_h \in \mathbb{R}^{\Nx \times \tilde{N}_x}\), \(\Sig_h^{(1)} \in \mathbb{R}^{N_{\mu_1} \times \tilde{N}_{\mu_1}}\), \(\Sig_h^{(2)} \in \mathbb{R}^{N_{\mu_2} \times \tilde{N}_{\mu_2}}\), and \(\mathbf{V}_h \in \mathbb{R}^{N_T \times \tilde{N}_T}\) are orthogonal factor matrices; \(\tilde{N}_x,\tilde{N}_{\mu_1},\tilde{N}_{\mu_2},\tilde{N}_T\) are reduced ranks (\(\tilde{N}_x \ll \Nx\)). The operator \(\times_m\) denotes mode-\(m\) tensor--matrix multiplication. The approximation error is controlled by a tolerance \(\epsilon_h\):
\begin{equation}\label{eq:TuckerApprox}
\|\tensor{Q}_h - \widetilde{\tensor{Q}}_h\|_F \le \epsilon_h \,\|\tensor{Q}_h\|_F.
\end{equation}  

\rev{For a matrix $\mathbf{A}$, we use the notation $\mathbf{A}[:, I]$ to denote the submatrix consisting of the columns of $\mathbf{A}$ indexed by $I$.}
The columns of the spatial factor matrix \(\mathbf{W}_h\) spans a universal subspace:
\begin{equation}\label{eq:univSpace}
  V_{\text{univ},h} = \spann\{\,\mathbf{W}_h\rev{[:,1:\tilde{N}_x]}\,\} \subset \mathbb{R}^{\Nx},
\end{equation}
representing information given \rev{by} all collected snapshots up to the LRTD error. 
The parameter factor matrices \(\Sig_h^{(1)},\Sig_h^{(2)}\) encode parametric dependencies, enabling the construction of local bases for unseen parameters. 

\rev{We would like to emphasize that time is treated as one of the dimensions of the full tensor $\tensor{Q}_h$, 
as well as the core tensor $\tensor{G}_h$. This imposes some restrictions on the sampling time-step and the total time in all simulations.
In particular, a uniform time step $\Delta t$ must be utilized in all simulations, and snapshots must be collected at exactly the same sub-sampling times $t_k$, $k=1,\ldots,N_T$. Furthermore, the number of snapshots, $N_T$ must be the same in all simulations. This implies that the total time, $T$, is also the same in all simulations.}

\subsubsection{1D SWE Interpolatory tROM}

The interpolatory tROM \cite{mamonov2022interpolatory, mamonov2024tensorial} generates parameter-specific bases by interpolating the parameter factor matrices. 
\rev{For a parameter \( \bmu^* = (\hL^*, \hR^*) \), not necessarily from the sampling set, the online phase constructs interpolation vectors using Lagrange interpolation of order \( p \):}
\begin{equation}\label{eq:interpVector}
\ChiL(\hL^*)_k = 
\begin{cases} 
\displaystyle
\prod_{\substack{m=1 \\ m \neq k}}^{p}
\frac{\hL^{(i_m)} - \hL^*}{\hL^{(i_m)} - \hL^{(k)}} & \text{if } k = i_k \in \{i_1, \ldots, i_p\}, \\[6pt]
0 & \text{otherwise},
\end{cases}
\end{equation}
and similarly for \( \ChiR(\hR^*) \), selecting the \( p \) closest sampled parameters \( \hL^{(i_m)}, \hR^{(i_m)} \). The local core matrix is formed by contracting the core tensor:
\begin{equation}\label{eq:interpCore}
\Ccore_h(\bmu^*) \;=\; 
\tensor{G}_h \times_2 \ChiL(\hL^*) \times_3 \ChiR(\hR^*) 
\in \R^{\tilde{N}_x \times \tilde{N}_T},
\end{equation}
followed by the SVD:
\( \Ccore_h(\bmu^*) = \widehat{\bu}_h \Sig_h \widehat{\bV}_h^{\!T} \).
The local basis is then given by \rev{the columns of}
\begin{equation}\label{eq:localBasis}
\bV_h^{\rh}(\bmu^*) = \mathbf{W}_h \,\rev{\widehat{\bU}_h[:,\,1:\rh],\quad 
\text{col}(\bV_h^{\rh}(\bmu^*))}\subset V_{\text{univ},h},
\end{equation}
where \( \rh \) is chosen based on a local energy threshold \( \epsloc \). The process is analogous for \( \bV_q^{\rqq}(\bmu^*) \).

\subsubsection{1D SWE Non-Interpolatory tROM}

The non-interpolatory tROM \cite{olshanskii2025approximating} constructs local bases by combining snapshot data from several nearest parameter grid points (four points are taken in our implementation), avoiding interpolation to enhance accuracy for discontinuous solutions. For \( \bmu^* \), the method identifies bracketing indices such that \( \hL^{(k_{\text{low}})} \leq \hL^* < \hL^{(k_{\text{high}})} \) and \( \hR^{(l_{\text{low}})} \leq \hR^* < \hR^{(l_{\text{high}})} \). Reduced matrices are computed:
\begin{equation}\label{eq:noninterpTensors}
\begin{aligned}
\Ccore_{h,1} &= \tensor{G}_h \times_2 \Sig_h^{(1)}[k_{\text{low}}, :] \times_3 \Sig_h^{(2)}[l_{\text{low}}, :], \\
\Ccore_{h,2} &= \tensor{G}_h \times_2 \Sig_h^{(1)}[k_{\text{low}}, :] \times_3 \Sig_h^{(2)}[l_{\text{high}}, :], \\
\Ccore_{h,3} &= \tensor{G}_h \times_2 \Sig_h^{(1)}[k_{\text{high}}, :] \times_3 \Sig_h^{(2)}[l_{\text{low}}, :], \\
\Ccore_{h,4} &= \tensor{G}_h \times_2 \Sig_h^{(1)}[k_{\text{high}}, :] \times_3 \Sig_h^{(2)}[l_{\text{high}}, :],
\end{aligned}
\end{equation}
and form a matrix \( \Ccore_h = [\Ccore_{h,1}, \Ccore_{h,2}, \Ccore_{h,3}, \Ccore_{h,4}] \). The SVD \( \Ccore_h = \widehat{\bu}_h \Sig_h \widehat{\bV}_h^{\!T} \) yields the basis as in \eqref{eq:localBasis}, with rank \( \rh \) set by \( \epsloc \). The process is repeated for \( q \). 

\subsubsection{1D SWE POD-ROM}

The POD-ROM constructs global bases using a matrix SVD, offering a simpler but less adaptive alternative. Snapshots are organized into matrices \( \bX_h, \bX_q \in \mathbb{R}^{N_x \times (N_{\mu_1} N_{\mu_2} N_T)} \), and a truncated SVD is applied:
\begin{equation}\label{eq:PODbasis}
\bX_h \approx \bU_h \Sig_h \bV_h^T, \quad \bX_q \approx \bU_q \Sig_q \bV_q^T,
\end{equation}
where \rev{the reduced bases} \( \bV_h^{\rh} = \bU_h[:, 1:\rh] \), \( \bV_q^{\rqq} = \bU_q[:, 1:\rqq] \) are the leading POD modes, selected to retain at least a 
\( 1-\epsilon_{\text{POD}} \) fraction of the total energy

\subsubsection{Implementation}

The tROM and POD-ROM employ an \rev{offline--online} strategy to balance computational effort.

\paragraph{Offline Stage}
\begin{enumerate}
\item[(O1)] \text{Snapshot Generation}: Solve the FOM over the parameter grid, storing snapshots at selected time steps in tensors \( \tensor{Q}_h, \tensor{Q}_q \) (tROM) or matrices \( \bX_h, \bX_q \) (POD).
\item[(O2)] \text{Compression}: For tROM, apply a truncated Higher-Order SVD (HOSVD) to \( \tensor{Q}_h, \tensor{Q}_q \), yielding core tensors and factor matrices with tolerances \( \epsilon_h, \epsilon_q \). For POD, perform an SVD on \( \bX_h, \bX_q \), obtaining global bases with threshold \( \epsilon_{\text{POD}} \).

\end{enumerate}

\paragraph{Online Stage}
\begin{enumerate}[label=(N\arabic*)]
  \item \text{Basis Construction:} Given a new parameter \(\bmu^* = (\hL^*, \hR^*)\),  
    \begin{itemize}
      \item \emph{Interpolatory tROM:}  
        Compute the local core via \eqref{eq:interpCore}, perform its thin SVD, and assemble 
        \(\bV_h^{\rh}(\bmu^*)\), \(\bV_q^{\rqq}(\bmu^*)\) using \eqref{eq:localBasis} with \(\epsloc\).  
      \item \emph{Non‐Interpolatory tROM:}  
        Form the four reduced cores from \eqref{eq:noninterpTensors}, concatenate them into \(\mathbf{C}_h\) (and \(\mathbf{C}_q\)), compute its thin SVD, and assemble 
        \(\bV_h^{\rh}(\bmu^*)\), \(\bV_q^{\rqq}(\bmu^*)\) via \eqref{eq:localBasis} with \(\epsloc\).  
      \item \emph{POD-ROM:}  
        Reuse the precomputed global bases \( \bV_h^{\rh} = \bU_h[:, 1:\rh] \), \( \bV_q^{\rqq} = \bU_q[:, 1:\rqq] \).  
    \end{itemize}

  \item \rev{\text{Flux Computation:}  Lift the solution in the full space and compute nonlinear fluxes on the mesh $x_i$, $i=1,\ldots,N_x$.}

  \item \text{Reduced Operator Projection:}  
    Project the computed flux operators onto the local bases \rev{given by the columns of} \(\bV_h^{\rh}(\bmu^*)\), \(\bV_q^{\rqq}(\bmu^*)\) (tROM) or  \(\bV_h^{\rh}\), \(\bV_q^{\rqq}\) (POD).

  \item \text{Time Integration:}  
    Integrate the reduced system \eqref{eq:rom_proj} (using Heun’s RK-2) in the chosen low-dimensional basis, and reconstruct 
    \(h^n = \bV_h^{\rh}(\bmu^*)\,\balpha_h^n\), 
    \(q^n = \bV_q^{\rqq}(\bmu^*)\,\balpha_q^n\) (or analogous POD form) when required.
\end{enumerate}

\rev{As an alternative to the step (N2) above, the Discrete Empirical Interpolation Method (DEIM) \cite{Sorensen2010DEIM} 
can be employed to develop a low-dimensional approximation for nonlinear terms. To this end, snapshots of the nonlinear term are collected, and another reduced basis is constructed. We do not use the DEIM method in this work to clearly distinguish between the error due to a fixed (POD-ROM) and adapted (tROM) parameter-dependent basis for the solution itself. The DEIM approach can be integrated with tROM \cite{mamonov2024tensorial}, and we expect that the DEIM approximation of nonlinear terms would not change the overall conclusions of this paper.}

\subsubsection{Computational cost}
\rev{
The computational cost of the offline stage is dominated by solving 
$N_{\mu_1} \times N_{\mu_2}$ full-order model problems for the parameter values located at the grid nodes of the parameter domain. 
This cost is the same for POD and tROM. 
For POD, some reduction of the offline cost can be achieved by using a greedy parameter-selection strategy, as in the Reduced Basis method~\cite{hesthaven2016certified}. 
For tROM, the number of required snapshots may be reduced by computing a low-rank approximation of the snapshot tensor from sparse parameter samples using tensor completion~\cite{mamonov2024model} 
or cross-interpolation techniques~\cite{budzinskiy2025lowrank}.

The snapshot collection is followed by a truncated SVD of an unfolding matrix of size 
$N_x \times (N_{\mu_1} N_{\mu_2} N_T)$ for the POD-ROM. 
In the case of two parameters and an HOSVD (as considered in this paper), truncated SVDs are applied to four unfolding matrices of sizes
\[
N_x \times (N_{\mu_1} N_{\mu_2} N_T),\qquad
N_{\mu_1} \times (N_x N_{\mu_2} N_T),\qquad
N_{\mu_2} \times (N_x N_{\mu_1} N_T),\qquad
N_T \times (N_{\mu_1} N_{\mu_2} N_x).
\]

The projection of the FOM model onto the reduced space has the same computational complexity for both ROMs.
In addition, the tensor ROM requires computing the local basis during the online stage. This is a computationally inexpensive procedure involving matrices and vectors of reduced dimensions, and it entails a cost of $O(\widetilde{N}_x\widetilde{N}_T\widetilde{N}_{\mu_1}\widetilde{N}_{\mu_2})$ operations for contracting the core tensor to obtain the  $\widetilde{N}_x\times \widetilde{N}_T$ core matrix and performing a truncated SVD of this matrix.
} 


\subsection{2D Shallow Water Equations}
\label{sec:2DSWE}
\rev{
The ideas outlined above can be extended to higher-dimensional PDEs and a larger number of parameters.
To illustrate this, we  consider the 2D shallow-water equations with topography and Manning friction on a rectangular domain 
\(\OmegaD=(0,L_x)\times(0,L_y)\) and time interval \(t\in[0,T]\).
The dependent variables are 
\(\bu = ( h, q_x,  q_y )^T\), 
where \(h(x,y,t)>0\) is the water depth and \((q_x, q_y) = (hu,\,hv)\) are discharges in the horizontal and vertical directions, respectively. 
The SWE with bottom topography and Manning friction reads
\begin{align}
\partial_t h + \partial_x (q_x) + \partial_y (q_y) &= 0, \label{eq:cont}\\[0.25em]
\partial_t (q_x) + \partial_x\!\Big(\frac{q_x^2}{h} + \tfrac{1}{2} g h^2\Big)
+ \partial_y \!\Big(\frac{q_x q_y}{h}\Big) 
&= - g h \,\partial_x Z_b - \tau_x, \label{eq:momx}\\[0.25em]
\partial_t (q_y) + \partial_x\!\Big(\frac{q_x q_y}{h}\Big)
+ \partial_y \!\Big(\frac{q_y^2}{h} + \tfrac{1}{2} g h^2\Big) 
&= - g h \,\partial_y Z_b - \tau_y, \label{eq:momy}
\end{align}
and we consider Gaussian hill bottom topography
\[
Z_b(x,y)= H\exp\!\Big(-\tfrac{(x-x_h)^2}{W_x^2}-\tfrac{(y-y_h)^2}{W_y^2}\Big),
\]
and Manning friction
\[
\tau_x = g n^2 \frac{u\sqrt{u^2+v^2}}{h^{4/3}},
\qquad
\tau_y = g n^2 \frac{v\sqrt{u^2+v^2}}{h^{4/3}},
\qquad\text{with}~~
u=\frac{q_x}{h},\qquad 
v=\frac{q_y}{h}.
\]
Here, the center of the hill is $(x_h,y_h)$, the width is determined by $(W_x,W_y)$, and $n$ is the bottom roughness coefficient.

The initial conditions are 
\[
\eta(x,y,0)=
\begin{cases}
\eta_L, & x<x_{\mathrm{dam}},\\[0.15em]
\eta_R, & x\ge x_{\mathrm{dam}}.
\end{cases}\qquad
h(x,y,0)=\eta(x,y,0)-Z_b(x,y),\quad
q_x=q_y=0,
\]
where the position of the dam is given by $x_\mathrm{dam} = L_x/2$.
Here $\eta(x,y,t) = h(x,y,t) + Z_b(x,y)$ is the absolute height of the water free surface.
We impose outflow (zero-gradient) boundary conditions at $x=0$ and $x=L_x$ and reflective boundary conditions at $y=0$ and $y=L_y$. 

We study the performance of the tROM reduced model with respect to variations of the 5D parameter vector
\[
\boldsymbol{\mu}=(\eta_L,\ \eta_R,\ n,\ x_h,\ y_h),
\]
with components belonging to bounded intervals \(\mathcal{D}_{\eta_L},\mathcal{D}_{\eta_R},\mathcal{D}_n,\mathcal{D}_{x_h},\mathcal{D}_{y_h}\) specified later.

For the Full Order Model, 
we discretize \eqref{eq:cont}--\eqref{eq:momy} using the Local Lax Friedrichs finite-volume method on a uniform Cartesian mesh with cell centers \((x_i,y_j)\), spacings \(\Delta x,\Delta y\), and one ghost layer on each side.

Define the physical fluxes 
\[
\blf(\bu)=\begin{bmatrix}
q_x\\[2pt]
\dfrac{q_x^2}{h}+\frac12 g h^2\\[6pt]
\dfrac{q_x q_y}{h}
\end{bmatrix},\qquad
\bg(\bu)=\begin{bmatrix}
q_y\\[2pt]
\dfrac{q_x q_y}{h}\\[6pt]
\dfrac{q_y^2}{h}+\frac12 g h^2
\end{bmatrix}.
\]
The eigenvalues of the Jacobian of  $\blf$ and $\bg$ are \(u-c,\ u,\ u+c\) and \(v-c,\ v,\ v+c\), respectively, with \(c=\sqrt{gh}\).
Hence the maximum wave speeds in \(x\) and \(y\) directions are \(\alpha_x=\max(|u|+c)\) and   is \(\alpha_y=\max(|v|+c)\), respectively.

We utilize the Local Lax--Friedrichs space discretization for fluxes.
At vertical face \((i+\tfrac12,j)\) with left/right states \(\bu_L,\bu_R\), the discrete flux in $x$-direction is defined as 
\[
{\blf}_{i+\frac12,j}
=\tfrac12\!\big(\blf(\bu_L)+\blf(\bu_R)\big)
-\tfrac12\,\alpha_{i+\frac12,j}\,(\bu_R-\bu_L),
\quad
\alpha_{i+\frac12,j}=\max_{L/R} \big(|u|+c\big).
\]
Numerical fluxes \({\bg}_{i,j+\frac12}\) are defined analogously. 

The source terms due to the bottom topography and friction are computed as follows.
Bed slopes are approximated by centered differences,
\[
(\partial_x Z_b)_{i,j}
=\frac{Z_{b\,i+1,j}-Z_{b\,i-1,j}}{2\Delta x},
\qquad
(\partial_y Z_b)_{i,j}
=\frac{Z_{b\,i,j+1}-Z_{b\,i,j-1}}{2\Delta y}.
\]
Velocities are reconstructed as 
\(
u_{i,j}={(q_x)_{i,j}}/{h_{i,j}}, 
\) and 
\(
v_{i,j}={(q_y)_{i,j}}/{h_{i,j}}.
\)
The resulting discrete source vector is
\[
\bs_{i,j} = 
\begin{bmatrix}
0 \\[6pt]
-\,g\,h_{i,j}(\partial_x Z_b)_{i,j}
\;-\;
g\,n^2\,\dfrac{u_{i,j}\sqrt{u_{i,j}^2+v_{i,j}^2}}{h_{i,j}^{4/3}}
\\[10pt]
-\,g\,h_{i,j}(\partial_y Z_b)_{i,j}
\;-\;
g\,n^2\,\dfrac{v_{i,j}\sqrt{u_{i,j}^2+v_{i,j}^2}}{h_{i,j}^{4/3}}
\end{bmatrix}.
\]

We use ghost layers to implement the boundary conditions.
If physical (interior) cells are indexed by \(i=1,\dots,N_x-2\) and \(j=1,\dots,N_y-2\), then
ghost cells are \(i=0,\,N_x-1\) and \(j=0,\,N_y-1\). 
Boundary conditions are imposed directly on the conservative variables. For the west/east boundaries, we use the outflow boundary conditions, and for the south/north boundaries, we use reflective boundary conditions. In discrete variables, the boundary conditions become \\
\emph{West/East (outflow):}
\[
\begin{aligned}
h_{0,j} &= h_{1,j}, 
& (q_x)_{0,j} &= (q_x)_{1,j},
& (q_y)_{0,j} &= (q_y)_{1,j},
\\[4pt]
h_{N_x-1,j} &= h_{N_x-2,j},
& (q_x)_{N_x-1,j} &= (q_x)_{N_x-2,j},
& (q_y)_{N_x-1,j} &= (q_y)_{N_x-2,j}.
\end{aligned}
\]
\emph{South/North (reflective):}
\[
\begin{aligned}
h_{i,0} &= h_{i,1},          
& (q_x)_{i,0} &= (q_x)_{i,1},          
& (q_y)_{i,0} &= -(q_y)_{i,1},
\\[4pt]
h_{i,N_y-1} &= h_{i,N_y-2},  
& (q_x)_{i,N_y-1} &= (q_x)_{i,N_y-2},  
& (q_y)_{i,N_y-1} &= -(q_y)_{i,N_y-2}.
\end{aligned}
\]


We use Heun (RK2) with fixed \(\Delta t\) and a 2D additive CFL check
\[
\mathrm{CFL} \;=\; \Delta t\!\left(
\frac{\max_{i,j}(|u|+c)_{i,j}}{\Delta x}
+
\frac{\max_{i,j}(|v|+c)_{i,j}}{\Delta y}
\right) \;<\;1.
\]
}

\subsection{2D Reduced-Order Modeling Framework}
\rev{
We build intrusive projection ROMs in conservative variables with three approaches: interpolatory tensor ROM (I-tROM), non-interpolatory tROM (NI-tROM), and a global POD-ROM baseline.

The offline phase involves performing FOM simulations and collecting snapshots for each parameter quintuple
$\boldsymbol{\mu}=(\eta_L,\eta_R,n,x_h,y_h)$
and for sampled times \(t_k\), collecting interior snapshots of \(h\), \(q_x\), \(q_y\), and vectorizing to length \(N_s=(N_x-2)(N_y-2)\). 
Therefore, we form \textbf{7-D} tensors (space, $\eta_L$, $\eta_R$, $n$, $x_h$, $y_h$, time):
\[
\mathcal{X}_h,\,\mathcal{X}_{q_x},\,\mathcal{X}_{q_y}
\;\in\; \R^{N_s \times N_{\eta_L}\times N_{\eta_R}\times N_n\times N_{x_h}\times N_{y_h}\times N_t}.
\]

The offline phase also involves performing a 
Tucker decomposition to perform tensor compression. 
Ranks are chosen based on ranks for the I-tROM computed with a given energy thresholding. 
Thus, the full tensor is decomposed as 
\[
\mathcal{X}_h \approx \mathcal{G}_h 
\times_1 W_h 
\times_2 A_h^{(L)} 
\times_3 A_h^{(R)} 
\times_4 A_h^{(n)} 
\times_5 A_h^{(x)} 
\times_6 A_h^{(y)}
\times_7 T_h,
\]
and likewise for \(q_x,q_y\). Here \(W_h\in\R^{N_s\times r_s^{(h)}}\) (spatial factor), \(A_h^{(\cdot)}\) are parameter-mode factors, and \(T_h\in\R^{N_t\times r_t^{(h)}}\).

Construction of local bases for new parameter values
$\mu^*$ (online) proceeds as follows.
For the Interpolatory tROM, we first construct the 
interpolation vectors
\[
\chi_L(\eta_L^*),\ \chi_R(\eta_R^*),\ \chi_n(n^*),\ \chi_x(x_h^*),\ \chi_y(y_h^*)
\] 
from the corresponding parameter factors using quadratic interpolation. Then, a matrix $C_h$ is computed using the core tensor and parameter factors
\[
C_h(\boldsymbol{\mu}^*) = \mathcal{G}_h
\times_2 \chi_L \times_3 \chi_R \times_4 \chi_n \times_5 \chi_x \times_6 \chi_y
\;\in\; \R^{r_s^{(h)}\times r_t^{(h)}},
\]
and SVD is performed \(C_h=U_h\Sigma_h V_h^T\), and the spatial basis is selected as 
\(V_h(\boldsymbol{\mu}^*)=W_h\,U_h(:,1:\ell_h)\), where \(\ell_h\) satisfies energy threshold (locally). The procedure is repeated for \(q_x\) and \(q_y\).

Constructing the non-interpolatory tROM for in-sample parameter values is straightforward since we take exactly the snapshots that correspond to these parameter values.
For an out-of-sample parameter value, 
bracket each parameter by its two neighboring training nodes. With 5 parameters this yields \(\mathbf{2^5=32}\) corner selections. For each corner, contract the core with the \emph{rows} of the parameter factors to obtain \(C_h^{(p)}\in\R^{r_s^{(h)}\times r_t^{(h)}}\), stack
\[
\widetilde{C}_h=\big[\,C_h^{(1)}\ C_h^{(2)}\ \cdots\ C_h^{(32)}\,\big],
\]
perform a thin SVD \(\widetilde{C}_h=\widetilde{U}_h\widetilde{\Sigma}_h \widetilde{V}_h^T\), and set \(V_h=W_h\,\widetilde{U}_h(:,1:\ell_h)\) by a \(\sigma^2\)-energy rule. Repeat for \(q_x,q_y\).

To obtain the reduced basis for POD ROM, we
matricize each tensor over space, parameters, and time, compute an SVD, and keep the leading left singular vectors. For fair comparisons, we match POD ranks with those of the interpolatory tROM at \(\boldsymbol{\mu}^*\).

We use Heun's method to perform the time-stepping in time.
At each Heun substep we
(i) reconstruct interior fields in full space; embed into ghost layers,
(ii) evaluate nonlinear terms in the full space, and
(iii) project nonlinear terms onto reduced bases.
Initialization uses projections of the parametric dam-IC.
We monitor the 2D CFL index
\[
\Delta t\left(\frac{\max(|u|+c)}{\Delta x}+\frac{\max(|v|+c)}{\Delta y}\right)\!,
\]
computed on reconstructed states. (In practice we use a fixed \(\Delta t\) with a guard.)
In particular, the time step is \(\dt=0.05~\mathrm{s}\), and the final time is \(T=7.2~\mathrm{s}\).


The spatial domain is 
\(\OmegaD = [0,L_x]\times[0,L_y]\subset\R^2\) with \(L_x=L_y=200~\mathrm{m}\).
We use a uniform, cell-centered rectangular mesh with \(\Nx=200\) and 
\(N_y=50\) points in $x$ and $y$ direction, respectively. Thus, the mesh size is
$\dx={L_x}/{\Nx}=1.0~\mathrm{m}$ and
$\Delta y={L_y}/{N_y}=4.0~\mathrm{m}.$
Since we consider wet-bed 2D simulations
($\eta_R > 0$), a shock develops in the $x$-direction and, therefore, the $x$-discretization requires a finer resolution compared with the $y$-direction. We also verified that increasing the resolution further does not result in visible differences in simulations of the full order model.

For the 2D case, 
we consider a 5D varying parameter vector \(\bmu=(\eta_L,\eta_R,n,x_h,y_h)\) with ranges
\[
\eta_L\in[8,14]~\mathrm{m},\quad
\eta_R\in[2,6]~\mathrm{m},\quad
n\in[0,0.4],\quad
x_h\in[100,150]~\mathrm{m},\quad
y_h\in[80,120]~\mathrm{m},
\]
and uniform grids for each parameter with
\(
\NL=\NR=N_n=N_{x_h}=N_{y_h}=5.
\)
Thus, the full tensor grid has \(5^5=3125\) parameter cases.
}

\section{1D SWE Numerical Experiments} 
\label{sec:experiments}

\begin{figure}
\centering
 \href{https://youtu.be/BzZbB1vsvpk}{
\begin{subfigure}{0.48\linewidth}
    \centering
    \includegraphics[width=\linewidth]{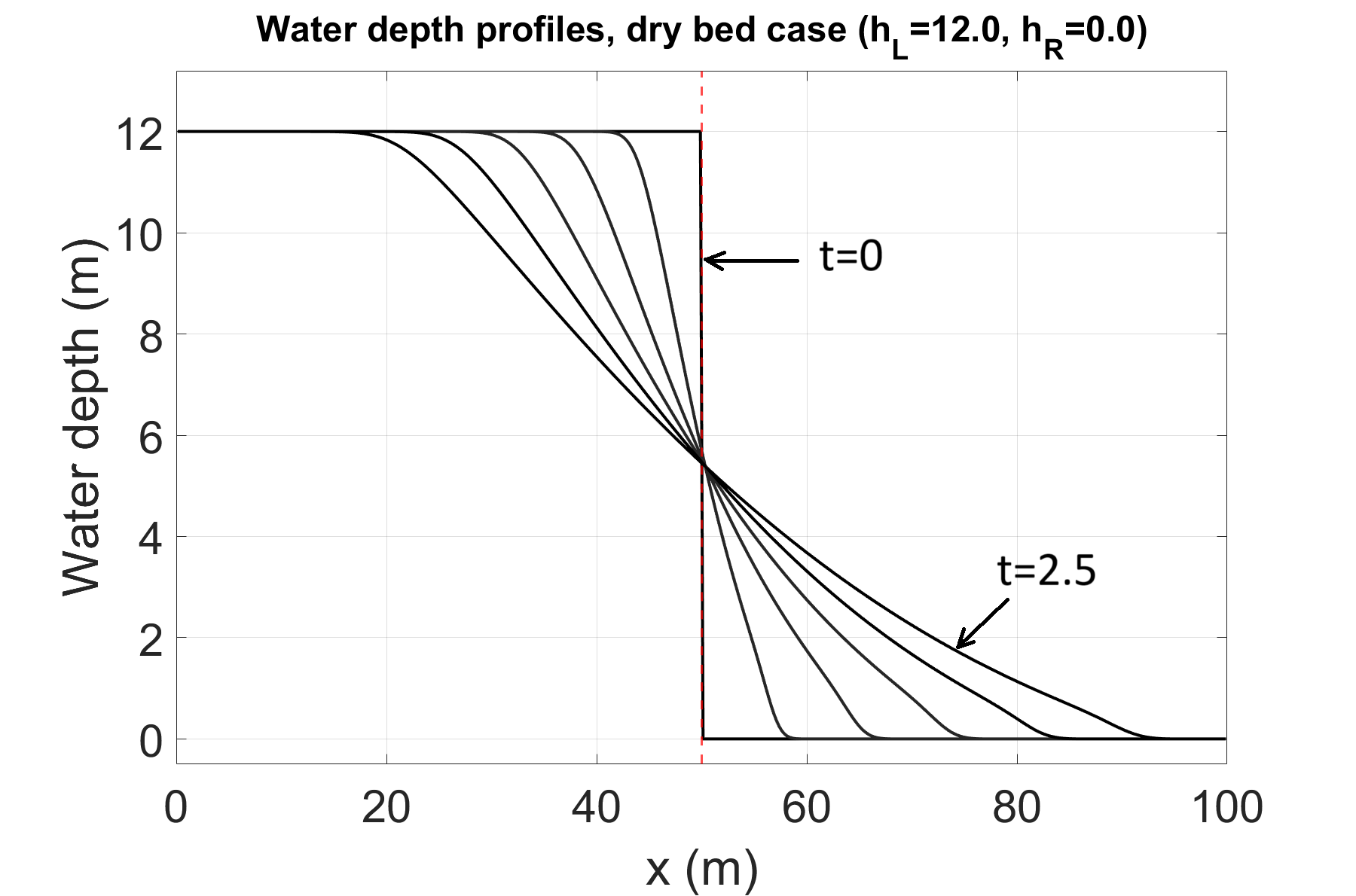}
    \subcaption{}
    \label{fig1:exp1_h_hR_0.14}
\end{subfigure}
\hfill
\begin{subfigure}{0.48\linewidth}
    \centering
    \includegraphics[width=\linewidth]{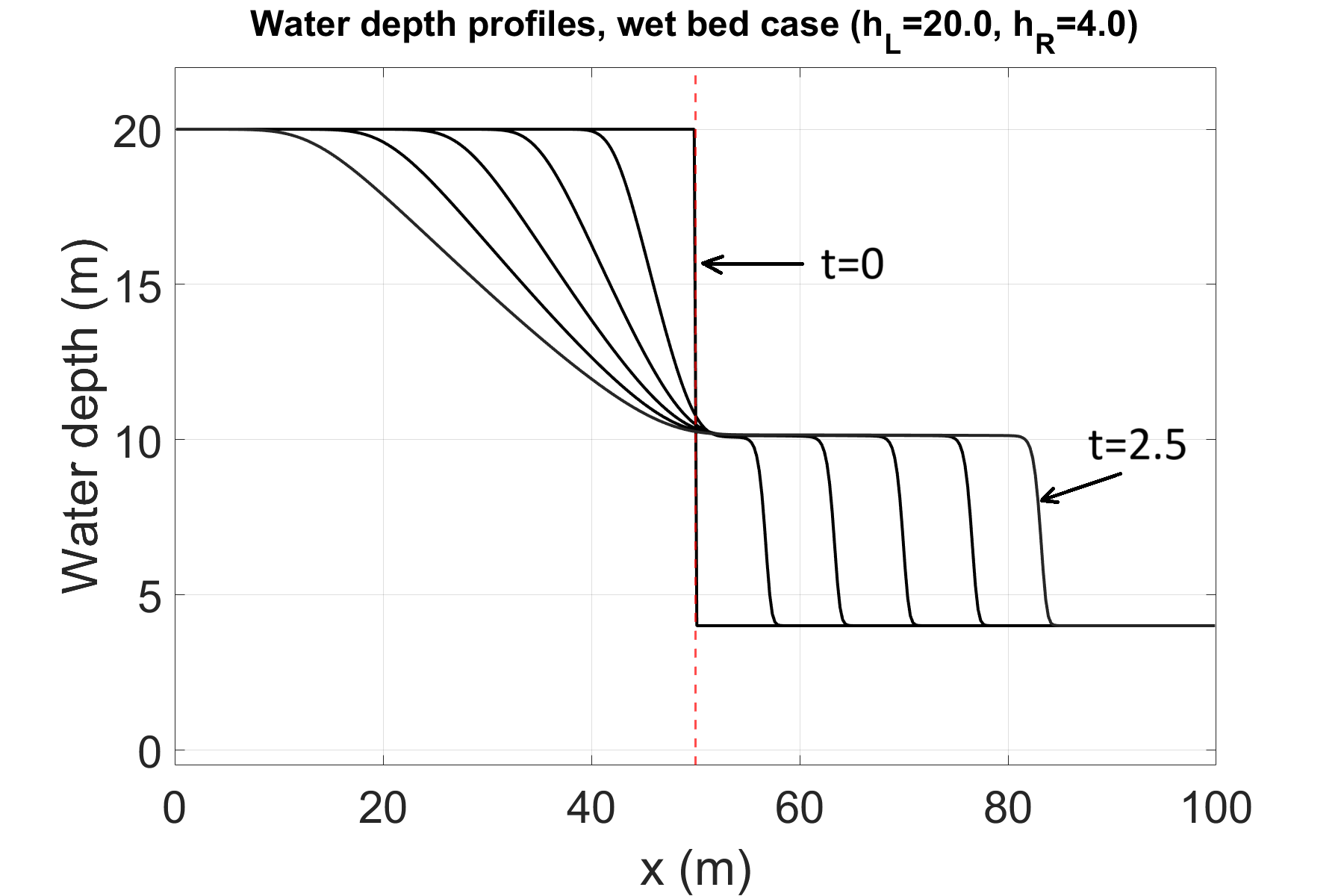}
    \subcaption{}
    \label{fig1:exp1_h_hR_3}
\end{subfigure}
}
\caption{Water depth profiles evolution in time for the cases of dry and wet-bed. Snapshots are shown every $\Delta t = 0.5$. Click the plots for full animation.\label{Fig1} }
\end{figure}

A typical evolution of the water depth for the cases of dry- and wet-beds is illustrated in Figure~\ref{Fig1}, which shows snapshots of the computed FOM solutions for $\hleft = 12$ (dry bed) and $\hleft = 20$, $\hright = 4$ (wet-bed).  
Note that the evolution differs significantly between the dry-bed and wet-bed cases. In the dry-bed case, we observe the formation of a \rev{dry--wet interface} that moves into the initially dry downstream region, with the water depth decreasing smoothly from the dam location to zero at the front edge. A rarefaction wave forms and propagates upstream. This FOM solution approximates the self-similar profiles given by~\eqref{eq:1086}, with the front speed given by $2 \sqrt{g \hleft}$; see Eq.~\eqref{eq:1134}.

In contrast, in the wet-bed case, the initial discontinuity between $\hleft$ and $\hright$ generates a shock wave. The height and speed of the wave depend on the ratio $\hleft/\hright$. Both rarefaction and shock structures evolve smoothly; their speed and position can be determined analytically by solving the system of algebraic equations \eqref{eq:umhm}--\eqref{eq:umhm2}.  

\rev{The computational parameters for the 1D SWE problems are $L_x = 100$, $N_x = 400$, $\Delta x = L_x/N_x = 0.25$, the computational time-step is $\delta t=10^{-4}$, the sub-sampling time-step to collect snapshots is $\Delta t=0.01$, and the number of snapshots is $N_T = 250$.
The spatial resolution is sufficient to resolve shocks in this problem. Comparison with a finer spatial resolution shows very small differences in the shock profile, and this resolution is sufficient to carry out the extensive numerical investigation of the suggested  ROMs.
}


\subsection{Tensor ROM performance}\label{sec:error-metrics}

Let $u_{\mathrm{FOM}}(x,t;\bmu)$ and  
$u_{\mathrm{ROM}}(x,t;\bmu)$ denote the full-order and reduced-order  
solutions of the dam-break problem on the spatial domain  
$\Omega = [0, L_x]$ and the time interval $[0, T]$, where  
$\bmu = (\hL, \hR)$ is a vector of physical parameters.

The error function for a given parameter value is defined as  
\begin{equation*}
  e(x,t;\bmu)
  =
  u_{\mathrm{ROM}}(x,t;\bmu)
  -
  u_{\mathrm{FOM}}(x,t;\bmu).
  \label{eq:error_def}
\end{equation*}

We are interested in the 
\emph{space--time} relative error in the $L^2(0,T;L^2(\Omega))$ norm, given by  
\begin{equation}
  E_{\LtwoLtwo}(\bmu)
  =
  \frac{
    \left(
      \displaystyle\int_{0}^{T}
      \lVert e(\cdot,t,\bmu)\rVert_{L^2(\Omega)}^{2}
      \,\mathrm{d}t
    \right)^{1/2}
  }{
    \left(
      \displaystyle\int_{0}^{T}
      \lVert u_{\mathrm{FOM}}(\cdot,t,\bmu)\rVert_{L^2(\Omega)}^{2}
      \,\mathrm{d}t
    \right)^{1/2}
  }.
  \label{eq:L2L2_error}
\end{equation}
Similarly, the relative error in the $L^2(0,T;H^1(\Omega))$ norm is denoted by  
$E_{\LtwoHone}(\bmu)$.

Let $I_{R} = [\hRmin, \hRmax]$ and $I_{L} = [\hLmin, \hLmax]$.  
For a fixed $\hleft$, the maximum and average errors over the wet-bed water height parameter $\hright$ are defined as  
\begin{equation}
    \Esup_{\LtwoLtwo}(\hleft) =
      \sup_{\hright \in I_{R}} E_{\LtwoLtwo}(\cdot,\hright),\quad
    \Eint_{\LtwoLtwo}(\hleft) =
      \int_{I_{R}} E_{\LtwoLtwo}(\cdot,\hright)\,\mathrm{d}\hright,
  \label{eq:L2L2_sup_int_simple}
\end{equation}
with analogous definitions for  
\rev{$\Esup_{\LtwoHone}(\hright)$ and $\Eint_{\LtwoHone}(\hright)$.}

Finally, the maximum and average errors over the entire parameter domain  
$\Omega_{\mu} = I_L \times I_R$ are given by  
\begin{equation}
  \Esup_{\LtwoLtwo}
  =
  \sup_{\bmu \in \Omega_{\mu}} E_{\LtwoLtwo}(\bmu),
  \qquad
  \Eint_{\LtwoLtwo}
  =
  \int_{\Omega_{\mu}} E_{\LtwoLtwo}(\bmu)\,\mathrm{d}\bmu.
  \label{eq:param_error_sup_int}
\end{equation}
The same definitions apply to the maximum and average errors in the $L^2(0,T;H^1(\Omega))$ norm.  
In practice, these quantities were computed using a Monte Carlo procedure with random sampling of  
$\Omega_{\mu}$.

\subsubsection{Effect of local threshold on  tROM performance}

We start the assessment of tROM by examining the effect of varying the local threshold
$\epsloc$ for parameter-specific basis generation in the online phase
of our tensor reduced-order model.

\begin{figure}
\centering
\begin{subfigure}{0.48\linewidth}
    \centering
    \includegraphics[width=\linewidth]{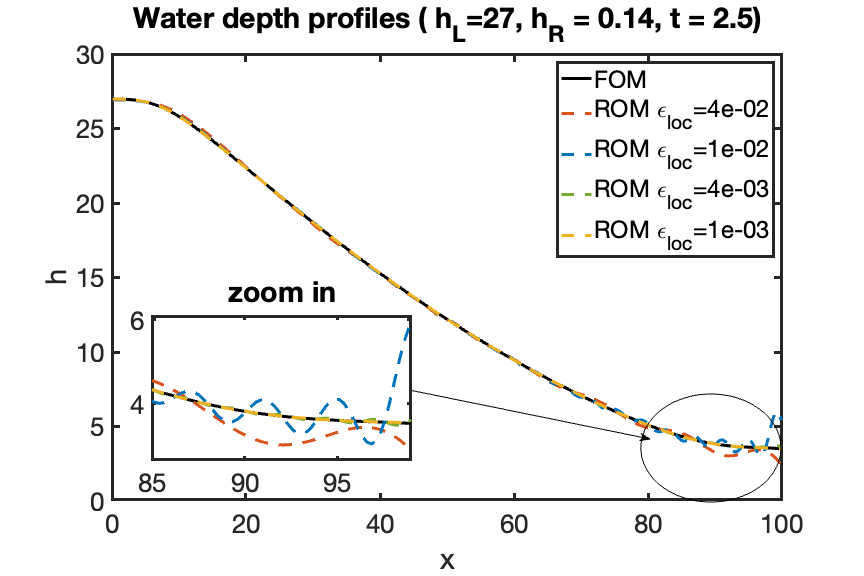}
    \subcaption{}
    \label{fig2:exp1_h_hR_0.14}
\end{subfigure}
\hfill
\begin{subfigure}{0.48\linewidth}
    \centering
    \includegraphics[width=\linewidth]{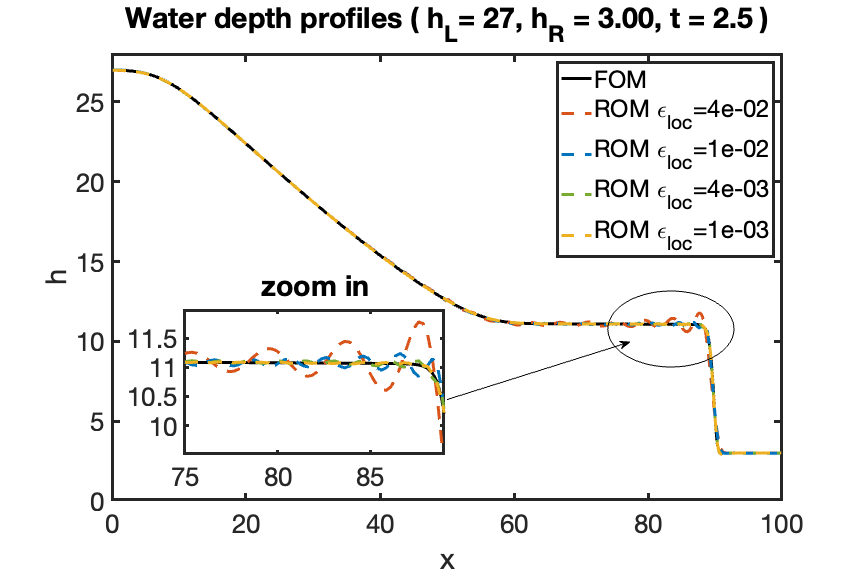}
    \subcaption{}
    \label{fig2:exp1_h_hR_3}
\end{subfigure}

\vspace{0.5em}

\begin{subfigure}{0.48\linewidth}
    \centering
    \includegraphics[width=\linewidth]{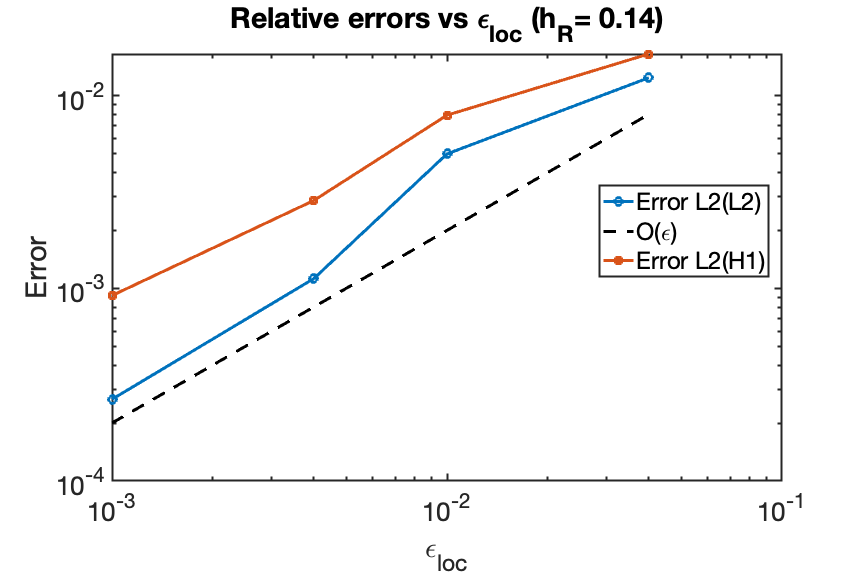}
    \subcaption{}
    \label{fig2:exp1_errors_hR_0.14}
\end{subfigure}
\hfill
\begin{subfigure}{0.48\linewidth}
    \centering
    \includegraphics[width=\linewidth]{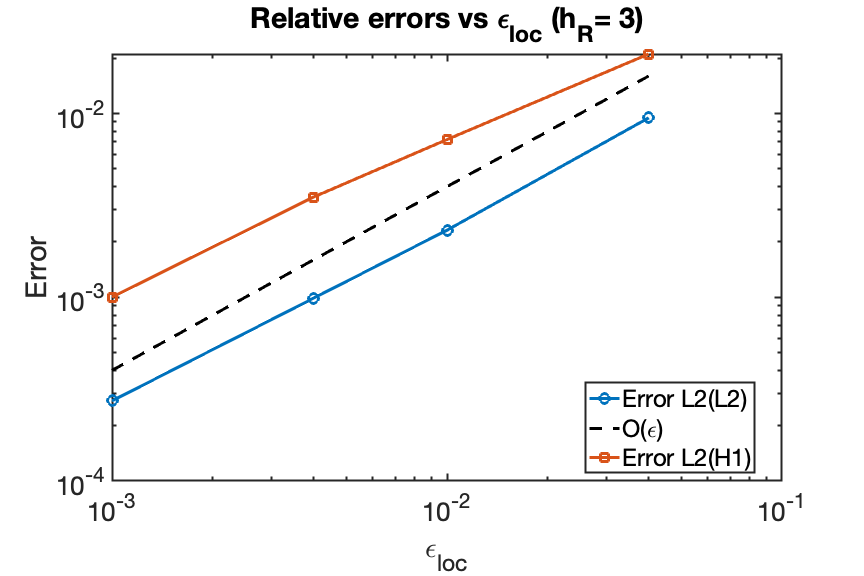}
    \subcaption{}
    \label{fig2:exp1_errors_hR_3}
\end{subfigure}

\vspace{0.5em}

\begin{subfigure}{0.48\linewidth}
    \centering
    \includegraphics[width=\linewidth]{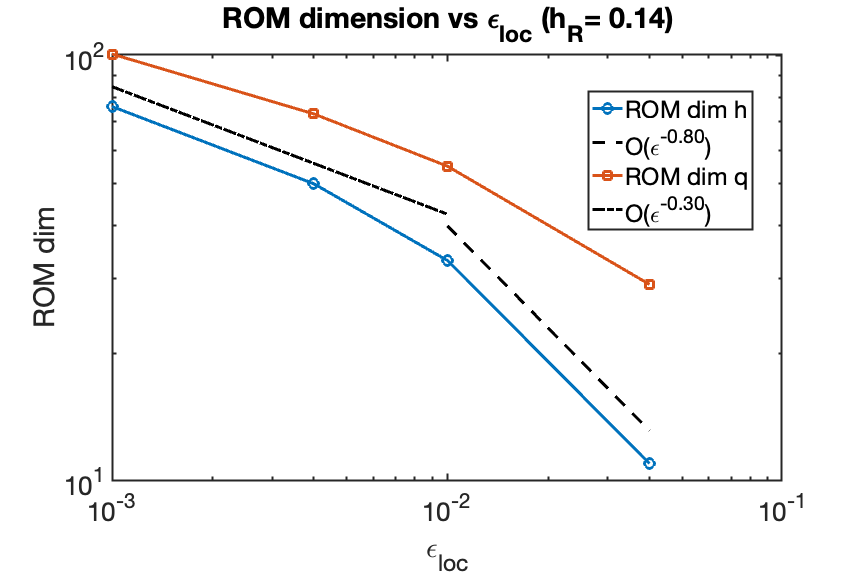}
    \subcaption{}
    \label{fig2:exp1_ranks_hR_0.14}
\end{subfigure}
\hfill
\begin{subfigure}{0.48\linewidth}
    \centering
    \includegraphics[width=\linewidth]{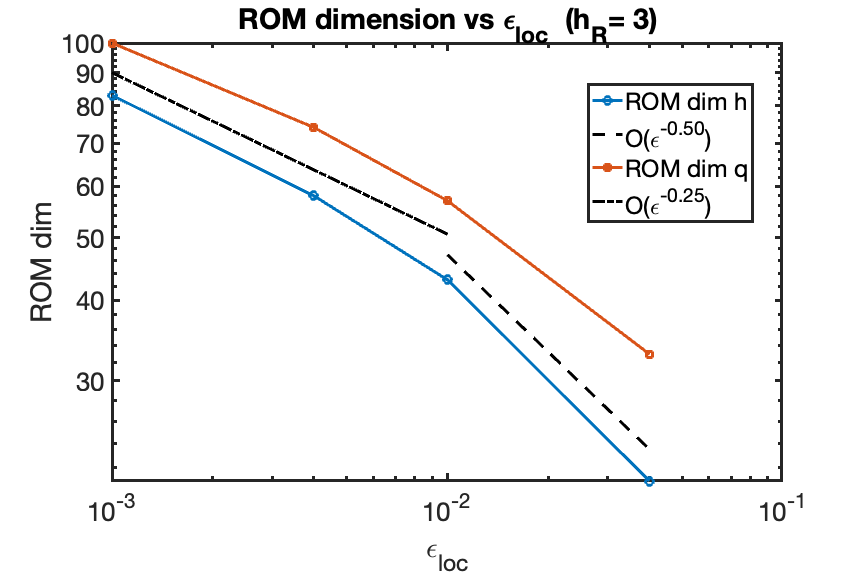}
    \subcaption{}
    \label{fig2:exp1_ranks_hR_3}
\end{subfigure}

\caption{Performance metrics for varying the local threshold \(\epsilon_{\text{loc}}\) (from \(4.0 \times 10^{-2}\) to \(1.0 \times 10^{-3}\)) in the non-interpolatory tROM (\(\hleft = 27 \, \text{m}\)): (a) water depth  profiles at final time for \(\hright = 0.14 \, \text{m}\), with a zoom-in of the step region showing improved shock front resolution as \(\epsilon_{\text{loc}}\) decreases; (b) water depth  profiles at final time for \(\hright = 3.00 \, \text{m}\), with a zoom-in of the step region indicating smoother transitions but similar trends; (c) relative errors in \(h\) versus the local threshold with  \(\hright = 0.14 \, \text{m}\); (d) relative errors in \(h\) versus the local threshold with \(\hright = 3.00 \, \text{m}\), exhibiting a similar error reduction trend; (e) ROM dimensions (\(\rh, \rqq\)) for \(\hright = 0.14 \, \text{m}\), increasing with smaller \(\epsilon_{\text{loc}}\) to retain more basis vectors; (f) ROM dimensions (\(\rh, \rqq\)) for \(\hright = 3.00 \, \text{m}\), also increasing with smaller \(\epsilon_{\text{loc}}\).
\label{fig:exp1_performance_metrics}
} 
\end{figure}

We test two cases with $\NL = 13$, $\NR = 17$,  
$\hL = 27\,\text{m}$, and  
$\hR \in \{0.14, 3.00\}\,\text{m}$,  
varying $\epsloc \in  
\{4.0 \times 10^{-2},\, 1.0 \times 10^{-2},\,  
  4.0 \times 10^{-3},\, 1.0 \times 10^{-3}\}$.  
Water-depth profiles, errors, and ROM dimensions for each $\hR$ are plotted in  
Figure~\ref{fig:exp1_performance_metrics}, with zoom-ins highlighting  
shock-front differences.

The results reveal an expected  trade-off: smaller $\epsloc$ enhances the  
non-interpolatory tROM’s accuracy by retaining more basis vectors,  
thereby reducing errors but increasing ROM dimensions and computational cost,  
with notably improved shock-front resolution for $\hR = \rev{3.00}\,\text{m}$  
due to its sharper discontinuity. From Figure~\ref{fig:exp1_performance_metrics},  
we observe that $\epsloc = 4.0 \times 10^{-3}$ performs reasonably well  
and offers a good balance between accuracy and efficiency (i.e., ROM dimension).  
Thus, we fix $\epsloc = 4.0 \times 10^{-3}$ for all remaining experiments (unless explicitly stated otherwise) to 
ensure consistent performance across varying conditions.

\subsubsection{Impact of wet-bed water depth}
\label{sec:312}

We now study the effect of varying $\hright$ on the accuracy and ROM dimensions of the tensor reduced-order model.
To this end, we fix $\hleft = 25\,\text{m}$, $\NL = 13$, $\NR = 17$, and $\epsloc = 4.0 \times 10^{-3}$ (as it was found to be close to optimal), and evaluate the $E_{\LtwoLtwo}(\hleft)$ and $E_{\LtwoHone}(\hleft)$ errors along with the ROM dimensions $\rh$ and $\rqq$. Next, we fix the ROM dimensions at $\rh = 30$, $\rqq = 50$, instead of using $\epsloc$, and assess performance for varying $\hright$.

Results for both cases (fixed $\epsloc$ and adaptive ROM dimensions and fixed ROM dimensions) are reported in Table~\ref{tab:exp2_combined}. For the adaptive ROM case, we observe that the ROM accuracy varies only slightly across the range of wet-bed water depths, with the resulting ROM dimensions increasing as $\hright$ decreases, reaching a maximum before dropping to their minima at $\hright = 0$.  
In the case of fixed ROM dimensions, the errors remain very stable over the entire range of $\hright$, with the highest accuracy achieved for the dry-bed scenario.
{These results indicate that tROMs can successfully handle the challenging case $\hright\to 0$. In particular, tROMs with both adaptive and fixed dimensions do not result in a sharp error increase as $\hright \to 0$.}

From the results we conclude that the tROM maintains consistent accuracy across varying $\hright$ values. Adaptive ROM dimensions tend to decrease as $\hright$ increases. \rev{One possible explanation is that the derivatives of the solution with respect to $\hright$ are unbounded in the limit $\hright \to 0$ and exhibit greater regularity for larger values of $\hright$; see Appendix~\ref{sec:ap} for details.}. In contrast, fixed ROM dimensions yield higher errors -- particularly for small $\hright$ -- indicating that adaptive basis truncation better captures solution variability. This suggests that adaptive ROM dimension selection is more effective across varying flow regimes and will be our default choice for the rest of the paper.

\renewcommand{\thetable}{2}

\sisetup{
  round-mode        = places,
  round-precision   = 2,
  retain-zero-exponent = true
}

\begin{table}\small
\centering
\begin{adjustbox}{max width=0.95\linewidth}
\begin{tabular}{@{}c
                S[table-format=1.2e-2]  S[table-format=1.2e-2]@{\extracolsep{10pt}}
                S[table-format=1.2e-2]  S[table-format=1.2e-2]@{\extracolsep{10pt}}
                c c c c@{}}
\toprule
\multirow{2}{*}{\(\hR\)}
  & \multicolumn{2}{c}{\(E_{\LtwoLtwo}(\bmu)\)}
  & \multicolumn{2}{c}{\(E_{\LtwoHone}(\bmu)\)}
  & \multicolumn{2}{c}{\(\rh\)}
  & \multicolumn{2}{c}{\(\rqq\)}\\
\cmidrule(lr){2-3} \cmidrule(lr){4-5} \cmidrule(lr){6-7} \cmidrule(l){8-9}
 & {\(\text{a}\)} & {\(\text{f}\)}
 & {\(\text{a}\)} & {\(\text{f}\)}
 & {\(\text{a}\)} & {\(\text{f}\)}
 & {\(\text{a}\)} & {\(\text{f}\)} \\
\midrule
0    & 1.594043e-03 & 1.275694e-04 & 2.517600e-03 & 2.593277e-04 & 15 & 30 & 23 & 50 \\
0.02 & 2.512907e-03 & 3.410967e-03 & 3.742781e-03 & 4.816119e-03 & 34 & 30 & 55 & 50 \\
0.07 & 6.865494e-04 & 7.579952e-03 & 1.713773e-03 & 1.103102e-02 & 51 & 30 & 74 & 50 \\
0.5  & 8.014219e-04 & 4.500751e-03 & 2.838054e-03 & 1.065671e-02 & 64 & 30 & 79 & 50 \\
1    & 8.054924e-04 & 5.082235e-03 & 2.836267e-03 & 1.264904e-02 & 65 & 30 & 80 & 50 \\
2    & 8.576266e-04 & 5.044983e-03 & 3.246915e-03 & 1.361348e-02 & 61 & 30 & 78 & 50 \\
3    & 9.607130e-04 & 4.791098e-03 & 3.492929e-03 & 1.307834e-02 & 56 & 30 & 71 & 50 \\
5    & 9.806757e-04 & 3.937232e-03 & 3.453308e-03 & 1.074114e-02 & 49 & 30 & 62 & 50 \\
7    & 9.715767e-04 & 3.231606e-03 & 2.952809e-03 & 8.199389e-03 & 45 & 30 & 57 & 50 \\
\bottomrule
\end{tabular}
\end{adjustbox}
\caption{Adaptive (\(^\text{a}\)) versus fixed-dimension (\(^\text{f}\)) ROM errors and dimensions for varying \(\hR\) (\(\hL = 25\), \(\epsloc = 4.0 \times 10^{-3}\); fixed dimensions: \(\rh = 30\), \(\rqq = 50\)).}
\label{tab:exp2_combined}
\end{table}

\subsubsection{Impact of Parameter Grid Distribution for $\hright$}
\label{sec:313}

The 1D dam-break problem with a wet-bed exhibits singular behavior as the wet-bed depth approaches zero; see the Appendix section. This motivates grid refinement near $\hright = 0$. To address this, we compare tROM performance using Chebyshev and uniformly distributed nodes over $[0, \hRmax]$.

We assess the effect of node distribution on tROM accuracy and ROM dimensions by varying the number of nodes $\NR \in \{5, 9, 17\}$ for $\hright \in [0, 8]$, while fixing $\NL = 13$ for $\hleft \in [10, 28]$. In the offline phase, the $(\NL, \NR)$ grid defines the training set.  In the online phase, we fix $\hleft = 25\,\text{m}$, $\epsloc = 4.0 \times 10^{-3}$, and evaluate performance over 80 $\hright$ values in $[0, 8]$.

\begin{figure}
\centering
\begin{subfigure}{0.48\linewidth}
    \centering
    \includegraphics[width=\linewidth]{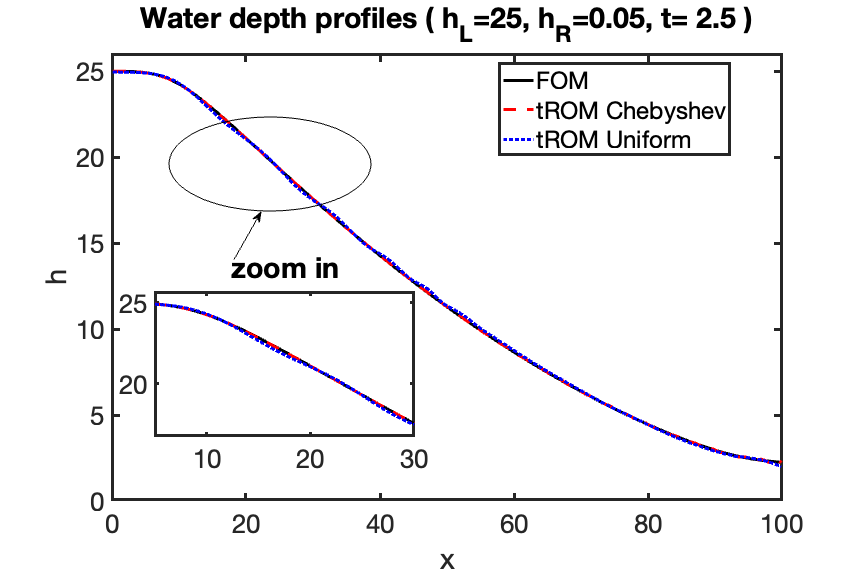}
    \subcaption{}
    \label{fig:exp3_h_profiles_0.05_combine}
\end{subfigure}
\hfill
\begin{subfigure}{0.48\linewidth}
    \centering
    \includegraphics[width=\linewidth]{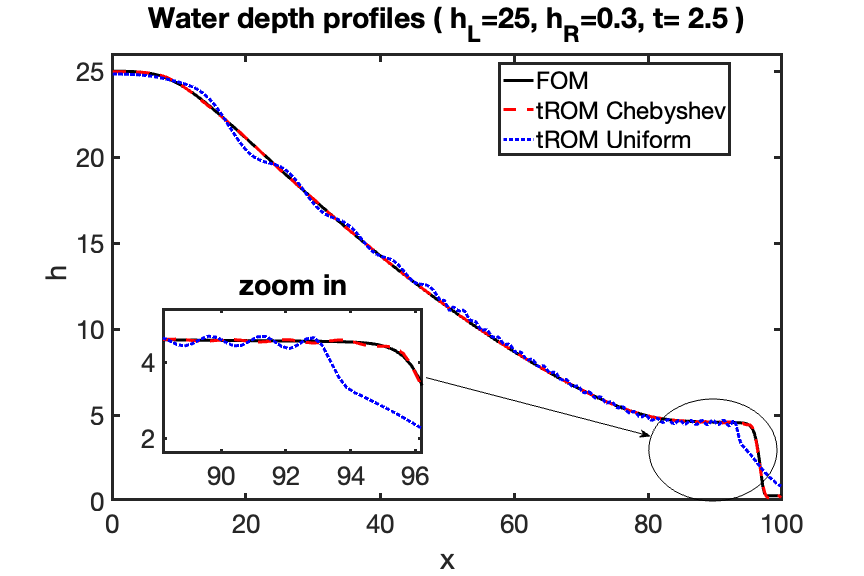}
    \subcaption{}
    \label{fig:exp3_h_profiles_0.3_combine}
\end{subfigure}

\vspace{0.5em}

\begin{subfigure}{0.48\linewidth}
    \centering
    \includegraphics[width=\linewidth]{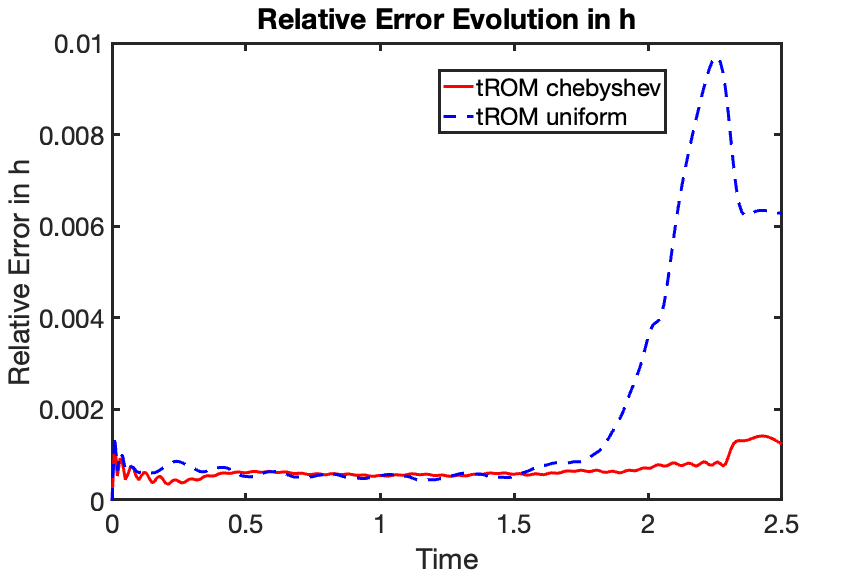}
    \subcaption{}
    \label{fig:exp3_h_errors_0.05}
\end{subfigure}
\hfill
\begin{subfigure}{0.48\linewidth}
    \centering
    \includegraphics[width=\linewidth]{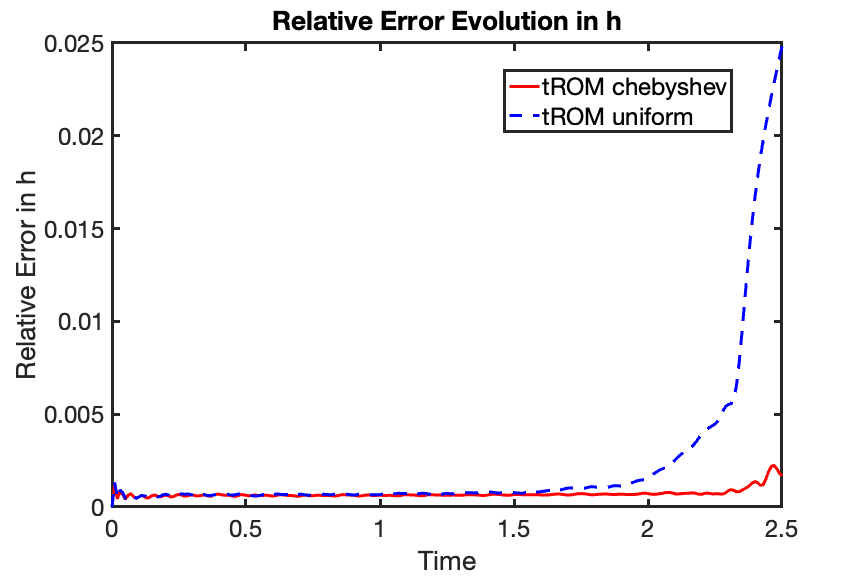}
    \subcaption{}
    \label{fig:exp3_h_errors_0.3}
\end{subfigure}

\caption{
Comparison of water depth profiles and relative $L^2$ errors for $\hleft = 25\,\text{m}$ using Chebyshev and uniform node distributions ($N_L = 13$, $N_R = 9$) in the non-interpolatory tROM:
(a) final-time depth profiles for $\hright = 0.05\,\text{m}$;
(b) final-time depth profiles for $\hright = 0.3\,\text{m}$;
(c) relative $L^2(0,T;L^2(\Omega))$ error over time for $\hright = 0.05\,\text{m}$;
(d) same error metric for $\hright = 0.3\,\text{m}$.
}
\label{fig:exp3_h_profiles_errors}
\end{figure}

We compute aggregate error metrics across these 80 points, along with mean and maximum ROM dimensions ($\rh$, $\rqq$). Results are summarized in Table~\ref{tab:exp3_nodes}. To highlight differences, we plot water-depth evolution and relative errors for selected cases with $\NL = 13$, $\NR = 9$, $\hleft = 25\,\text{m}$, and $\hright \in \{0.05, 0.3\}\,\text{m}$ (see Figure~\ref{fig:exp3_h_profiles_errors}).

The key finding is that Chebyshev nodes consistently outperform uniform ones, yielding lower errors across all $\NR$. Uniform grids underperform near small $\hright$, where solution discontinuities dominate, leading to increased errors (Figures~\ref{fig:exp3_h_errors_0.05} and \ref{fig:exp3_h_errors_0.3}). Chebyshev nodes, by clustering points near $\hright = 0$, better resolve these features and improve both accuracy and stability in the non-interpolatory tROM.

\renewcommand{\thetable}{4}
\begin{table}
\centering
\begin{adjustbox}{max width=\linewidth}
\begin{tabular}{@{}c l
                S[table-format=1.2e-2]
                S[table-format=1.2e-2]
                S[table-format=1.2e-2]
                S[table-format=1.2e-2]
                c c c c@{}}
\toprule
\(\NR\) & tROM type &
\multicolumn{1}{c}{\(\Esup_{\LtwoLtwo}(\hleft)\)} &
\multicolumn{1}{c}{\(\Eint_{\LtwoLtwo}(\hleft) \)} &
\multicolumn{1}{c}{\(\Esup_{\LtwoHone}(\hleft)\)} &
\multicolumn{1}{c}{\(\Eint_{\LtwoHone}(\hleft)\)} &
Mean \(\rh\) & Mean \(\rqq\) & Max \(\rh\) & Max \(\rqq\) \\
\midrule
5  & Chebyshev & 3.48027e-03 & 1.14203e-02 & 5.91490e-03 & 2.82674e-02 & 57 & 73 & 70 & 85 \\
5  & Uniform   & 6.67021e-03 & 1.32519e-02 & 1.01911e-02 & 3.65460e-02 & 55 & 71 & 64 & 83 \\
\addlinespace[0.3em]
9  & Chebyshev & 1.59533e-03 & 7.61879e-03 & 3.72022e-03 & 2.50259e-02 & 55 & 70 & 68 & 89 \\
9  & Uniform   & 5.14644e-03 & 9.15591e-03 & 7.94153e-03 & 3.17468e-02 & 54 & 69 & 67 & 85 \\
\addlinespace[0.3em]
17 & Chebyshev & 1.59532e-03 & 7.51281e-03 & 3.53661e-03 & 2.52534e-02 & 52 & 67 & 65 & 86 \\
17 & Uniform   & 2.92422e-03 & 8.05167e-03 & 5.09794e-03 & 3.08203e-02 & 52 & 67 & 65 & 85 \\
\bottomrule
\end{tabular}
\end{adjustbox}
\caption{Comparison of Chebyshev vs.\ uniform node distributions for varying \(\NR\)
         (\(\hL = 25\), \(\epsloc = 4.0 \times 10^{-3}\)).}
\label{tab:exp3_nodes}
\end{table}

\subsubsection{Interpolatory vs. Non-Interpolatory tROM}
\label{sec:314}

In this experiment, we compare the performance of interpolatory and non-interpolatory tROMs using Chebyshev nodes for the 1D dam-break problem. We fix $\hL = 25 \,\text{m}$, $\epsloc = 4.0 \times 10^{-3}$, and $\NL = 13$, and vary $\NR \in {5, 9, 17}$. In the offline phase, $\NL$ and $\NR$ define the parameter grid for $\hL$ and $\hR$, respectively. For each $(\NL,\NR)$ pair, we generate snapshot data using the FOM, which is then stored for use in the online phase. During the online phase, we evaluate the ROMs over a fine grid in $\hR \in [0, 8]$, computing aggregated error metrics (as defined in~\eqref{eq:error_def}) along with the mean and maximum ROM dimensions ($\rh$, $\rqq$).

\begin{figure}
\centering
\begin{subfigure}{0.48\linewidth}
    \centering
    \includegraphics[width=\linewidth]{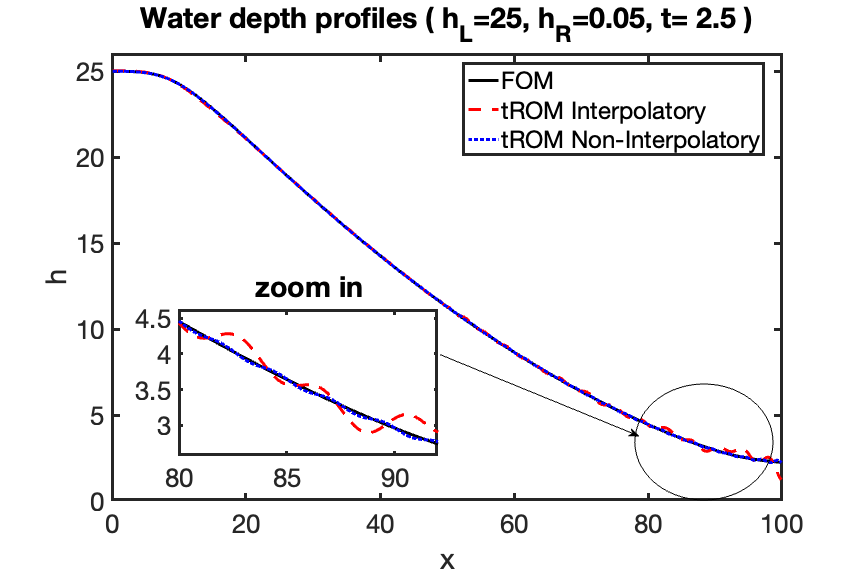}
    \subcaption{}
    \label{fig:exp4_h_evolution_hR_0.05_combine}
\end{subfigure}
\hfill
\begin{subfigure}{0.48\linewidth}
    \centering
    \includegraphics[width=\linewidth]{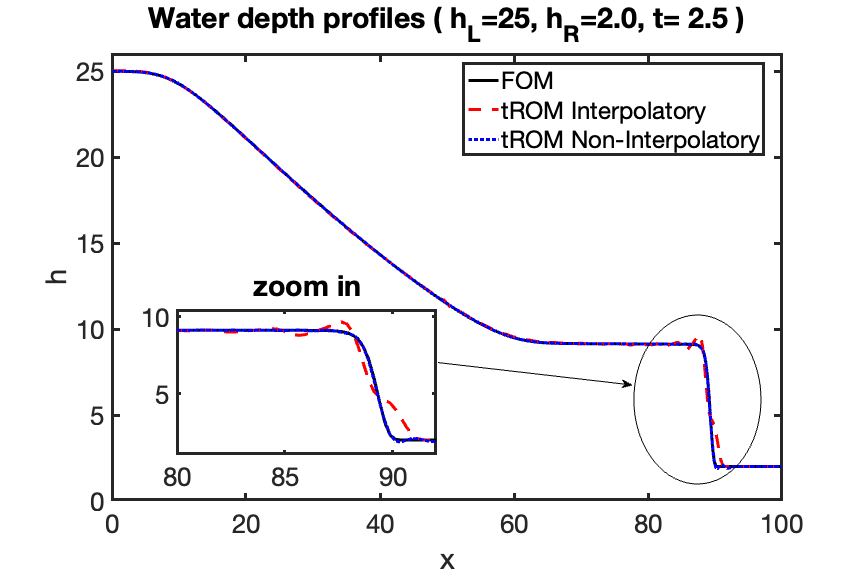}
    \subcaption{}
    \label{fig:exp4_h_evolution_hR_2_combine}
\end{subfigure}

\vspace{0.5em}

\begin{subfigure}{0.48\linewidth}
    \centering
    \includegraphics[width=\linewidth]{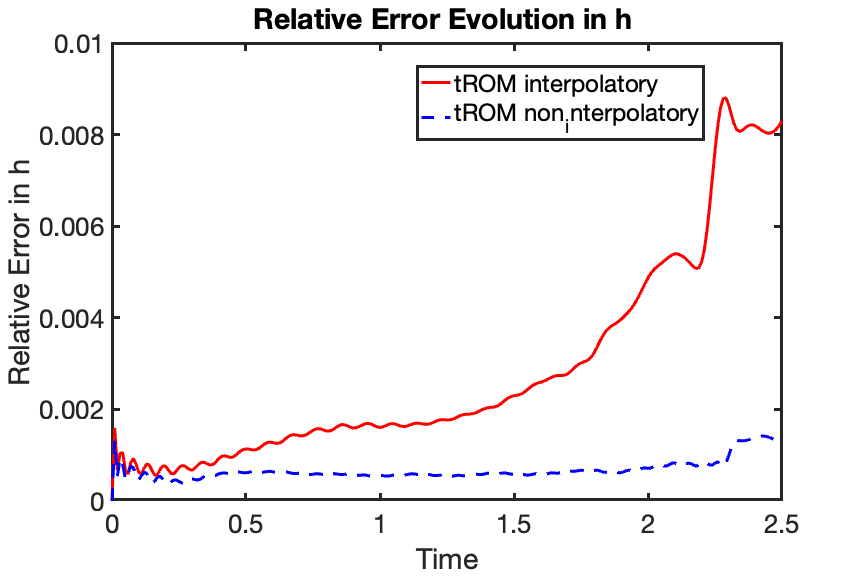}
    \subcaption{}
    \label{fig:exp4_h_errors_hR_0.05}
\end{subfigure}
\hfill
\begin{subfigure}{0.48\linewidth}
    \centering
    \includegraphics[width=\linewidth]{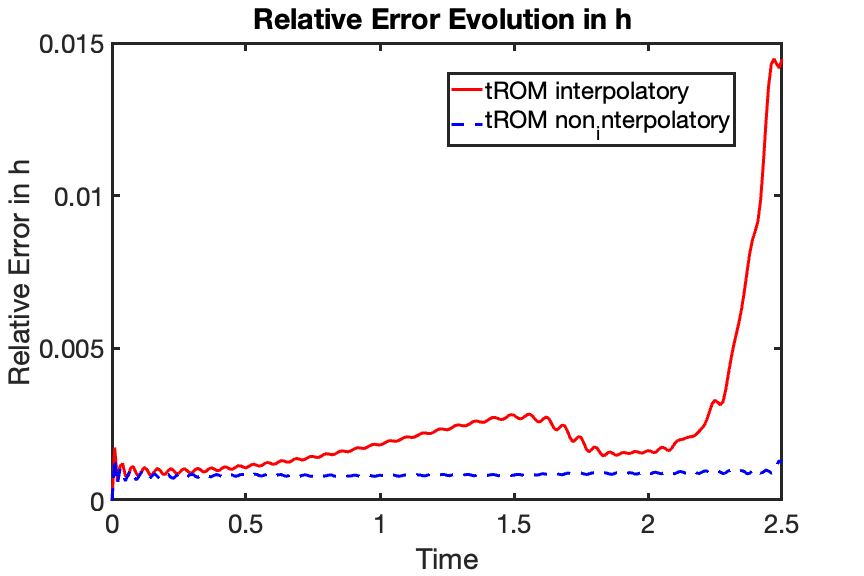}
    \subcaption{}
    \label{fig:exp4_h_errors_hR_2}
\end{subfigure}

\caption{Water depth profiles and error evolution over time for specific parameter pairs (\(N_L = 13, N_R = 9, \hleft = 25 \, \text{m}\)), comparing interpolatory and non-interpolatory tROMs : (a) water depth  profiles at final time  for \(\hright = 0.05 \, \text{m}\); (b) water depth  profiles at final time  for \(\hright = 2 \, \text{m}\); (c) relative error evolution in over time for \(\hright = 0.05 \, \text{m}\); (d)  relative error evolution   over time for \(\hright = 2.00 \, \text{m}\).}
\label{fig:exp4_h_evolution}
\end{figure}

The results are summarized in Table~\ref{tab:exp4_troms}. To further illustrate performance differences, we present simulation plots for the specific case $\NL = 13$, $\NR = 9$, $\hL = 25 \, \text{m}$, and $\hR \in {0.05, 2}\, \text{m}$, showing the evolution of water depth and relative errors in $h$ over time (see Figures~\ref{fig:exp4_h_evolution}).

We find that the non-interpolatory tROM consistently outperforms the interpolatory tROM across all tested values of $\NR$. Although it results in higher ranks, the increased dimensionality enables better resolution of solution variability.

\renewcommand{\thetable}{5}
\begin{table}
\centering
\begin{adjustbox}{max width=\linewidth}
\begin{tabular}{@{}c l
                S[table-format=1.2e-2]
                S[table-format=1.2e-2]
                S[table-format=1.2e-2]
                S[table-format=1.2e-2]
                c c c c@{}}
\toprule
\(\NR\) & tROM type &
\multicolumn{1}{c}{\(\Esup_{\LtwoLtwo}(\hleft)\)} &
\multicolumn{1}{c}{\(\Eint_{\LtwoLtwo}(\hleft) \)} &
\multicolumn{1}{c}{\(\Esup_{\LtwoHone}(\hleft)\)} &
\multicolumn{1}{c}{\(\Eint_{\LtwoHone}(\hleft)\)} &
Mean \(\rh\) & Mean \(\rqq\) & Max \(\rh\) & Max \(\rqq\) \\
\midrule
5  & Interpolatory     & 2.08076e-02 & 2.86495e-02 & 3.25643e-02 & 5.74810e-02 & 43 & 56 & 53 & 70 \\
5  & Non-interpolatory & 3.48027e-03 & 1.14203e-02 & 5.91490e-03 & 2.82674e-02 & 57 & 73 & 70 & 85 \\
\addlinespace[0.3em]
9  & Interpolatory     & 1.01665e-02 & 1.47304e-02 & 1.78460e-02 & 4.09531e-02 & 43 & 56 & 53 & 70 \\
9  & Non-interpolatory & 1.59533e-03 & 7.61879e-03 & 3.72022e-03 & 2.50259e-02 & 55 & 70 & 68 & 89 \\
\addlinespace[0.3em]
17 & Interpolatory     & 3.39686e-03 & 1.05460e-02 & 8.03443e-03 & 3.42998e-02 & 44 & 56 & 53 & 70 \\
17 & Non-interpolatory & 1.59532e-03 & 7.51281e-03 & 3.53661e-03 & 2.52534e-02 & 52 & 67 & 65 & 85 \\
\bottomrule
\end{tabular}
\end{adjustbox}
\caption{Interpolatory vs.\ non-interpolatory tROMs  for varying \(\NR\) (\(\hL = 25\), \(\epsloc = 4.0 \times 10^{-3}\)).}
\label{tab:exp4_troms}
\end{table}

\subsubsection{tROM vs POD ROM}

In this experiment, we compare the performance of the non-interpolatory tROM to that of the POD ROM \rev{for two scenarios. First, we fix the threshold 
$\epsloc$, and, second, we use the reduced bases with the same fixed dimensions.}

\textbf{Comparison with fixed $\epsloc$.}
We fix $\epsloc = 4.0 \times 10^{-3}$ and vary the parameter grid sizes with $\NL \in {4, 7, 13, 20}$ and $\NR \in {5, 9, 17, 25}$. Therefore, ROM dimensions may vary depending on the parameter grid used in each simulation. 
\rev{To ensure a consistent comparison, we use the same ROM dimensions for the POD ROM and tROM. The dimensions of the two bases are obtained from the tROM online phase (with $\epsloc = 4.0 \times 10^{-3}$), and the same dimensions are then used for POD ROM.}

\begin{figure}
\centering
\begin{subfigure}{0.48\linewidth}
    \centering
    \includegraphics[width=\linewidth]{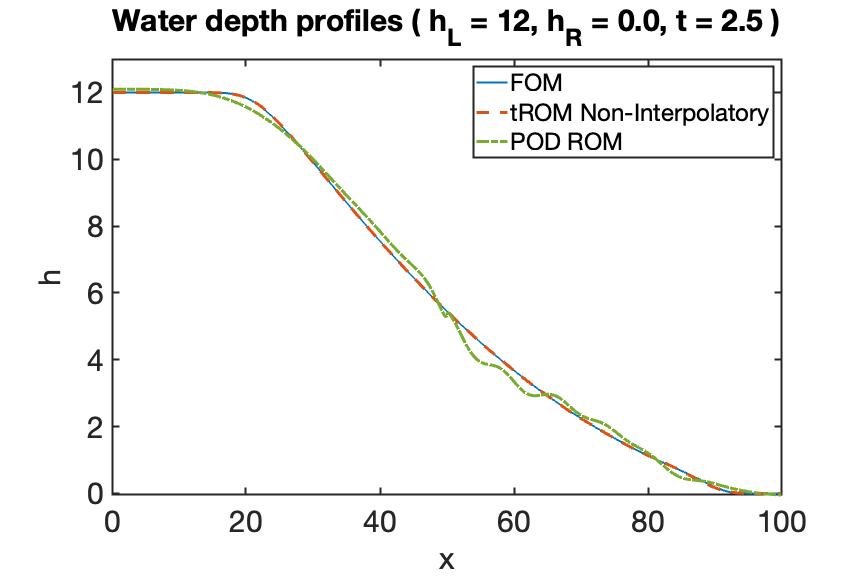}
    \subcaption{}
    \label{fig:exp7_hcomp_hL_12_hR_0}
\end{subfigure}
\hfill
\begin{subfigure}{0.48\linewidth}
    \centering
    \includegraphics[width=\linewidth]{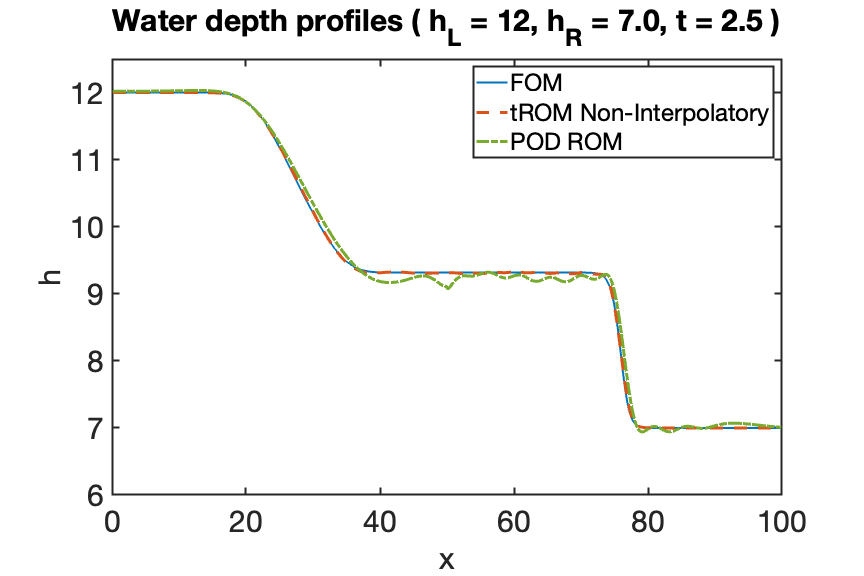}
    \subcaption{}
    \label{fig:exp7_hcomp_hL_12_hR_7}
\end{subfigure}

\vspace{0.5em}

\begin{subfigure}{0.48\linewidth}
    \centering
    \includegraphics[width=\linewidth]{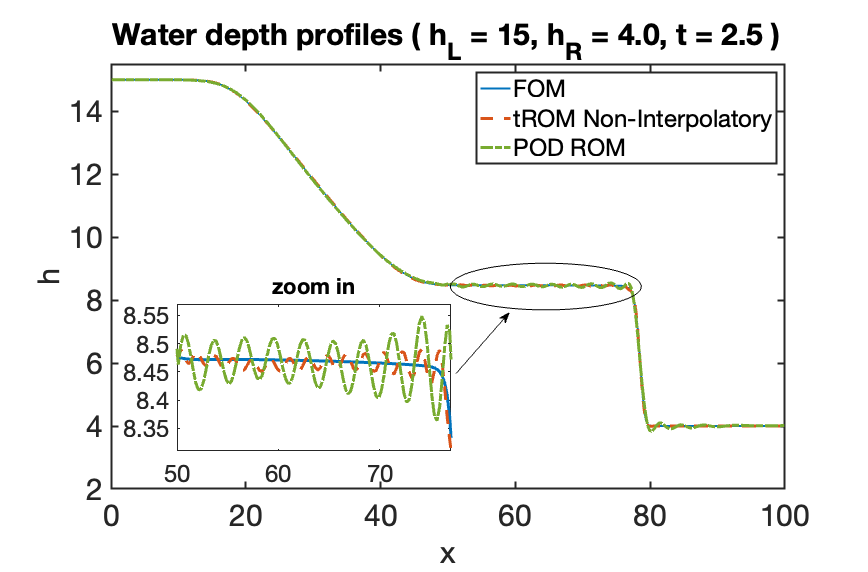}
    \subcaption{}
    \label{fig:exp7_hcomp_hL_15_hR_4}
\end{subfigure}
\hfill
\begin{subfigure}{0.48\linewidth}
    \centering
    \includegraphics[width=\linewidth]{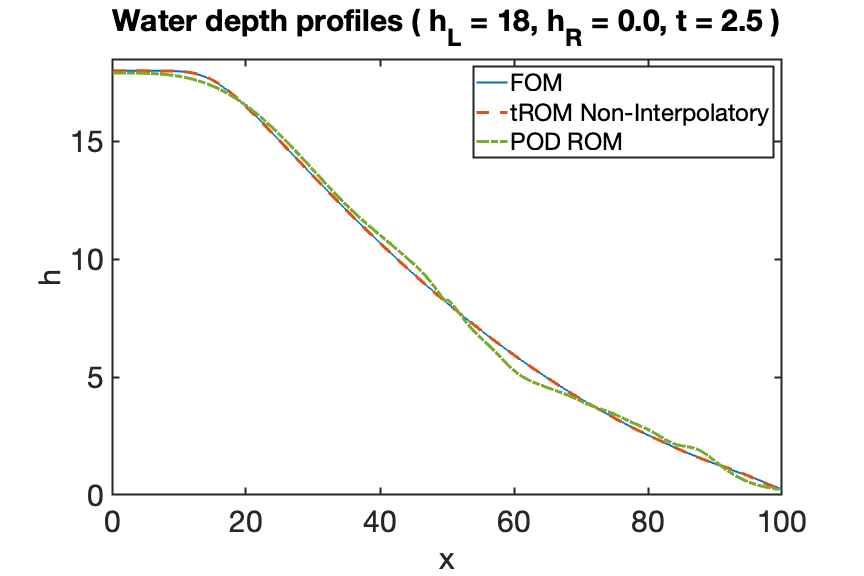}
    \subcaption{}
    \label{fig:exp7_hcomp_hL_18_hR_0}
\end{subfigure}

\vspace{0.5em}

\begin{subfigure}{0.48\linewidth}
    \centering
    \includegraphics[width=\linewidth]{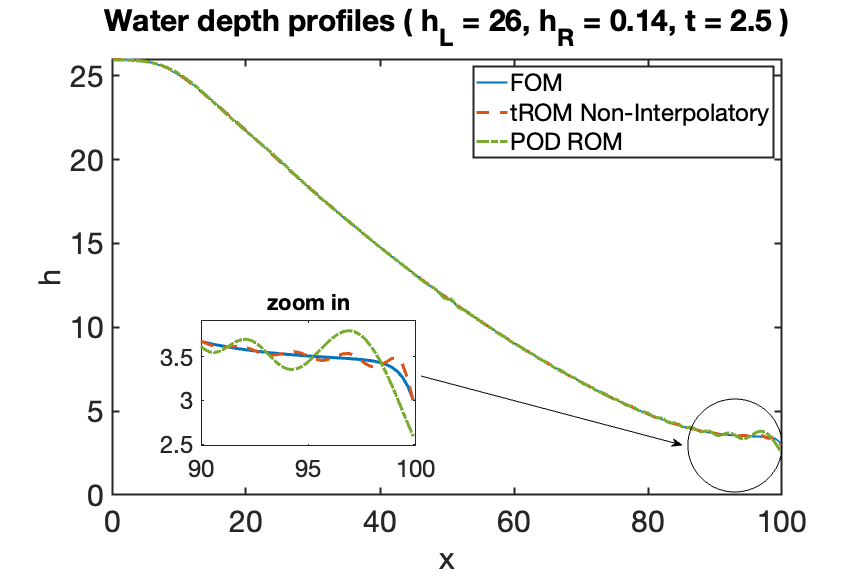}
    \subcaption{}
    \label{fig:exp7_hcomp_hL_26_hR_0.14}
\end{subfigure}
\hfill
\begin{subfigure}{0.48\linewidth}
    \centering
    \includegraphics[width=\linewidth]{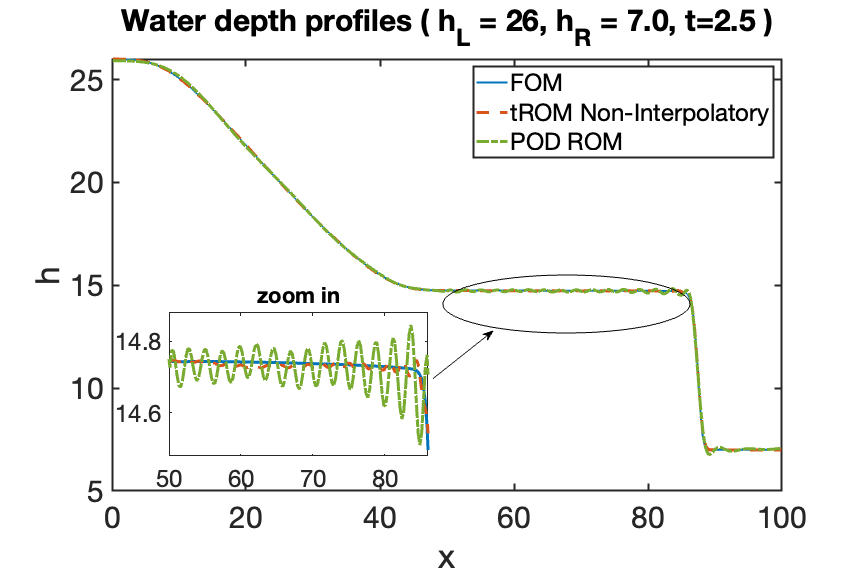}
    \subcaption{}
    \label{fig:exp7_hcomp_hL_26_hR_7}
\end{subfigure}

\caption{Water depth profiles at final time for different parameter pairs, comparing non-interpolatory tROM and POD ROM (with equal ROM dimensions derived from tROM using \(\epsilon_{\text{loc}} = 4.0 \times 10^{-3}\)) with Chebyshev nodes: (a) for \(\hleft = 12, \hright = 0\); (b) for \(\hleft = 12, \hright = 7\); (c) for \(\hleft = 15, \hright = 4\); (d) for \(\hleft = 18, \hright = 0\); (e) for \(\hleft = 26, \hright = 0.14\); (f) for \(\hleft = 26, \hright = 7\).}
\label{fig:exp7_h_comparison}
\end{figure}

\begin{figure}
\centering
\begin{subfigure}{0.48\linewidth}
    \centering
    \includegraphics[width=\linewidth]{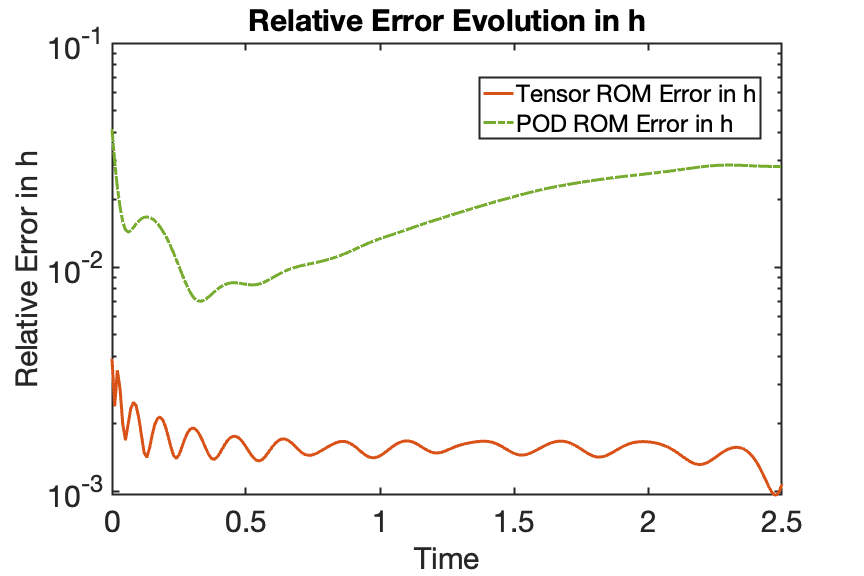}
    \subcaption{}
    \label{fig:exp7_error_hL_12_hR_0}
\end{subfigure}
\hfill
\begin{subfigure}{0.48\linewidth}
    \centering
    \includegraphics[width=\linewidth]{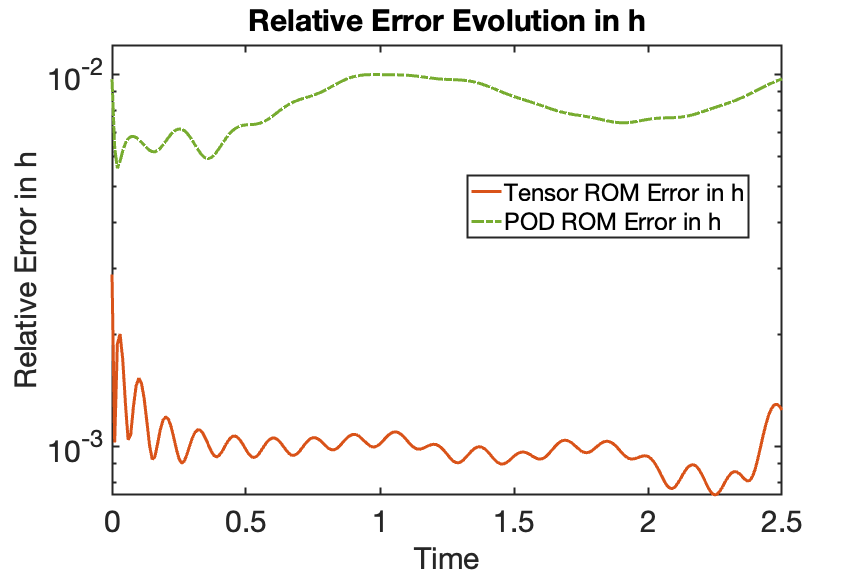}
    \subcaption{}
    \label{fig:exp7_error_hL_12_hR_7}
\end{subfigure}

\vspace{0.5em}

\begin{subfigure}{0.48\linewidth}
    \centering
    \includegraphics[width=\linewidth]{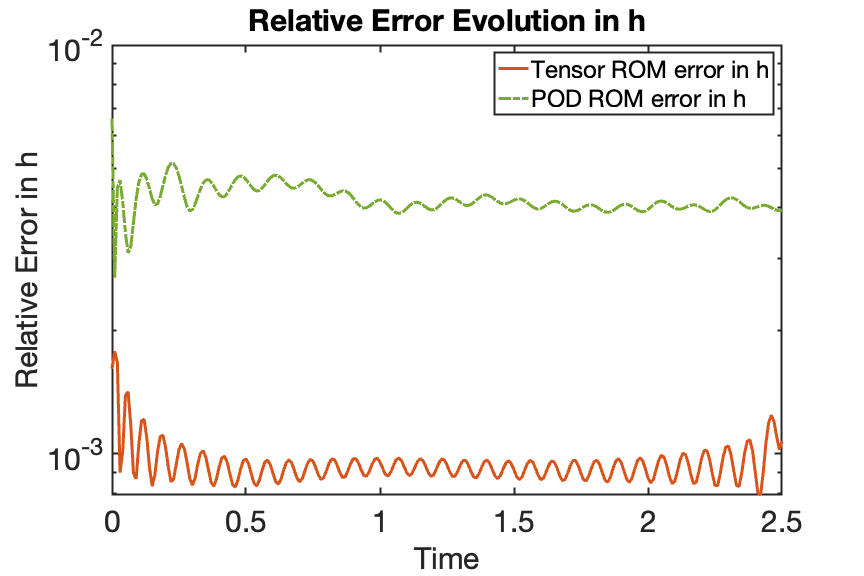}
    \subcaption{}
    \label{fig:exp7_error_hL_15_hR_4}
\end{subfigure}
\hfill
\begin{subfigure}{0.48\linewidth}
    \centering
    \includegraphics[width=\linewidth]{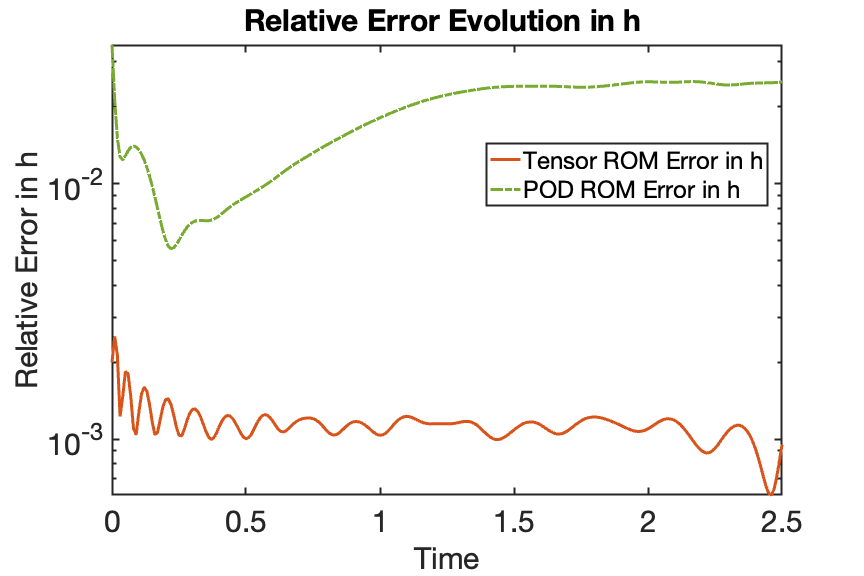}
    \subcaption{}
    \label{fig:exp7_error_hL_18_hR_0}
\end{subfigure}

\vspace{0.5em}
\begin{subfigure}{0.48\linewidth}
    \centering
    \includegraphics[width=\linewidth]{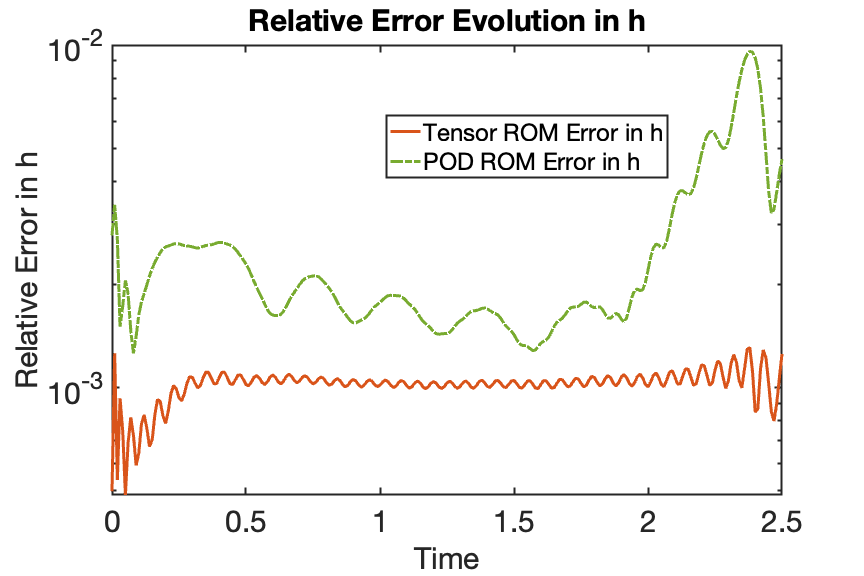}
    \subcaption{}
    \label{fig:exp7_error_hL_26_hR_0.14}
\end{subfigure}
\hfill
\begin{subfigure}{0.48\linewidth}
    \centering
    \includegraphics[width=\linewidth]{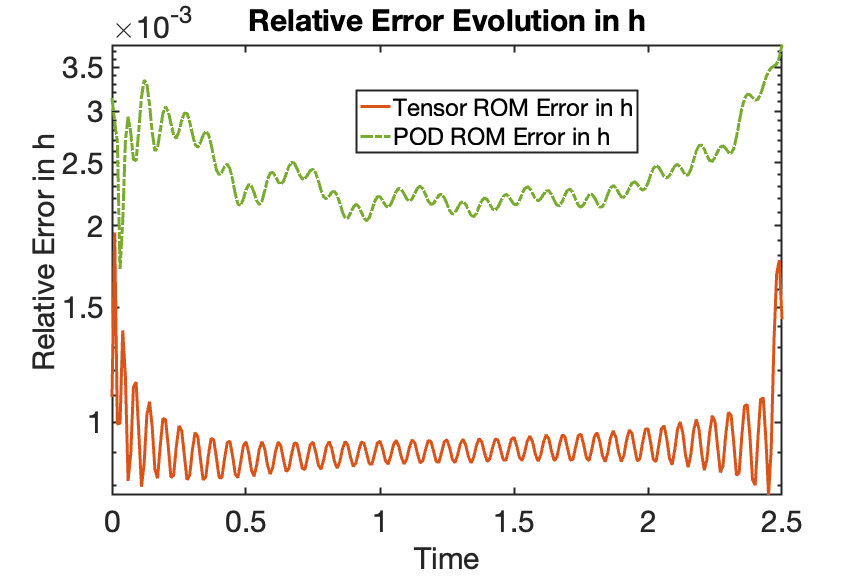}
    \subcaption{}
    \label{fig:exp7_error_hL_26_hR_7}
\end{subfigure}

\caption{Evolution of relative errors in water depth   for different parameter pairs, comparing non-interpolatory tROM and POD ROM (with equal ROM dimensions derived from tROM using \(\epsilon_{\text{loc}} = 4.0 \times 10^{-3}\)) with Chebyshev nodes : (a) for \(\hleft = 12, \hright = 0\); (b) for \(\hleft = 12, \hright = 7\); (c) for \(\hleft = 15, \hright = 4\); (d) for \(\hleft = 18, \hright = 0\); (e) for \(\hleft = 26, \hright = 0.14\); (f) for \(\hleft = 26, \hright = 7\).}
\label{fig:exp7_error_evolution}
\end{figure}

\rev{Water depth and discharge profiles are presented in Figures \ref{fig:exp7_h_comparison} and \ref{fig:exp7_hu_comparison}.
The time-evolution of relative errors is presented in Figures \ref{fig:exp7_error_evolution} and 
\ref{fig:exp7_error_hu_evolution}.}
Numerical errors are also summarized in Table~\ref{tab:exp5_trom_pod}. Our results show that the non-interpolatory tROM consistently outperforms the POD ROM across all tested values of $\NL$ and $\NR$. This improvement is attributed to the tensor \rev{structure, which enables the tROM} to better adapt to the solution variability across the parameter domain.

\renewcommand{\thetable}{6} 
\begin{table}
\centering
\begin{adjustbox}{max width=\linewidth}
\begin{tabular}{@{}cc l
                S[table-format=1.2e-2]
                S[table-format=1.2e-2]
                S[table-format=1.2e-2]
                S[table-format=1.2e-2]
                c
                c
                c
                c@{}}
\toprule
\(\NL\) & \(\NR\) & ROM Type &
\multicolumn{1}{c}{\(\Esup_{\LtwoLtwo}\)} &
\multicolumn{1}{c}{\(\Eint_{\LtwoLtwo}\)} &
\multicolumn{1}{c}{\(\Esup_{\LtwoHone}\)}  &
\multicolumn{1}{c}{\(\Eint_{\LtwoHone}\)} &
\multicolumn{1}{c}{Mean \(\rh\)} &
{Max \(\rh\)} &
\multicolumn{1}{c}{Mean \(\rqq\)} &
{Max \(\rqq\)} \\
\midrule
4  & 5  & tROM & 3.32482e-03 & 1.40384e-01 & 6.40727e-03 & 4.15977e-01 & 50 & 72 & 67 & 92 \\
4  & 5  & POD  & 4.51831e-03 & 3.28964e-01 & 1.08839e-02 & 8.20292e-01 & 50 & 72 & 67 & 92 \\
\addlinespace[0.3em]
7  & 9  & tROM & 1.47438e-03 & 1.40981e-01 & 4.12165e-03 & 4.32671e-01 & 45 & 70 & 60 & 90 \\
7  & 9  & POD  & 6.79836e-03 & 4.50186e-01 & 1.10868e-02 & 1.03866e+00 & 45 & 70 & 60 & 90 \\
\addlinespace[0.3em]
13 & 17 & tROM & 1.53420e-03 & 1.42511e-01 & 3.81702e-03 & 4.48196e-01 & 42 & 66 & 56 & 85 \\
13 & 17 & POD  & 1.02681e-02 & 5.44312e-01 & 1.36055e-02 & 1.18492e+00 & 42 & 66 & 56 & 85 \\
\addlinespace[0.3em]
20 & 25 & tROM & 1.35122e-03 & 1.44950e-01 & 3.72353e-03 & 4.55263e-01 & 41 & 65 & 54 & 84 \\
20 & 25 & POD  & 1.17029e-02 & 5.87390e-01 & 1.43979e-02 & 1.25797e+00 & 41 & 65 & 54 & 84 \\
\bottomrule
\end{tabular}
\end{adjustbox}
\caption{Comparison of tROM and POD ROM with Monte Carlo sampling for varying \(\NL\) and \(\NR\) (\(\epsloc = 4.0 \times 10^{-3}\)). \rev{Note that ROM dimensions are the same for POD ROM and tROM for a fixed value of $(N_L,N_R)$.}}
\label{tab:exp5_trom_pod}
\end{table}

\textbf{Comparison with fixed ROM dimensions.}
We further compare the non-interpolatory tROM to the POD ROM using fixed ROM dimensions \(\rh = 30\) and \(\rqq = 50\) for both methods, instead of adapting dimensions based on the threshold \(\epsloc\). We vary the parameter grid sizes with \(\NL \in \{4, 7, 13, 20\}\) and \(\NR \in \{5, 9, 17, 25\}\). Results, summarized in Table~\ref{tab:exp6_trom_pod}, confirm that the non-interpolatory tROM consistently outperforms the POD ROM across all \(\NL\) and \(\NR\) even with fixed ROM dimensions.

\renewcommand{\thetable}{7}
\begin{table}\small 
\centering
\begin{adjustbox}{max width=\linewidth}
\begin{tabular}{@{}cc l
                S[table-format=1.2e-2]
                S[table-format=1.2e-2]
                S[table-format=1.2e-2]
                S[table-format=1.2e-2]@{}}
\toprule
\(\NL\) & \(\NR\) & ROM Type &
\multicolumn{1}{c}{\(\Esup_{\LtwoLtwo}\)} &
\multicolumn{1}{c}{\(\Eint_{\LtwoLtwo}\)} &
\multicolumn{1}{c}{\(\Esup_{\LtwoHone}\)} &
\multicolumn{1}{c}{\(\Eint_{\LtwoHone}\)} \\
\midrule
4  & 5  & tROM & 7.35522e-03 & 5.52622e-01 & 1.70112e-02 & 1.33668e+00 \\
4  & 5  & POD  & 9.73226e-03 & 9.18864e-01 & 1.82700e-02 & 1.93031e+00 \\
\addlinespace[0.3em]
7  & 9  & tROM & 6.47034e-03 & 4.45161e-01 & 1.54951e-02 & 1.16569e+00 \\
7  & 9  & POD  & 9.55300e-03 & 9.41798e-01 & 1.81970e-02 & 1.92800e+00 \\
\addlinespace[0.3em]
13 & 17 & tROM & 5.67702e-03 & 3.92979e-01 & 1.44389e-02 & 1.06244e+00 \\
13 & 17 & POD  & 9.74292e-03 & 9.43969e-01 & 1.80491e-02 & 1.90962e+00 \\
\addlinespace[0.3em]
20 & 25 & tROM & 5.47933e-03 & 3.71994e-01 & 1.42556e-02 & 1.02878e+00 \\
20 & 25 & POD  & 9.81101e-03 & 9.43520e-01 & 1.80041e-02 & 1.90319e+00 \\
\bottomrule
\end{tabular}
\end{adjustbox}
\caption{Comparison of tROM and POD ROM with Monte Carlo sampling for varying \(\NL\) and \(\NR\), fixed ROM dimensions \(\rh = 30\), \(\rqq = 50\).}
\label{tab:exp6_trom_pod}
\end{table}

Finally,  Figures~\ref{fig:exp7_h_comparison} and \ref{fig:exp7_error_evolution} visualize the water depth profiles at the final time and the evolution of relative errors for parameter pairs \((\hleft, \hright) \in \{(12, 0), (12, 7),(15, 4), (18, 0), (26, 0.14), (26, 7)\}\), representing a range of flow regimes. Again, we observe that the tROM outperforms the POD ROM in terms of solution accuracy and suppression of spurious oscillations, given that both methods use the same reduced dimensions.

\rev{
\begin{figure}
\centering
\begin{subfigure}{0.48\linewidth}
    \centering
    \includegraphics[width=\linewidth]{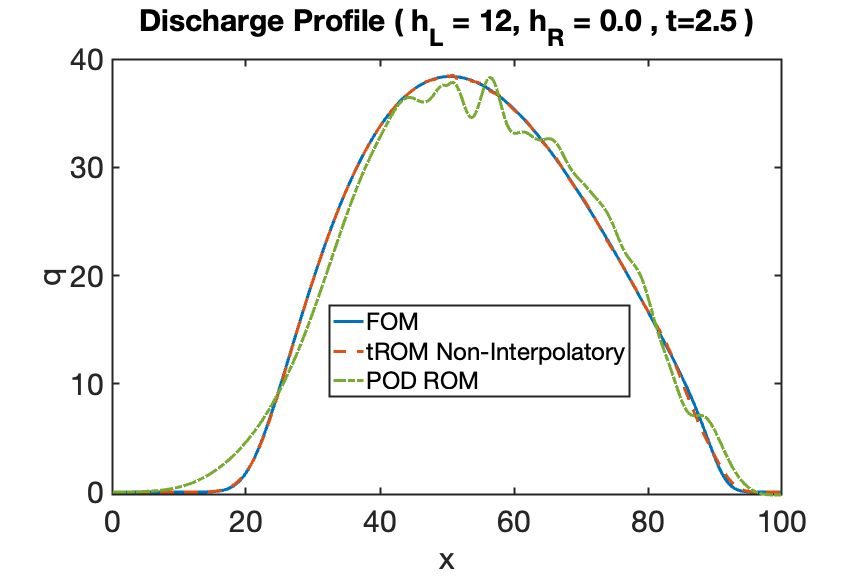}
    \subcaption{}
    \label{fig:exp7_hucomp_hL_12_hR_0}
\end{subfigure}
\hfill
\begin{subfigure}{0.48\linewidth}
    \centering
    \includegraphics[width=\linewidth]{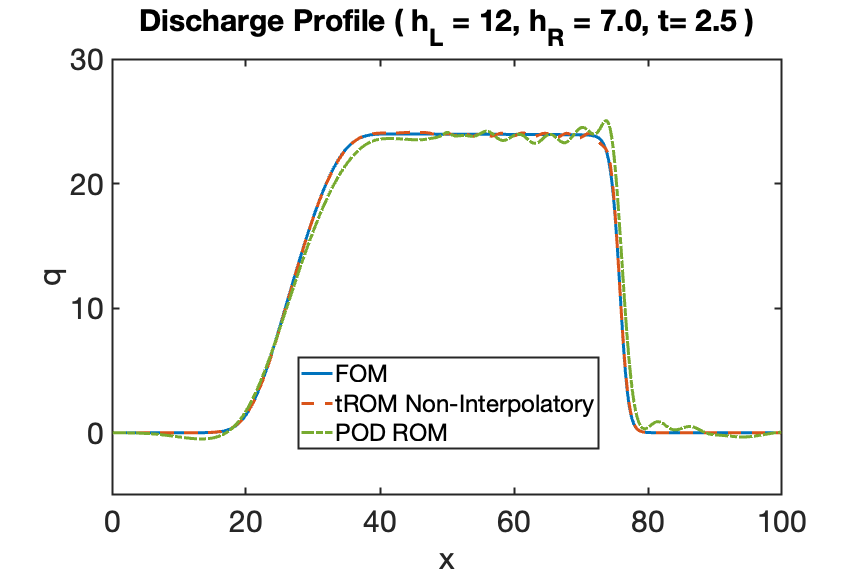}
    \subcaption{}
    \label{fig:exp7_hucomp_hL_12_hR_7}
\end{subfigure}

\vspace{0.5em}

\begin{subfigure}{0.48\linewidth}
    \centering
    \includegraphics[width=\linewidth]{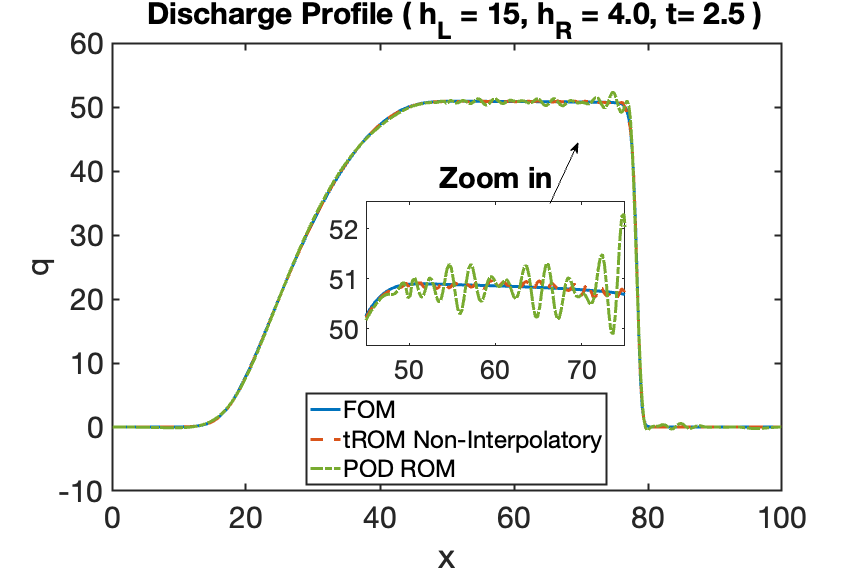}
    \subcaption{}
    \label{fig:exp7_hucomp_hL_15_hR_4}
\end{subfigure}
\hfill
\begin{subfigure}{0.48\linewidth}
    \centering
    \includegraphics[width=\linewidth]{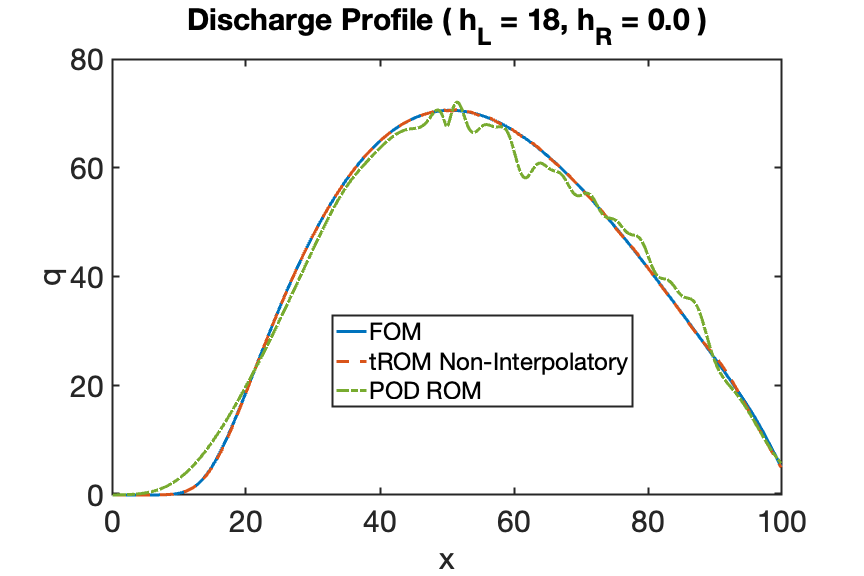}
    \subcaption{}
    \label{fig:exp7_hucomp_hL_18_hR_0}
\end{subfigure}

\vspace{0.5em}

\begin{subfigure}{0.48\linewidth}
    \centering
    \includegraphics[width=\linewidth]{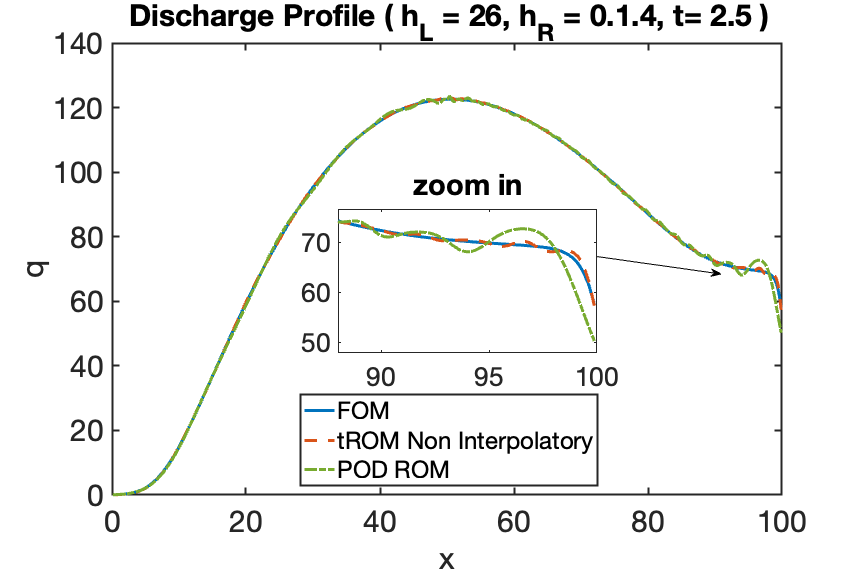}
    \subcaption{}
    \label{fig:exp7_hucomp_hL_26_hR_0.14}
\end{subfigure}
\hfill
\begin{subfigure}{0.48\linewidth}
    \centering
    \includegraphics[width=\linewidth]{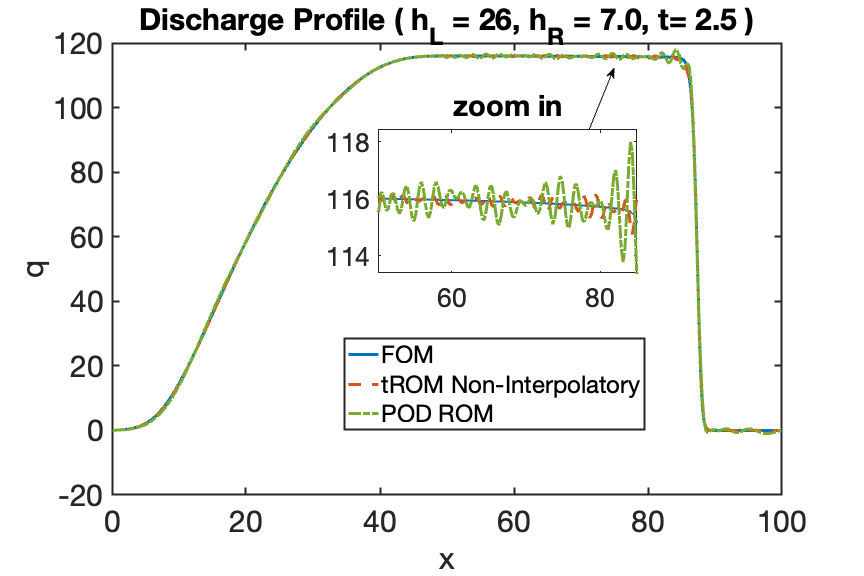}
    \subcaption{}
    \label{fig:exp7_hucomp_hL_26_hR_7}
\end{subfigure}

\caption{Discharge profiles at final time for different parameter pairs, comparing non-interpolatory tROM and POD ROM (with equal ROM dimensions derived from tROM using \(\epsilon_{\text{loc}} = 4.0 \times 10^{-3}\)) with Chebyshev nodes: (a) for \(\hleft = 12, \hright = 0\); (b) for \(\hleft = 12, \hright = 7\); (c) for \(\hleft = 15, \hright = 4\); (d) for \(\hleft = 18, \hright = 0\); (e) for \(\hleft = 26, \hright = 0.14\); (f) for \(\hleft = 26, \hright = 7\).}
\label{fig:exp7_hu_comparison}
\end{figure}

\begin{figure}
\centering
\begin{subfigure}{0.48\linewidth}
    \centering
    \includegraphics[width=\linewidth]{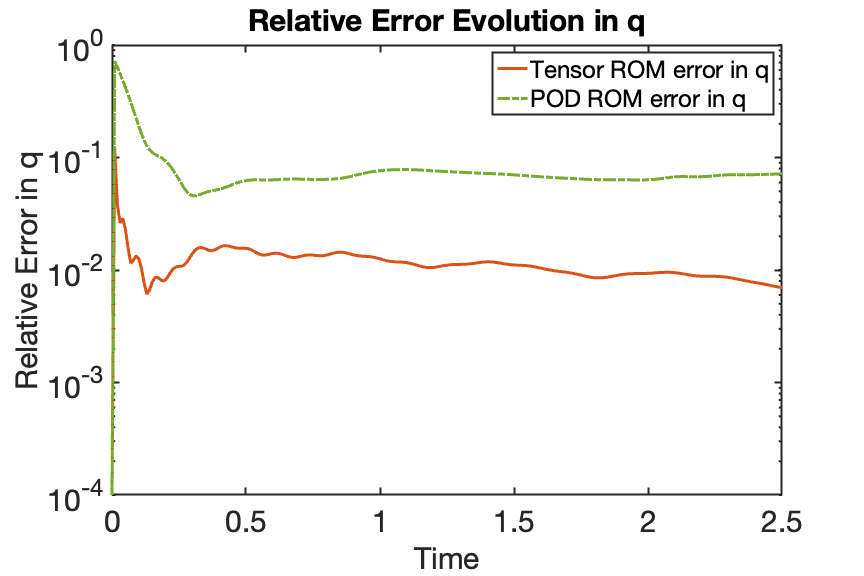}
    \subcaption{}
    \label{fig:exp7_error_hu_hL_12_hR_0}
\end{subfigure}
\hfill
\begin{subfigure}{0.48\linewidth}
    \centering
    \includegraphics[width=\linewidth]{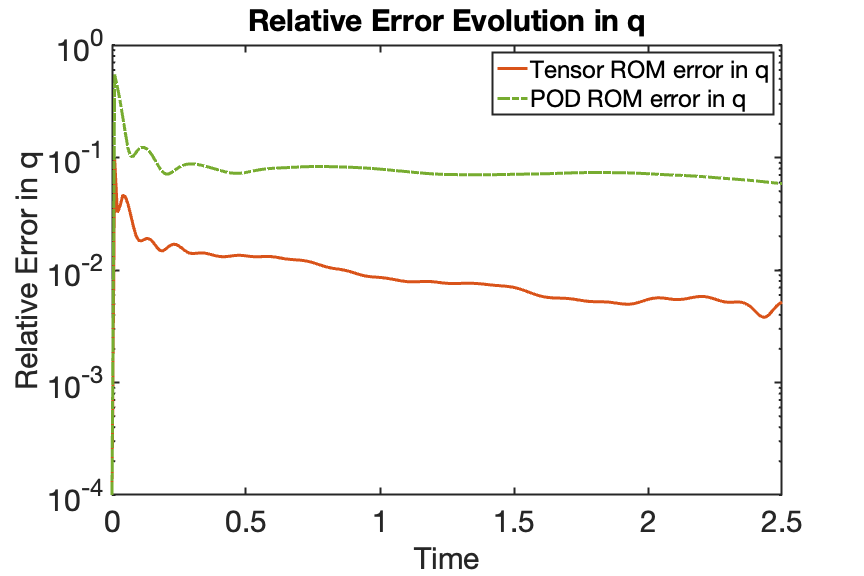}
    \subcaption{}
    \label{fig:exp7_error_hu_hL_12_hR_7}
\end{subfigure}

\vspace{0.5em}

\begin{subfigure}{0.48\linewidth}
    \centering
    \includegraphics[width=\linewidth]{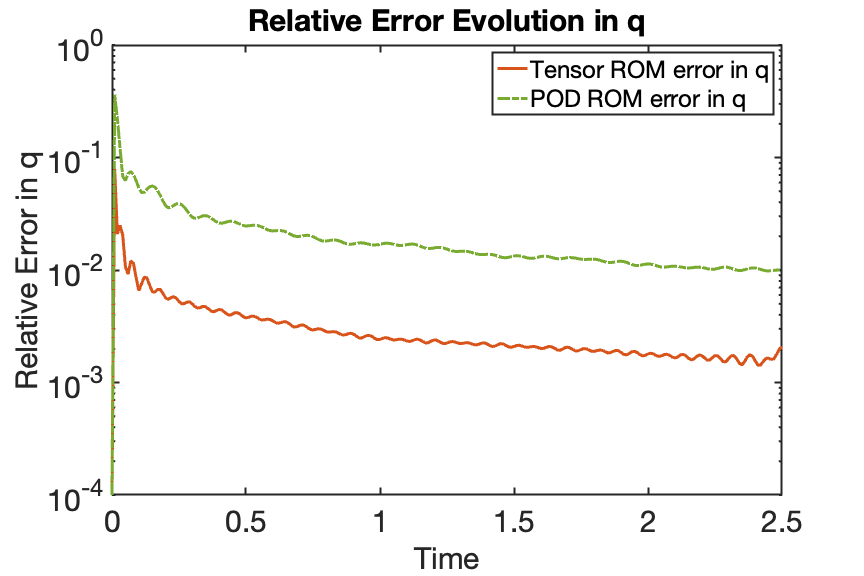}
    \subcaption{}
    \label{fig:exp7_error_hu_hL_15_hR_4}
\end{subfigure}
\hfill
\begin{subfigure}{0.48\linewidth}
    \centering
    \includegraphics[width=\linewidth]{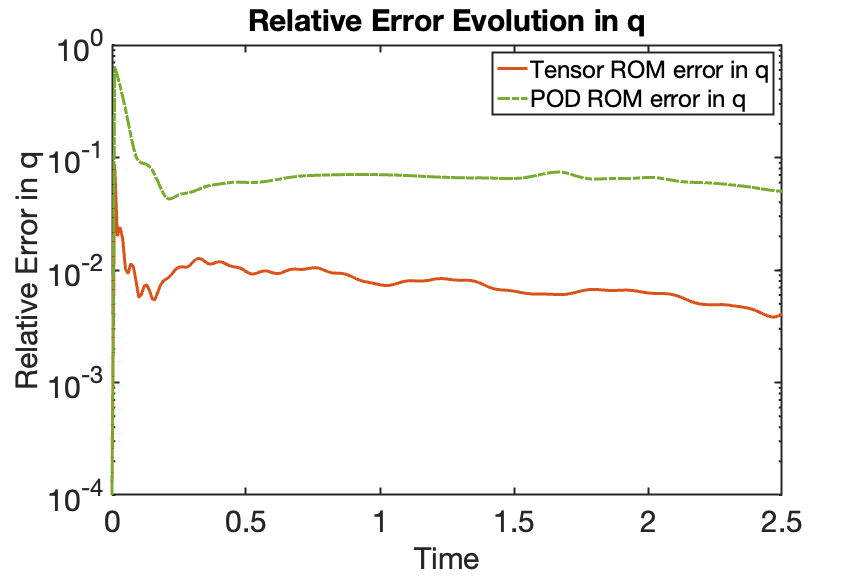}
    \subcaption{}
    \label{fig:exp7_error_hu_hL_18_hR_0}
\end{subfigure}

\vspace{0.5em}
\begin{subfigure}{0.48\linewidth}
    \centering
    \includegraphics[width=\linewidth]{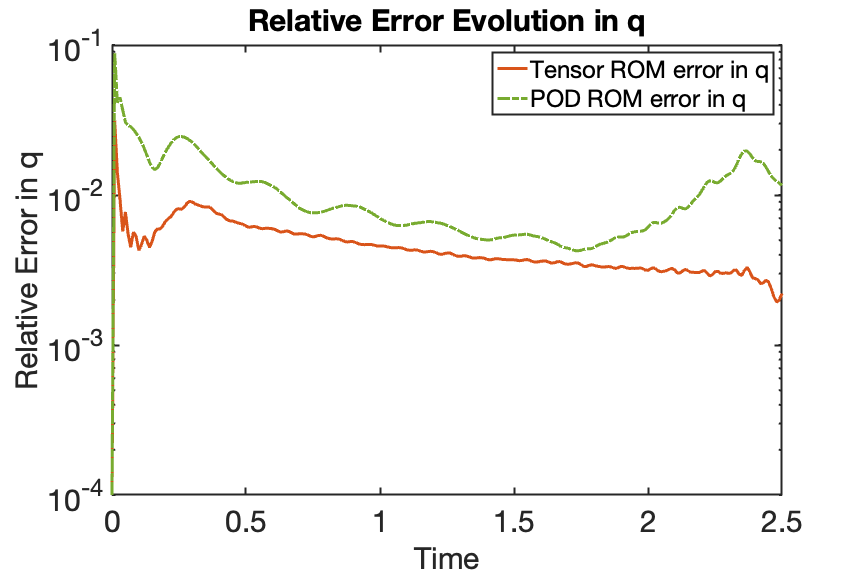}
    \subcaption{}
    \label{fig:exp7_error_hu_hL_26_hR_0.14}
\end{subfigure}
\hfill
\begin{subfigure}{0.48\linewidth}
    \centering
    \includegraphics[width=\linewidth]{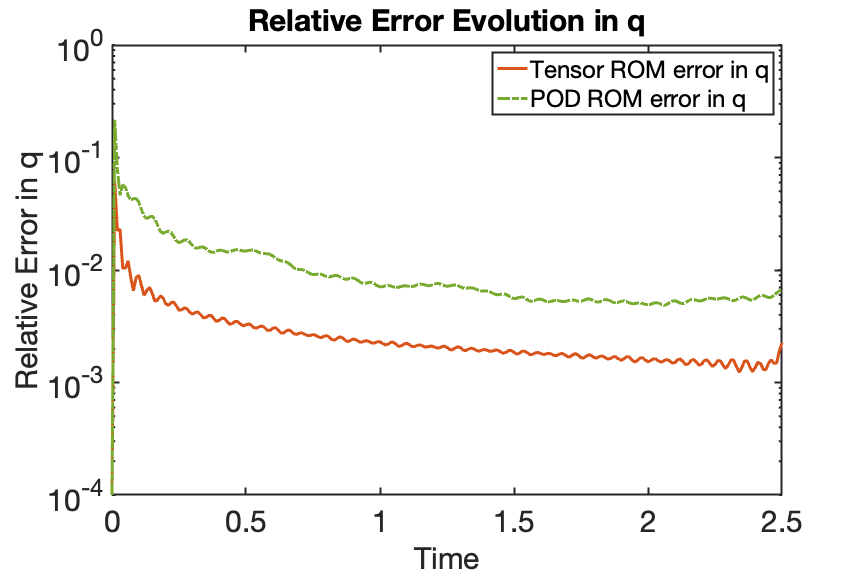}
    \subcaption{}
    \label{fig:exp7_error_hu_hL_26_hR_7}
\end{subfigure}

\caption{Evolution of relative errors in the discharge for different parameter pairs, comparing non-interpolatory tROM and POD ROM (with equal ROM dimensions derived from tROM using \(\epsilon_{\text{loc}} = 4.0 \times 10^{-3}\)) with Chebyshev nodes : (a) for \(\hleft = 12, \hright = 0\); (b) for \(\hleft = 12, \hright = 7\); (c) for \(\hleft = 15, \hright = 4\); (d) for \(\hleft = 18, \hright = 0\); (e) for \(\hleft = 26, \hright = 0.14\); (f) for \(\hleft = 26, \hright = 7\).}
\label{fig:exp7_error_hu_evolution}
\end{figure}}


\section{2D SWE Numerical Experiments} 
\label{sec:2d_experiments}
\begin{figure}[h]
\centering
\begin{subfigure}{0.48\linewidth}
    \centering   \includegraphics[width=\linewidth]{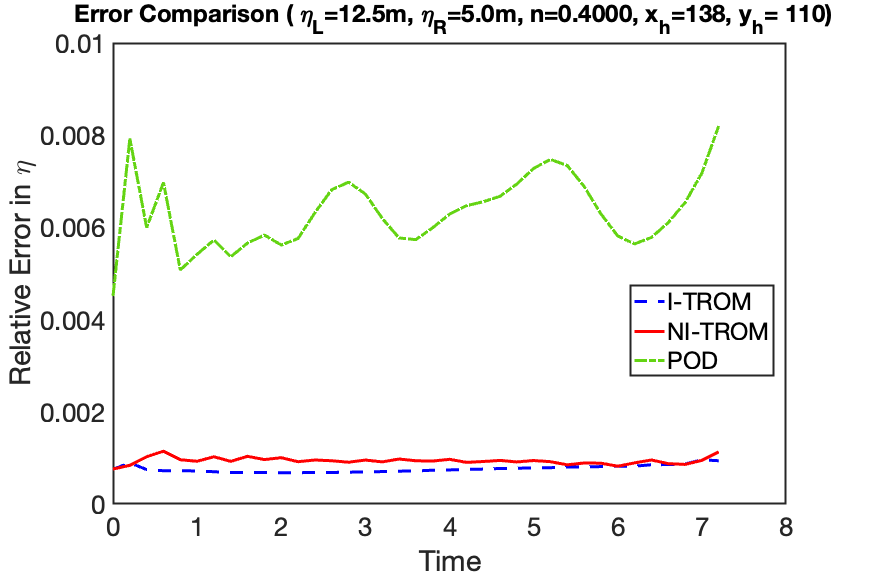}
    \subcaption{}    \label{fig:insample_error_case1}
\end{subfigure}
\hfill
\begin{subfigure}{0.48\linewidth}
    \centering    \includegraphics[width=\linewidth]{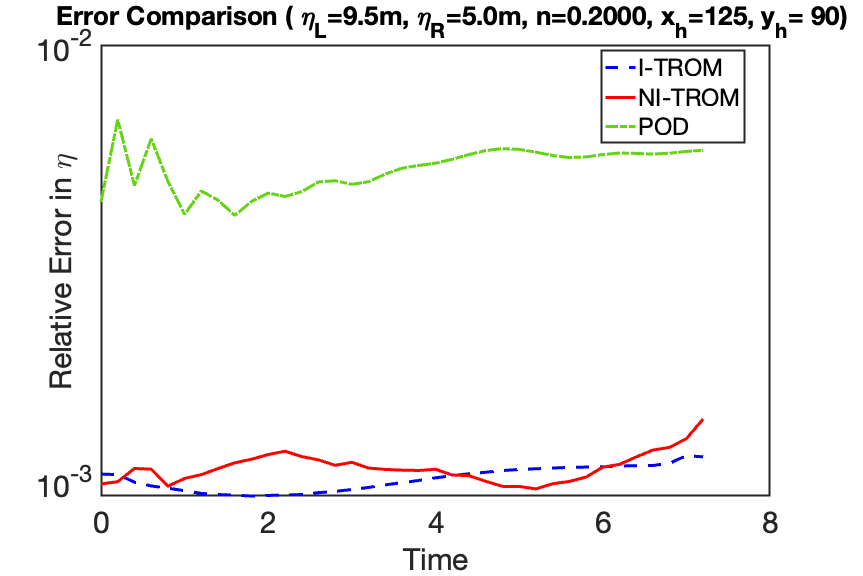}
    \subcaption{}    \label{fig:insample_error_case2}
\end{subfigure}
\vspace{0.6em}
\begin{subfigure}{0.48\linewidth}
    \centering    \includegraphics[width=\linewidth]{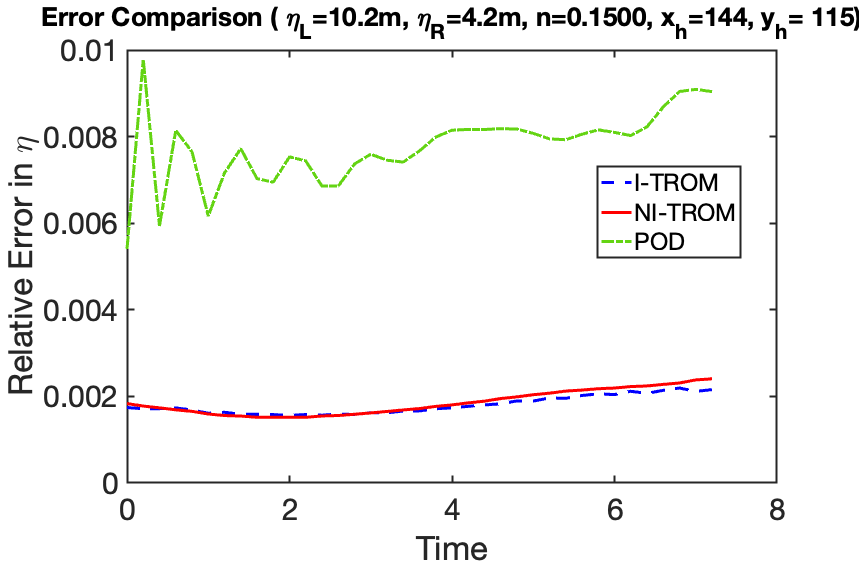}
    \subcaption{}    \label{fig:outsample_error_case3}
\end{subfigure}
\hfill
\begin{subfigure}{0.48\linewidth}
    \centering    \includegraphics[width=\linewidth]{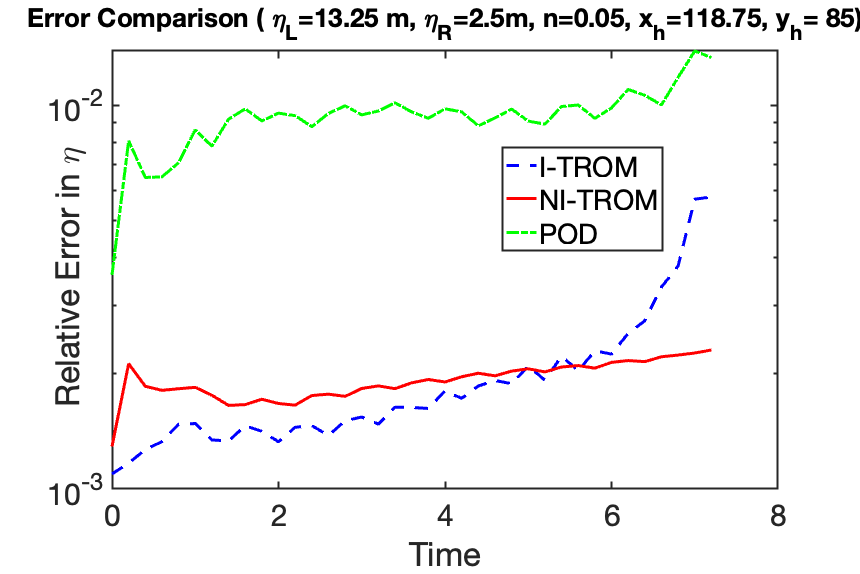}
    \subcaption{}    \label{fig:outsample_error_case4}
\end{subfigure}
\caption{Evolution of the relative error for the free-surface elevation $\eta$ 
in simulations of I-tROM, NI-tROM, and POD ROM.
Top row: (a) and (b) show the in-sample cases~1 and~2;  
bottom row: (c) and (d) show the out-of-sample cases~3 and~4.  
For each case, all reduced-order models use the same reduced bases dimensions derived from the I-tROM basis with tolerance $\epsilon_{\text{loc}} = 10^{-3}$.}
\label{fig:error_insample_outsample_2d}
\end{figure}

\begin{figure}[h]
\centering
\begin{subfigure}{0.32\linewidth}
\centering
\includegraphics[width=\linewidth]{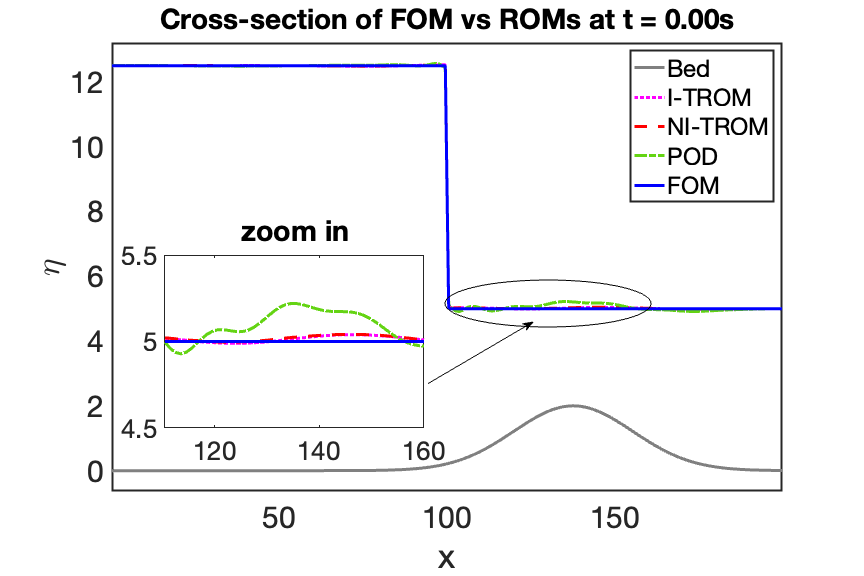}
\subcaption{}
\label{fig:insample_cs_initial_case1}
\end{subfigure}\hfill
\begin{subfigure}{0.32\linewidth}
\centering
\includegraphics[width=\linewidth]{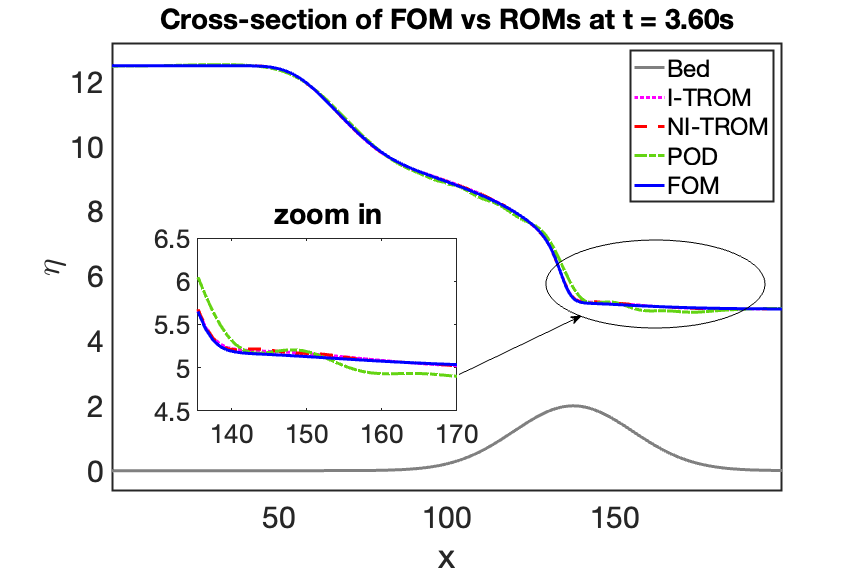}
\subcaption{}
\label{fig:insample_cs_middle_case1}
\end{subfigure}\hfill
\begin{subfigure}{0.32\linewidth}
\centering
\includegraphics[width=\linewidth]{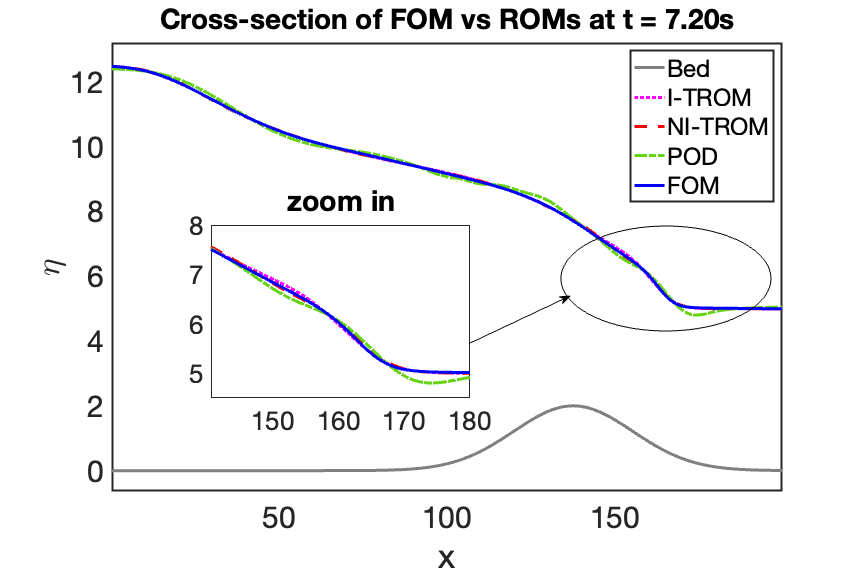}
\subcaption{}
\label{fig:insample_cs_final_case1}
\end{subfigure}
\vspace{0.6em}

\begin{subfigure}{0.32\linewidth}
\centering
\includegraphics[width=\linewidth]{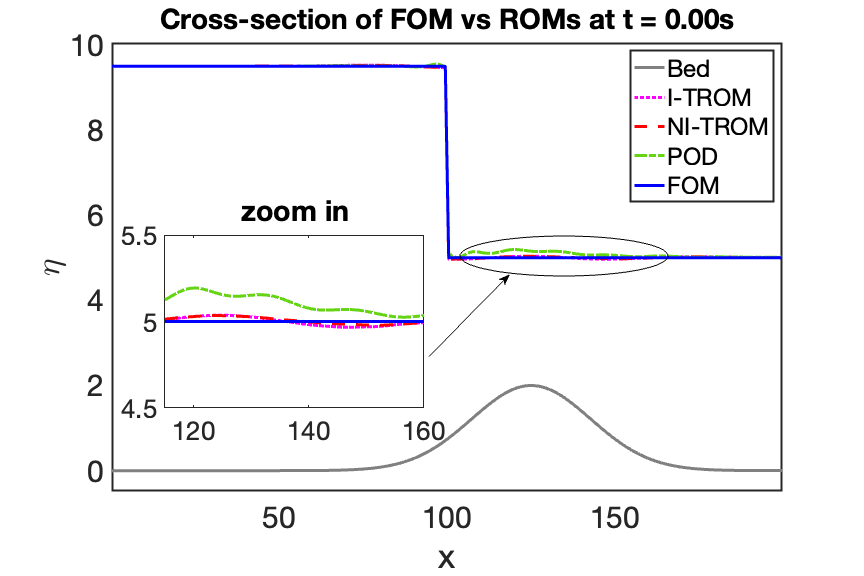}
\subcaption{}
\label{fig:insample_cs_initial_case2}
\end{subfigure}\hfill
\begin{subfigure}{0.32\linewidth}
\centering
\includegraphics[width=\linewidth]{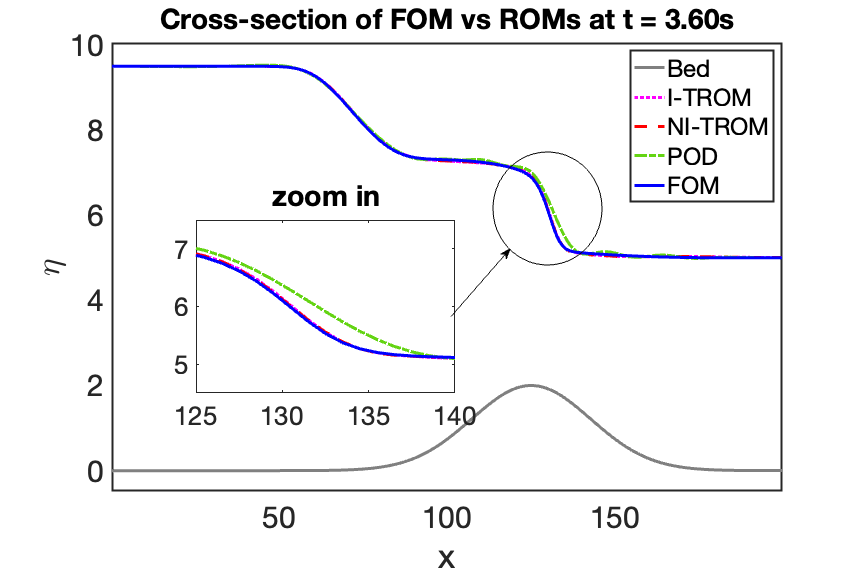}
\subcaption{}
\label{fig:insample_cs_middle_case2}
\end{subfigure}\hfill
\begin{subfigure}{0.32\linewidth}
\centering
\includegraphics[width=\linewidth]{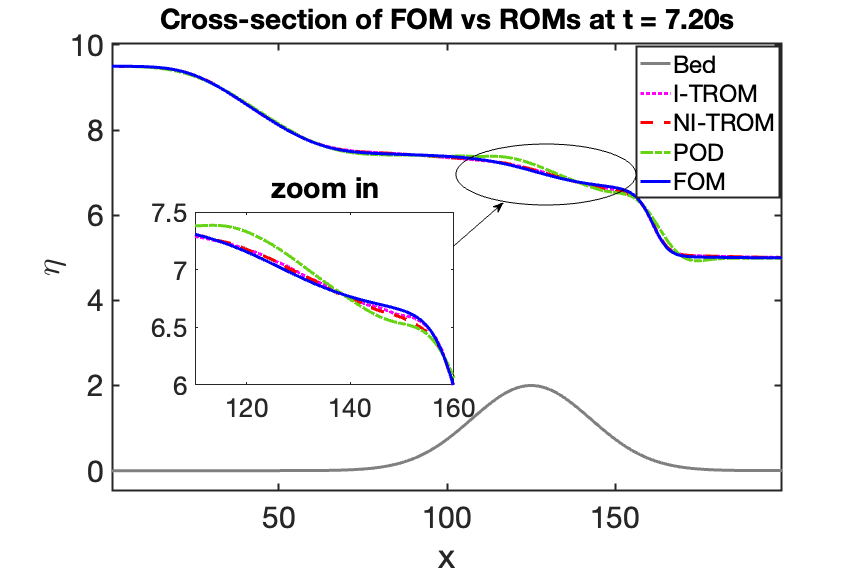}
\subcaption{}
\label{fig:insample_cs_final_case2}
\end{subfigure}
\caption{Snapshots of the free surface height cross-section $\eta(x, y=L_y/2, t)$ for $t=0, 3.6, 7.2$ in ROM simulations with the in-sample parameter values for Case 1: (a)--(c) and Case 2: (d)--(f).}
\label{fig:cross_section_insample_2d}
\end{figure}

\begin{figure}[h]
\centering

\begin{subfigure}{0.32\linewidth}
\centering
\includegraphics[width=\linewidth]{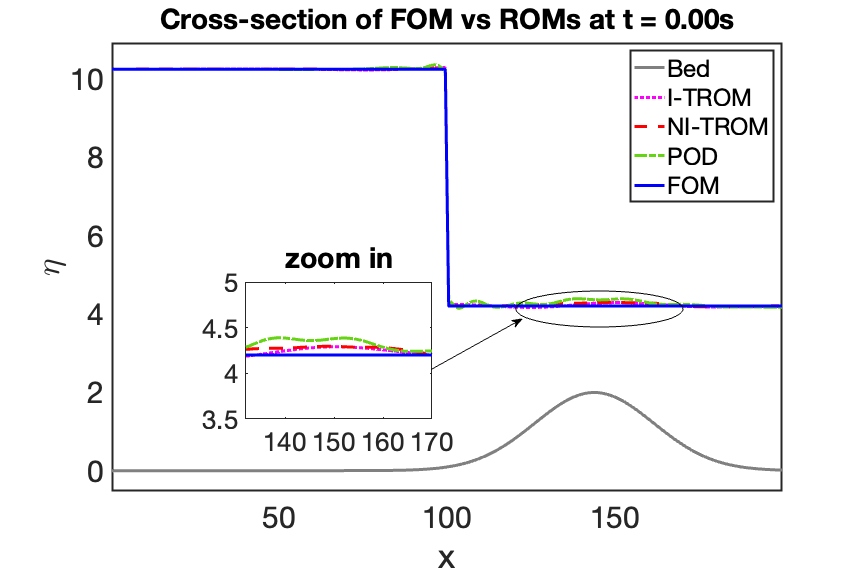}
\subcaption{}
\label{fig:outsample_cs_initial_case3}
\end{subfigure}\hfill
\begin{subfigure}{0.32\linewidth}
\centering
\includegraphics[width=\linewidth]{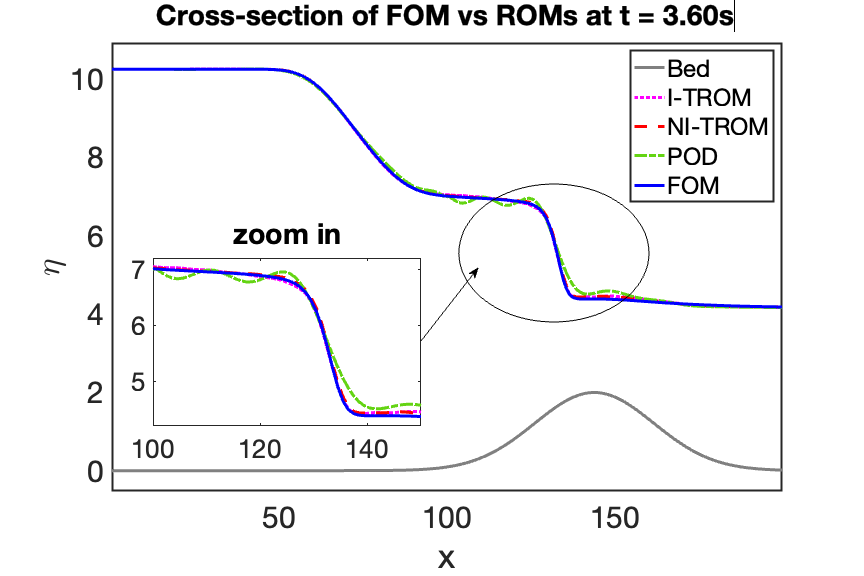}
\subcaption{}
\label{fig:outsample_cs_middle_case3}
\end{subfigure}\hfill
\begin{subfigure}{0.32\linewidth}
\centering
\includegraphics[width=\linewidth]{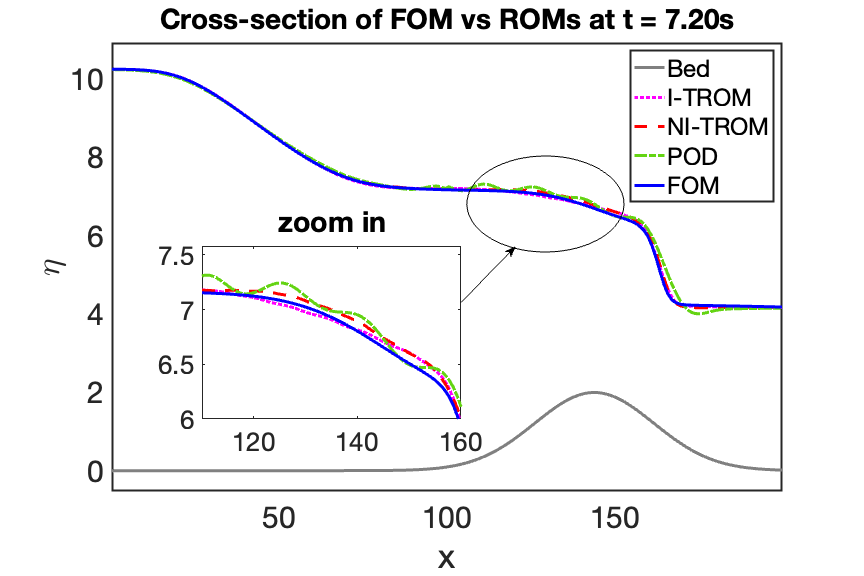}
\subcaption{}
\label{fig:outsample_cs_final_case3}
\end{subfigure}

\vspace{0.6em}

\begin{subfigure}{0.32\linewidth}
\centering
\includegraphics[width=\linewidth]{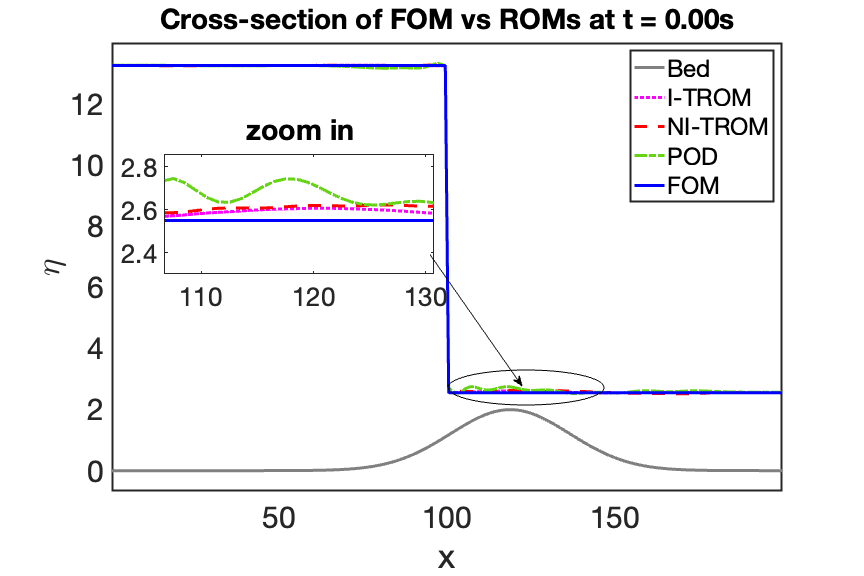}
\subcaption{}
\label{fig:outsample_cs_initial_case4}
\end{subfigure}\hfill
\begin{subfigure}{0.32\linewidth}
\centering
\includegraphics[width=\linewidth]{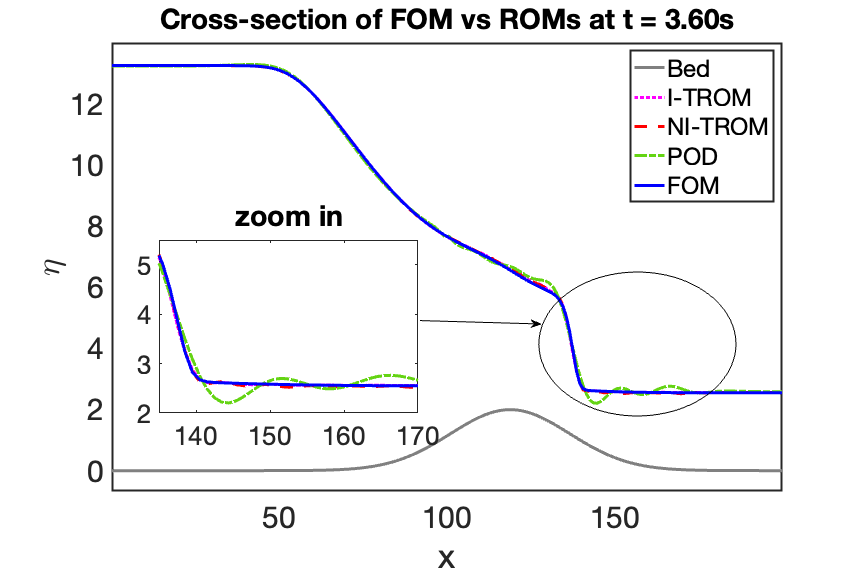}
\subcaption{}
\label{fig:outsample_cs_middle_case4}
\end{subfigure}\hfill
\begin{subfigure}{0.32\linewidth}
\centering
\includegraphics[width=\linewidth]{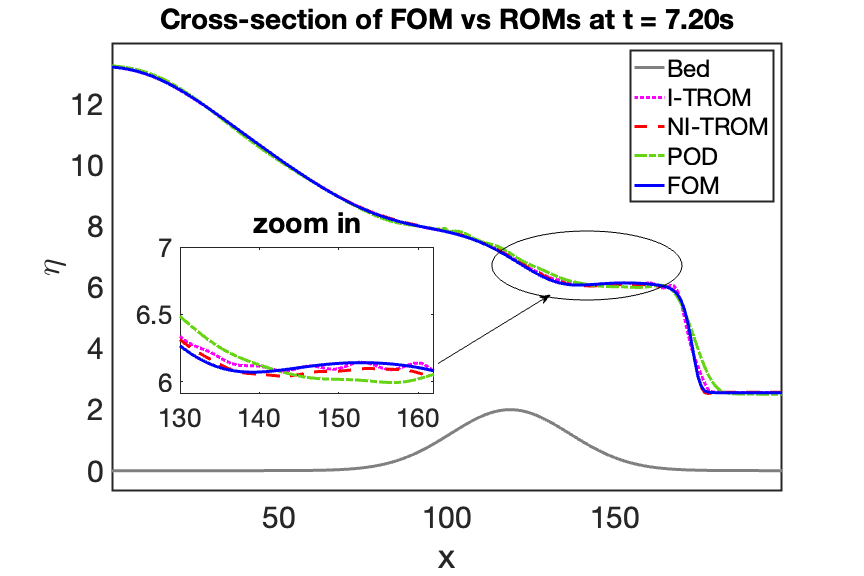}
\subcaption{}
\label{fig:outsample_cs_final_case4}
\end{subfigure}

\caption{Snapshots of the free surface height cross-section $\eta(x, y=L_y/2, t)$ for $t=0, 3.6, 7.2$ in ROM simulations with the out-of-sample parameter values for Case 3: (a)--(c) and Case 4: (d)--(f).}
\label{fig:cross_section_outsample_2d}
\end{figure}

\begin{figure}[h]
\centering
\includegraphics[width=\textwidth]{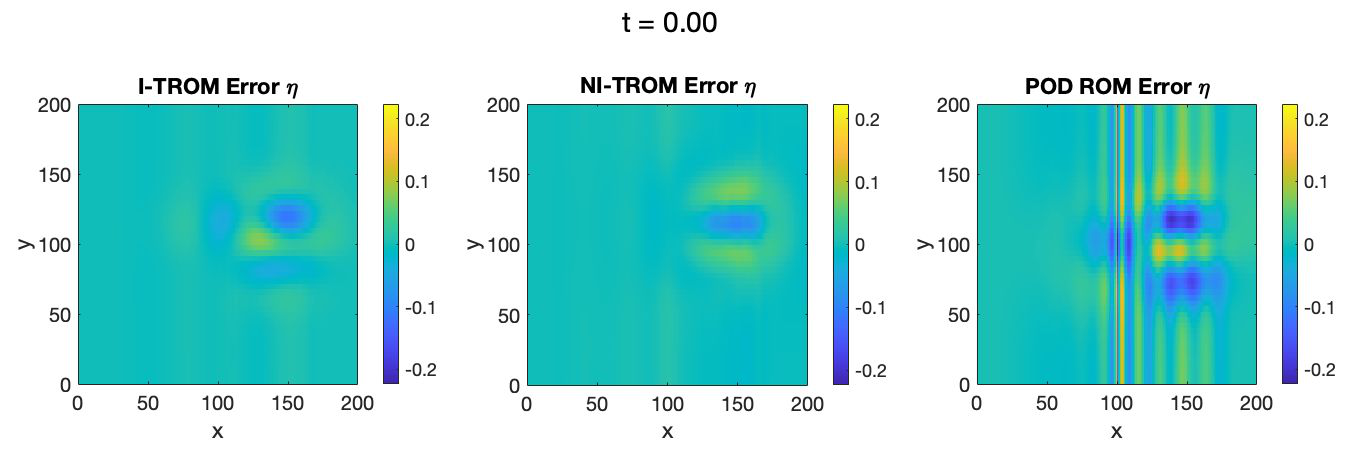}
\includegraphics[width=\textwidth]{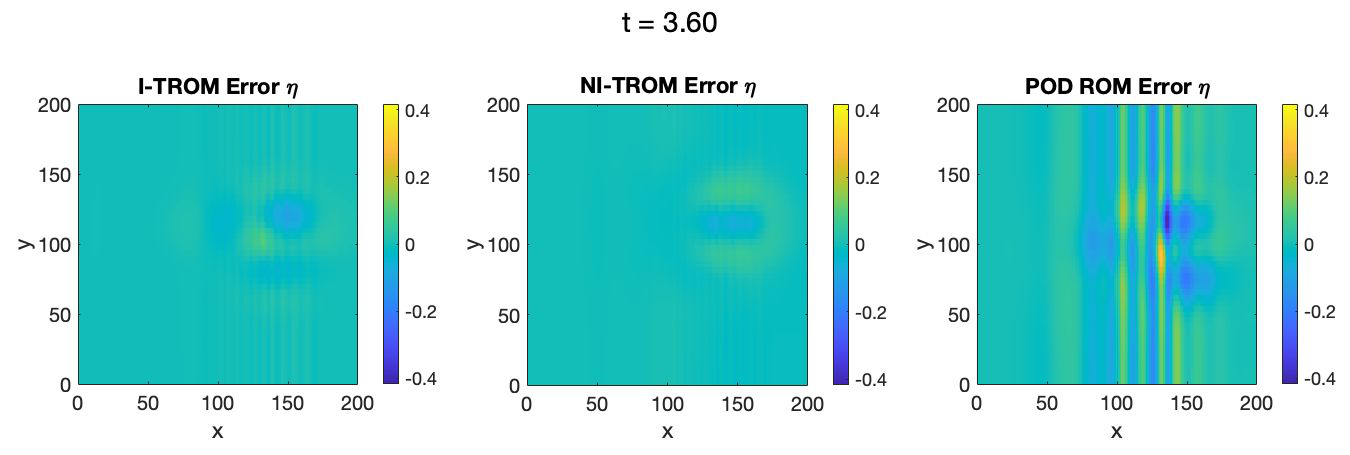}
\includegraphics[width=\textwidth]{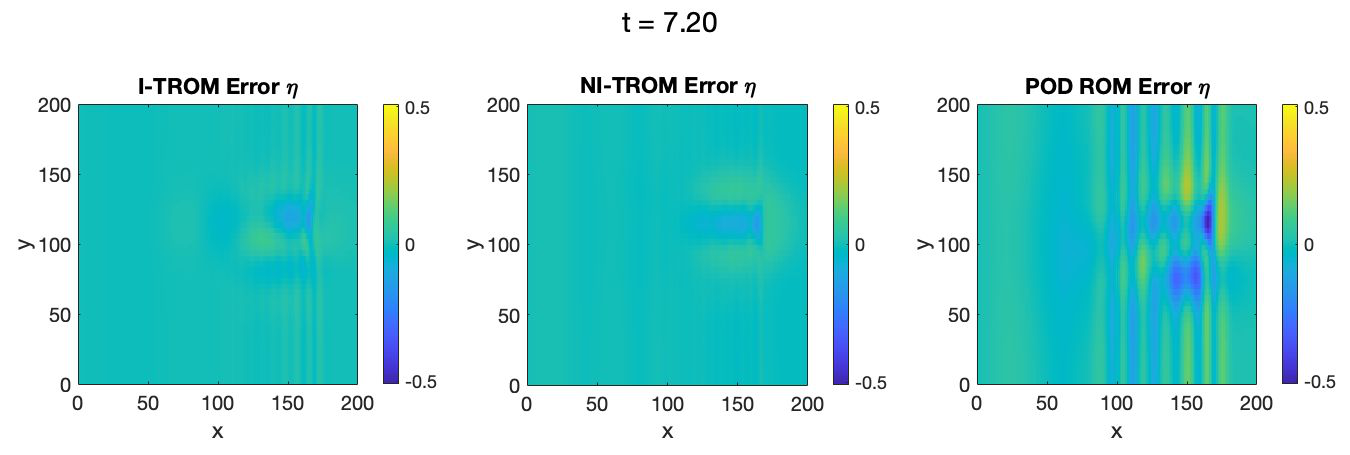}
    \caption{Spatial distribution of the error for the surface height
    $\eta^{FOM}(x,y,t) - \eta^{ROM}(x,y,t)$ for $t=0$ (top row), $t=3.6$ (middle row), and $t=7.2$ (bottom row)
    in simulations of Case 3. Left - I-tROM, Middle - NI-tROM, right - POD ROM.}
    \label{fig:eta_FOM-ROM_error_all}
\end{figure}

\begin{figure}[h]
\centering
\includegraphics[width=\textwidth]{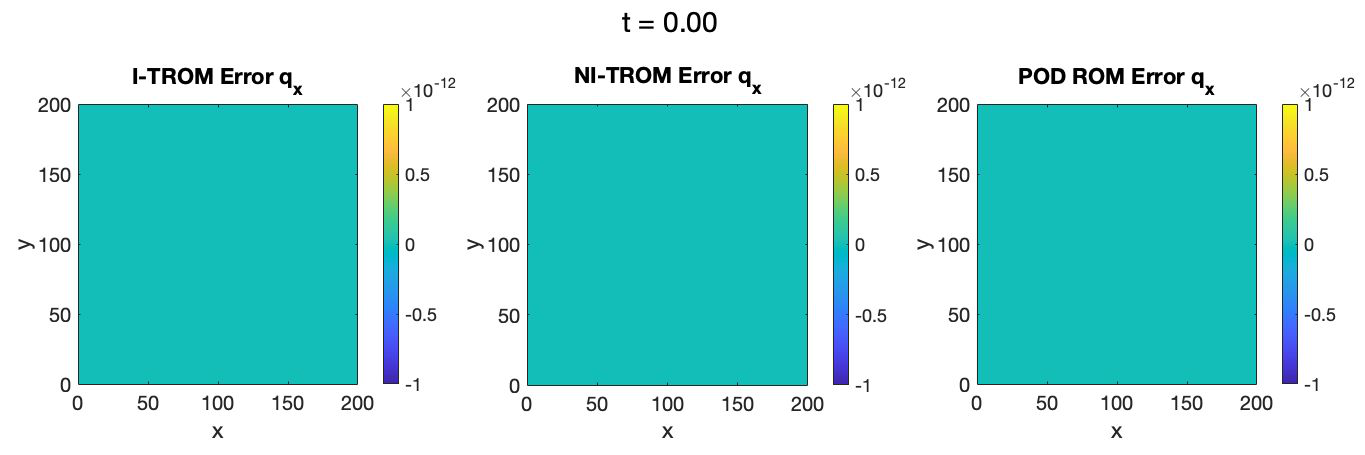}
\includegraphics[width=\textwidth]{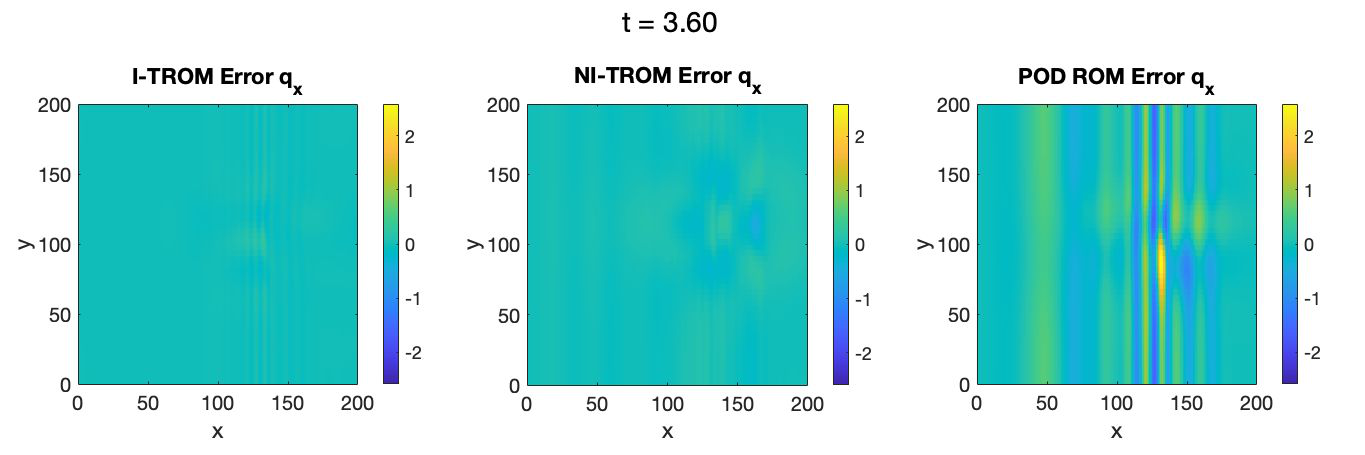}
\includegraphics[width=\textwidth]{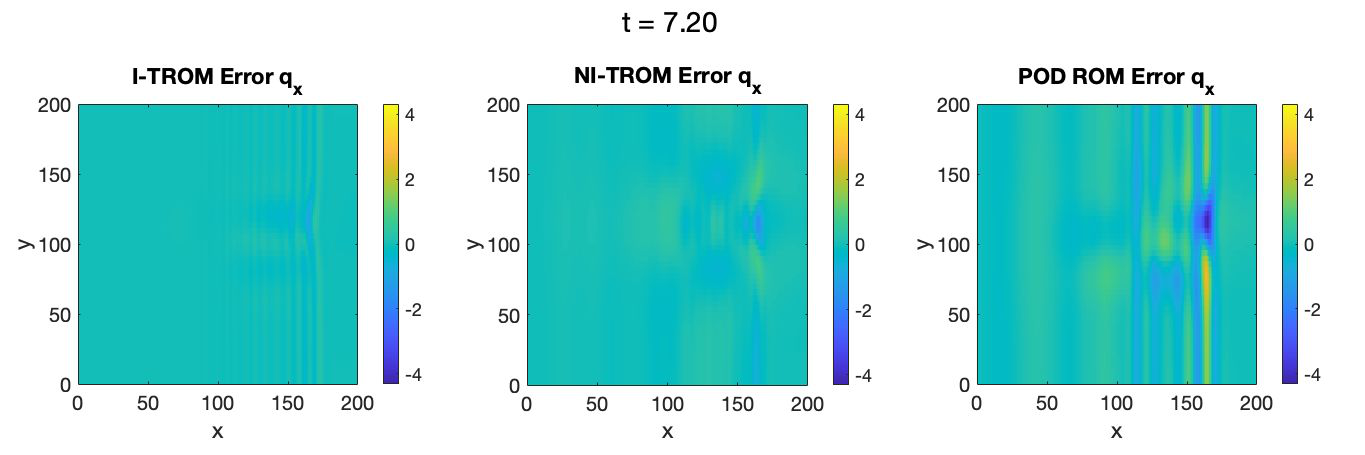}
\caption{Spatial distribution of the   error for the discharge
    $q_x^{FOM}(x,y,t) - q_x^{ROM}(x,y,t)$ for $t=0$ (top row), $t=3.6$ (middle row), and $t=7.2$ (bottom row)
    in simulations of Case 3. Left - I-tROM, Middle - NI-tROM, right - POD ROM.}
    \label{fig:qx_FOM-ROM_error_all}
\end{figure}

\rev{
In this section, we present numerical experiments for the 2D SWE model to evaluate the performance of the tensorial reduced-order models in a high-dimensional setting. The number of varied parameters is also higher compared with the 1D SWE example. 
Both interpolatory (I-TROM) and non-interpolatory (NI-TROM) tROM formulations are compared with the classical proper orthogonal decomposition (POD) ROM and the full-order model (FOM). 
Two in-sample and two out-of-sample test cases are considered to examine the generalization of the reduced-order models across different parameter configurations.

Each case is defined by the left and right free-surface elevations $(\eta_L, \eta_R)$, the Manning roughness coefficient $n$, and the hill center coordinates $(x_h, y_h)$. 
The specific parameter values considered here are:
\[
\begin{aligned}
\text{In-sample:} \quad 
&(\eta_L, \eta_R, n, x_h, y_h) = (12.5,\, 5.0,\, 0.4,\, 137.5,\, 110.0) \ \text{(Case~1)},\\
&(\eta_L, \eta_R, n, x_h, y_h) = (9.5,\, 5.0,\, 0.2,\, 125.0,\, 90.0) \ \text{(Case~2)},\\[4pt]
\text{Out-of-sample:} \quad 
&(\eta_L, \eta_R, n, x_h, y_h) = (10.25,\, 4.2,\, 0.15,\, 143.75,\, 115.0) \ \text{(Case~3)},\\
&(\eta_L, \eta_R, n, x_h, y_h) = (13.25,\, 2.50,\, 0.05,\, 118.75,\, 85.0) \ \text{(Case~4)}.
\end{aligned}
\]

For each particular case above (fixed values of parameters), 
all ROMs use identical reduced dimensions derived from the I-TROM basis with tolerance $\epsilon_{\text{loc}} = 4*10^{-3}$. However, these dimensions might change from one case to another. In particular, the reduced dimensions for $(h, q_x, q_y)$ are
case 1: (22, 26, 15),
case 2: (23, 27, 15),
case 3: (20, 25, 14), and
case 4: (24, 29, 17).
%

The bottom topography plays an important role in 2D simulations. Reduced bases are sensitive to the position of the hill, $(x_h, y_h)$, and if the reduced basis does not represent the bottom topography adequately, it can lead to larger errors even at time $t=0$.
We define the $L^2$ relative error and and MAX error as
\[
L^2_{rel}(t) = \frac{\left(\iint\limits_{\Omega} \left(h^{FOM}(x,y,t) - h^{ROM}(x,y,t)\right)^2 \, dxdy\right)^{1/2}}{\left(\iint\limits_{\Omega} \left(h^{FOM}(x,y,t) \right)^2 \, dxdy\right)^{1/2}}
\]
and $E_{max}(t) = \max_{x,y} |h^{FOM}(x,y,t) - h^{ROM}(x,y,t)|$, and similarly for $q_x$ and $q_y$. Errors are computed in discrete sense.
We consider numerical projection errors at time $t=0$ in Table \ref{tab:error2D}.
For both in-sample and out-of-sample parameter values, initial projection errors for POD ROMs are about 3-4 times larger than those of the I-tROMs and NI-tROMs.
This degrades the accuracy of POD ROMs even for earlier times in all simulations.
Figure \ref{fig:error_insample_outsample_2d} depicts the temporal evolution of relative $L^2$ errors for four test cases. Relative errors for POD ROMs are approximately 3-4 times larger for both in-sample and out-of-sample parameter values for the whole duration of simulations $t\in[0,7.2]$.
These larger errors are manifested strongly from the very beginning of simulations ($t=0$) are are due to the inability of POD ROM to adequately project the initial conditions for the free water surface.
Numerical errors increase for I-tROM towards the end of simulations for case 4.
This indicates that the NI-tROM is slightly more robust in simulations with out-of-sample parameter values.

We also show snapshots of the free surface height cross-section $\eta(x, y=L_y/2, t)$ for $t=0, 3.6, 7.2$ in Figures
\ref{fig:cross_section_insample_2d} 
and
\ref{fig:cross_section_outsample_2d}
for the in-sample and out-of-sample parameter values, respectively.
Larger numerical errors for POD ROM are visible for all times for both in-sample and out-of-sample simulations.

We depict the spatial distribution of absolute errors for $\eta(x,y,t)$ and $q_x(x,y,t)$ for Case 3 in Figures 
\ref{fig:eta_FOM-ROM_error_all} 
and
\ref{fig:qx_FOM-ROM_error_all},
respectively. Numerical errors for $q_y$ follow a similar trend. 
Figures 
\ref{fig:eta_FOM-ROM_error_all}
and 
\ref{fig:qx_FOM-ROM_error_all}
Demonstrate that numerical errors for POD ROM are larger compared to both I-tROM and NI-tROM. Moreover, numerical errors have a wave-like structure in space
and affect the whole computational domain downstream and upstream of the dam.
}

\renewcommand{\thetable}{8} 
\begin{table}[h]
\centering
\renewcommand{\arraystretch}{1.2}
\setlength{\tabcolsep}{6pt}
\begin{tabular}{llcc}
\hline
\text{Case Type} & \text{ROM Type} & \text{Relative $L_2$ Error} & \text{Max Error} \\
\hline
\multicolumn{4}{l}{\text{In-sample Case 1:} $(\eta_L,\eta_R,n,x_h,y_h) = (12.5,\,5.0,\,0.4,\,137.5,\,110.0)$} \\ 
\hline
 & Interpolatory tROM     & $7.48\times10^{-4}$ & $5.26\times10^{-2}$ \\
 & Non-Interpolatory tROM & $7.54\times10^{-4}$ & $5.34\times10^{-2}$ \\
 & POD ROM               & $4.512\times10^{-3}$ & $2.2\times10^{-1}$ \\[3pt]
\hline
\multicolumn{4}{l}{\text{In-sample Case 2:} $(\eta_L,\eta_R,n,x_h,y_h) = (9.5,\,5.0,\,0.2,\,125.0,\,90.0)$} \\ 
\hline
 & Interpolatory tROM     & $1.04\times10^{-3}$ & $5.206\times10^{-2}$ \\
 & Non-Interpolatory tROM & $9.8\times10^{-4}$ & $4.93\times10^{-2}$ \\
 & POD ROM               & $4.34\times10^{-3}$ & $1.93\times10^{-1}$ \\[3pt]
\hline
\multicolumn{4}{l}{\text{Out-of-sample Case 3:} $(\eta_L,\eta_R,n,x_h,y_h) = (10.25,\,4.2,\,0.15,\,143.75,\,115.0)$} \\ 
\hline
 & Interpolatory tROM     & $1.73\times10^{-3}$ & $1.1\times10^{-1}$ \\
 & Non-Interpolatory tROM & $1.82\times10^{-3}$ & $9.57\times10^{-2}$ \\
 & POD ROM               & $5.40\times10^{-3}$ & $2.59\times10^{-1}$ \\[3pt]
\hline
\multicolumn{4}{l}{\text{Out-of-sample Case 4:} $(\eta_L,\eta_R,n,x_h,y_h) = (13.25,\,2.50,\,0.05,\,118.75,\,85.0)$} \\ 
\hline
 & Interpolatory tROM     & $1.09\times10^{-3}$ & $8.6\times10^{-2}$ \\
 & Non-Interpolatory tROM & $1.28\times10^{-3}$ & $7.24\times10^{-2}$ \\
 & POD ROM               & $3.7\times10^{-3}$ & $1.95\times10^{-1}$ \\
\hline
\end{tabular}
\caption{Comparison of relative $L_2$, and maximum errors at $t=0$ (initial condition) for in-sample and out-of-sample cases for different ROM types.}
\label{tab:error2D}
\end{table}

%

\begin{figure}[h]
\centering
\includegraphics[width=\textwidth]{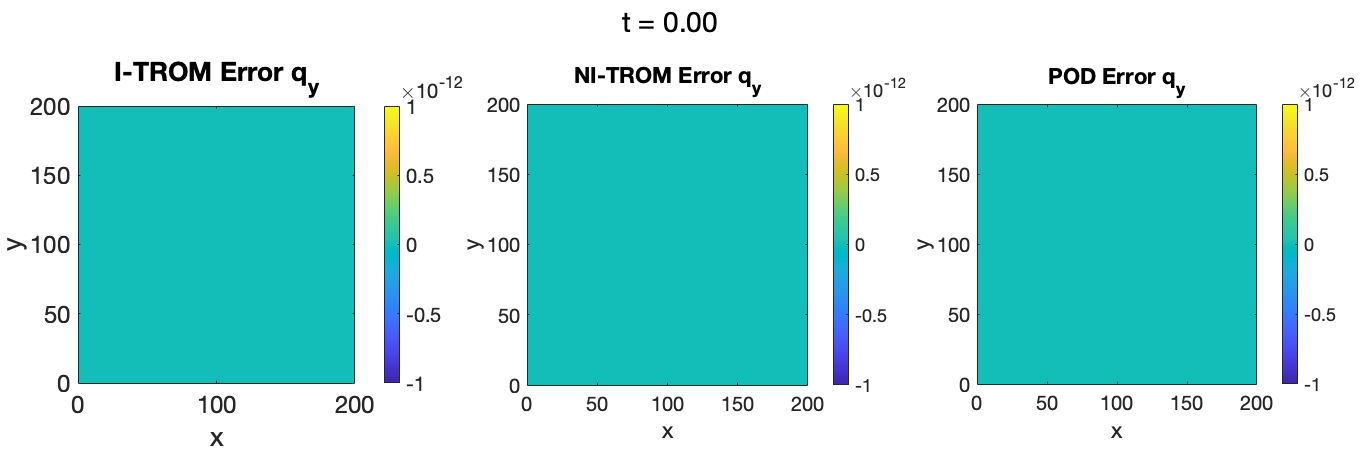}
\includegraphics[width=\textwidth]{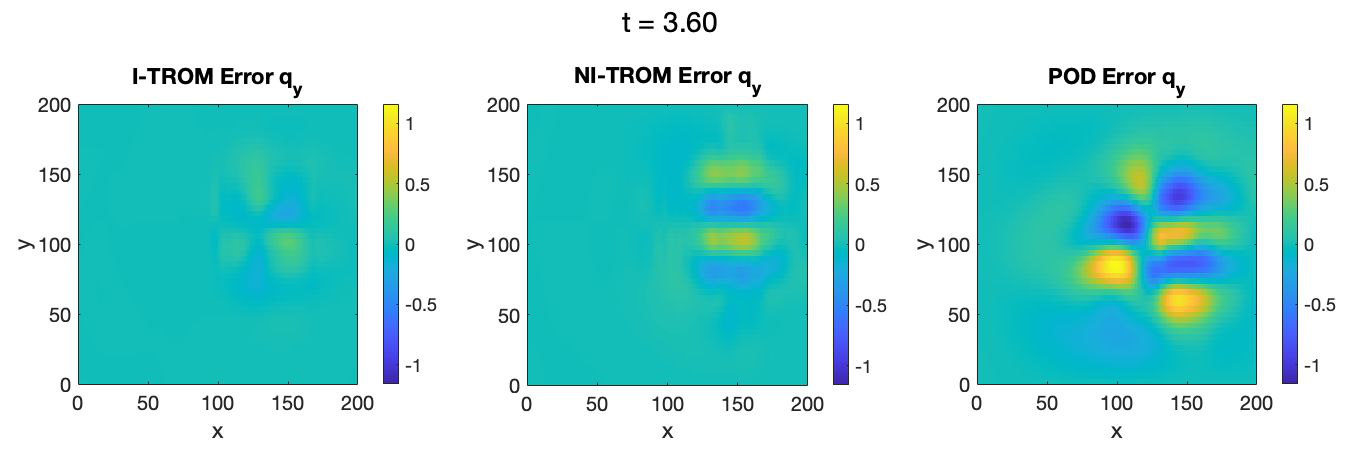}
\includegraphics[width=\textwidth]{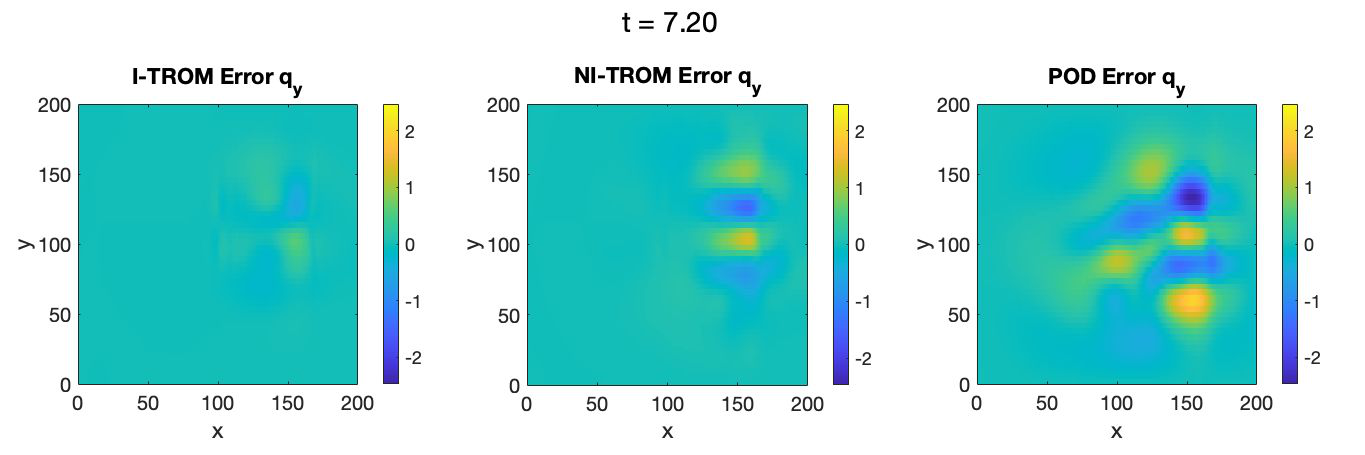}
\caption{Spatial distribution of the   error for the discharge
    $q_y^{FOM}(x,y,t) - q_y^{ROM}(x,y,t)$ for $t=0$ (top row), $t=3.6$ (middle row), and $t=7.2$ (bottom row)
    in simulations of Case 3. Left - I-tROM, Middle - NI-tROM, right - POD ROM.}
    \label{fig:qy_FOM-ROM_error_all}
\end{figure}

\FloatBarrier

\section{Conclusions}
\label{sec:conc}

In this paper, we developed a non-interpolatory variant of the tROM for the parameterized shallow-water dam-break problem. We also investigated the dependence of solutions on the parameters \((\hleft, \hright)\).
This setup presents several challenges specific to hyperbolic problems. First, the dimension of the reduced basis is expected to be relatively high compared to parabolic problems. This is directly linked to the slow decay of the Kolmogorov \(N\)-width in hyperbolic settings. Second, the solution may develop shock waves. Combined with parametric dependence, this necessitates the development of advanced numerical methods for model reduction. We demonstrated that a properly designed tROM can successfully address this challenge and provide accurate approximations of high-fidelity simulations.

Another difficulty is particular to the problem at hand and stems from the behavior of the SWE solutions for small values of \(\hright\). As shown in the appendix, the derivative of the solution with respect to the second parameter is not uniformly bounded in the neighborhood of \(\hright = 0\). Consequently, constructing a parametric ROM for the wet-bed problem is especially challenging. The non-interpolatory tROM with adaptive sampling --- using Chebyshev nodes to cluster samples near \(\hright = 0\) --- effectively addresses this issue and offers an efficient framework for projection-based ROM development. As demonstrated in Section~\ref{sec:312}, tROMs with both adaptive and fixed local dimensions can successfully handle the limiting case \(\hright \to 0\).

Our results indicate that the novel tROM methodology can successfully handle hyperbolic problems with varying parameters. In particular, tROMs developed in this paper perform very well for a test problem in the regime representing typical behavior of hyperbolic systems. Therefore, we expect that the tROM approach discussed in this paper will also be applicable \rev{to other hyperbolic systems.}

\rev{
Application of tensor ROMs to hyperbolic systems still has room for improvement.
In particular, the Principal Interval Decomposition method has been successful in improving POD-ROMs by employing a more localized-in-time construction of the reduced basis~\cite{borggaard2007interval}. This approach has been effectively applied to the shallow water equations in~\cite{zokagoa2018pod,solan2023development,solan2025combination,gomez2025well} within the POD–ROM framework.
The same localization ideas can be extended to the tROM framework, yielding reduced bases that are localized both in time and in the parameter domain.
}

However, one of the major challenges in the ROM community is developing accurate and efficient projection-based ROMs for problems that involve a relatively high-dimensional vector of parameters. As discussed in our paper, for hyperbolic problems, it might be necessary to perform clustered sampling near certain parameter values; this would present a major challenge if the vector of parameters is high-dimensional. Thus, further development of tensor completion and interpolation techniques might be necessary for problems that require clustered sampling.
This will be addressed in elsewhere.

\section*{Acknowledgment}

The authors R.B.M. and M.O. were supported in part by the U.S. National Science Foundation under awards
DMS-2309197 and DMS-2408978.

\appendix

\section{Analytical Solution of the Dam-Break Problem}
\label{sec:ap}
In this section, we recall analytical solutions for the 1D dam-break problem and study the behavior of the wet-bed case solution in the proximity of the dry-bed case. 
\rev{For an infinite domain,} the dam-break problem can be solved analytically for both cases, using a combination of Riemann invariants, \rev{Rankine--Hugoniot} conditions, and Newton iteration to determine the middle state.

\subsubsection*{Dry-Bed Case (\(\hright = 0\))}
For the dry-bed case, where the initial \rev{downstream} water height \(\hright = 0\), the solution consists of a rarefaction wave spreading from the dam location \(\xdam\). Given an initial \rev{upstream} water height \(\hleft\), the solution at time \(t\) and position \(x\) is determined by the characteristic speeds. Let \(\cleft = \sqrt{g \hleft}\) be the gravity wave speed on the left side, and \(x_\ell = \xdam - \cleft t\) and \(x_r = \xdam + 2 \cleft t\) define the left and the right ends (i.e. the tail and the head) of the rarefaction wave, respectively. The analytical solution for water depth \(h(x, t)\) and velocity \(u(x, t)\) is given by
\begin{equation}
\label{eq:1086}
\left[ u(x,t), \, h(x, t)  \right] = \begin{cases} 
\left[0, \, \hleft \right] & \text{if } x < x_\ell, \\
\left[
\dfrac{2}{3} \left( \cleft + \delta \right) , \,
\dfrac{4}{9g} \left( \cleft - \dfrac{\delta}{2} \right)^2
\right]
& \text{if } x_\ell \leq x \leq x_r, \\
\left[0, \, 0\right] & \text{if } x > x_r,
\end{cases}
\end{equation}
where \(\delta = (x - \xdam)/t\).
This solution reflects a rarefaction wave expanding into the dry region, with the water depth decreasing smoothly from \(\hleft\) to 0 across the wave.

\subsubsection*{Wet-Bed Case (\(\hright > 0\))}
For the wet-bed case, where \(\hright > 0\), the solution involves a rarefaction wave on the left, a constant middle state \((h_m, u_m)\), and a shock wave on the right. The middle state and shock speed \(s\) are determined by solving a system of equations using the Riemann invariant across the left rarefaction wave and the \rev{Rankine--Hugoniot} conditions.
The middle state \((h_m, u_m)\) and shock speed \(s\) satisfy the nonlinear equation
\begin{equation}
\label{eq:umhm}
u_m + 2\sqrt{g h_m} - 2\sqrt{g \hleft} = 0,
\end{equation}
where
\begin{equation}
\label{eq:umhm2}
u_m = s - \frac{g \hright}{4s} \left( \sqrt{1 + \frac{8 s^2}{g \hright}} + 1 \right),
\qquad
h_m = \frac{\hright}{2} \left( \sqrt{1 + \frac{8 s^2}{g \hright}} - 1 \right).
\end{equation}
The $s$ solving the above equation can be computed numerically using, for instance, Newton's method.
Newton's iteration starts with an initial guess \(s = \sqrt{g \hleft}\) and iteratively updates \(s\) until convergence.

Given the middle state, the solution at time \(t\) and position \(x\) is determined by the positions of the rarefaction wave and shock. The gravity wave speed  $\cleft$, $\delta$, and the tail of the wave $x_1=x_\ell$ are the same as in the dry case. Let \(\cm = \sqrt{g h_m}\), \(x_2 = \xdam + (u_m - \cm)t\), and \(x_3 = \xdam + s t\). The analytical solution for the wet case is then given by
\begin{equation}\label{eq:1134}
\left[ u(x,t), \, h(x, t)  \right] =
\begin{cases} 
\left[0, \, \hleft \right] & \text{if } x \leq x_1, \\
\left[ \dfrac{2}{3} \left( \cleft + \delta \right), \, 
\dfrac{4}{9g} \left( \cleft - \dfrac{\delta}{2} \right)^2 \right] & \text{if } x_1 < x < x_2, \\
\left[u_m, \, h_m \right] & \text{if } x_2 \leq x \leq x_3, \\
\left[ 0, \, \hright \right] & \text{if } x > x_3.
\end{cases}
\end{equation}
This solution captures the rarefaction wave, the constant middle state, and the shock wave propagating to the right.

\subsubsection*{Convergence and regularity for \(\hright \to 0\)}
We are now interested in the behavior of the wet-bed solution when $\hright \to 0$.
It is easy to see from the second equation in \eqref{eq:umhm2}
that $h_m \sim \sqrt{\hright} \to 0$. Then from \eqref{eq:umhm} we obtain $u_m \to 2\sqrt{g\hleft} = 2c$ also with the speed $O(\sqrt{\hright})$. On the other hand, from the first equation in \eqref{eq:umhm2}, $u_m \to s$. Therefore, $(h_m, u_m, s) \to (0, 2c, 2c)$
as $\hright \to 0$. We also have $x_3 \to x_r$ and $x_2 \to x_r$ as $\hright \to 0$.
Therefore, the solution of the wet-bed converges in $L^\infty(0,L)$ to the solution of the dry-bed as $\hR\to0$ uniformly in time 
$t \in [0, T]$ if $\xdam + 2cT < L$. However, the solution of the wet-bed is not differentiable for any $t>0$ because the shock appears immediately, since $x_3 - x_2 > 0$ for any $t>0$. This is evident since $s - u_m > 0$ for any $\hright>0$ (the first equation in \eqref{eq:umhm2}) and $c_m$ is also strictly positive for $\hright>0$.
Thus, although there is convergence in $L^\infty(0,L)$, the  wet-bed solutions  fail to converge to the  dry-bed solution in $H^1(0,L)$ norm for $\hright\to0$ and any fixed $t>0$.  

We are also interested in the regularity of the wet-bed case solution with respect the $\hright$ parameter in the proximity of $\hright=0$.
First note that $s\to c> 0$ as $\hright\to 0$ and so the propagation speed of the shock is uniformly bounded away from zero. For $x\in(x_2,x_3)$ the water depth is given by $h=h_m$ and from \eqref{eq:umhm2} we obtain $\partial h / \partial \hright \sim O(\hright^{-1/2})$ on $(x_2,x_3)$ for small $\hright$.
The interval length is   $x_3 - x_2 = (s - u_m + c_m) t \sim O(\hright^{1/4}) t$, since $s-u_m \ge 0$, $c_m = \sqrt{g h_m}$ and 
$h_m \sim O(\hright^{1/2})$. 
We conclude that  the derivative of $\partial h(x,t) / \partial \hright$ 
is unbounded on $(0,L)$ when $\hright\to 0$ and for any $t>0$. More generally, straightforward calculations show that the $L^p(0,L)$ norm of $\partial h / \partial \hright$ has $O(\hright^{\frac{1}{4p}-\frac12})$ asymptotic for $\hright\to0$.  Furthermore, similar calculations shows the same asymptotics for   $\partial u(x,t) / \partial \hright$ and for $\partial q(x,t) / \partial \hright$ in the neighborhood of the dry-bed case.

\section{Sampled Parameter Values}
\label{sec:ap2}
For the 1D simulations, we sample parameters $h_L$ and $h_R$ according to \eqref{eq:hL_nodes}
and \eqref{eq:hR_nodes}, respectively.
In all the Figures for the 1D case, we construct tROM using $N_L=13$ and the following values of 
\[
h_L =\{10.0, 11.5, 13.0, 14.5, 16.0, 17.5, 19.0, 20.5, 22.0, 23.5, 25.0, 26.5, 28.0\}.
\]
In Figures \ref{fig:exp1_performance_metrics},
\ref{fig:exp7_h_comparison},
\ref{fig:exp7_error_evolution},
\ref{fig:exp7_hu_comparison},
and
\ref{fig:exp7_error_hu_evolution},
tROM is constructed using Chebyshev nodes with $N_R=17$, resulting in the following sampled values (approximately)
\[
\hR = \{ 0, 0.038, 0.154, 0.345, 0.61, 0.945, 1.35, 1.82, 2.34, 2.93, 3.56, 4.23, 4.94, 5.68, 6.44, 7.22, 8\}.
\]
In Figures 
\ref{fig:exp3_h_profiles_errors} and
\ref{fig:exp4_h_evolution}
we consider tROM constructed with Chebyshev nodes with $\NR=9$, i.e., 
\[
\hR =
\{0, 0.154, 0.609, 1.348, 2.343, 3.555, 4.939, 6.439, 8\}.
\]


\end{document}